\documentclass[3p]{elsarticle}

\usepackage{amsfonts}
\usepackage{multirow}
\usepackage[dvipsnames]{xcolor}
\usepackage{tikz}
\usepackage{mathtools}
\usepackage{stackengine}
\usepackage{arydshln}
\newtheorem{thm}{Theorem}
\newtheorem{lem}{Lemma} 
\newtheorem{prop}{Proposition}
\newtheorem{problem}{Problem}
\newdefinition{eg}{Example}
\newproof{pf}{Proof}
\newcommand\undermat[2]{%
  \makebox[0pt][l]{$\smash{\underbrace{\phantom{%
    \begin{matrix}#2\end{matrix}}}_{\text{$#1$}}}$}#2}
\begin{document}

\begin{frontmatter}

\title{$2$-limited broadcast domination on grid graphs}

\author[uvicmath]{Aaron Slobodin\corref{cor1}\fnref{footnote1}}
\ead{aslobodin@uvic.ca}

\author[uvicmath]{Gary MacGillivray\fnref{footnote2}}
\ead{gmacgill@uvic.ca}

\author[uviccomp]{Wendy Myrvold\fnref{footnote2}}
\ead{wendym@cs.uvic.ca}

\affiliation[uvicmath]{
organization={University of Victoria},
department={Mathematics and Statistics},
addressline={P.O. BOX 1700 STN CSC},
city={Victoria},
province={B.C.},
postcode={V8W 2Y2},
country={Canada}}

\affiliation[uviccomp]{
organization={University of Victoria},
department={Computer Science},
addressline={P.O. BOX 1700 STN CSC},
city={Victoria},
province={B.C.},
postcode={V8W 2Y2},
country={Canada}}

\cortext[cor1]{Corresponding author}
\fntext[footnote1]{Supported by CGS-D funding from Natural Sciences and Engineering Research of Council of Canada (NSERC).}
\fntext[footnote2]{Supported by funding from NSERC.}

\begin{abstract}
We establish upper and lower bounds for the $2$-limited broadcast domination number of various grid graphs, in particular the Cartesian product of two paths, a path and a cycle, and two cycles. The upper bounds are derived by explicit constructions. The lower bounds are obtained via linear programming duality by finding lower bounds for the fractional $2$-limited multipacking numbers of these graphs. 
\end{abstract}

\begin{keyword}
Graph Theory \sep Domination \sep Broadcast Domination \sep Limited Broadcast Domination \sep Multipacking \sep Limited Multipacking
\end{keyword}

\end{frontmatter}

\section{Introduction}\label{intro}
Consider a city partitioned into neighbourhoods, each of which has a communication tower capable of transmitting a radio broadcast of some non-negative strength. Given this premise, a natural question is: \textit{how can we design a broadcast scheme such that broadcasts are heard by each neighbourhood while reducing the number of communication towers transmitting?} The discrete version of this question leads to  broadcast domination, which was introduced by Erwin in 2001 \cite{ogerwin} (also see \cite{erwin}).

In a $k$-limited broadcast, each vertex $v$ of a graph $G$ is assigned a broadcast strength $f(v) \in \left\{ 0,1,\dots, k \right\} $, where $k\leq$ rad$(G)$ and rad$(G)$ is the radius of $G$. We say that a vertex $u$ \textit{hears} the broadcast from $v$ if $d(u,v) \leq f(v)$, where $d(x,v)$ is the distance between $x$ and $v$ in $G$. A broadcast $f$ is \textit{dominating} if each vertex of $G$ hears the broadcast from some vertex. The \textit{cost} of a broadcast $f$ is $\sum_{v \in V(G)} f(v)$. The \textit{$k$-limited broadcast domination number} $\gamma_{b,k}(G)$ of $G$ is the minimum cost of a $k$-limited dominating broadcast. The $k$-limited broadcast domination number can also be formulated as an integer linear program (ILP). Let $G$ be a graph and fix $1 \leq k \leq rad(G)$. For each vertex $i \in V(G)$ and $\ell \in \left\{ 1,2,\dots, k \right\} $, let $x_{i,\ell}=1$ if vertex $i$ is broadcasting at strength $\ell$ and $0$ otherwise. The $k$-limited broadcast domination number  of a graph $G$ is the cost of an optimal solution to \ref{ILP}.

\renewcommand{\theequation}{ILP 1}
\begin{equation}\label{ILP}
\begin{aligned}
\textnormal{Minimize:} \quad &\sum_{\ell=1}^{k}\sum_{i \in V(G)} \ell\cdot x_{i,\ell}\\
\textnormal{Subject to:} \quad &(1) \quad\sum_{\ell=1}^{k}\hspace{0.1cm} \sum_{\text{$i \in V(G) \textnormal{ s.t. }$} \text{$d(i,j) \leq \ell$}} x_{i,\ell} \geq 1, \quad \textnormal{for each vertex } j \in V(G),  \\
                                   &(2) \quad  x_{i,\ell} \in \left\{ 0,1 \right\}  \quad \textnormal{for each vertex } i \in V(G) \textnormal{ and } \ell \in \left\{ 1,2,\dots,k \right\}  . \\
\end{aligned}
\end{equation}
\renewcommand{\theequation}{\arabic{equation}}

The $k$-limited broadcast domination number of a graph was first mentioned, although not explored, in  \cite{ogerwin}. Given a tree $T$ on $n$ vertices, the best possible upper bound $ \gamma_{b,k} (T) \leq \left\lceil n(k+2)/3(k+1)\right\rceil $ is shown in \cite{thing2}. It gives a bound for all graphs via spanning trees. However,  these bounds are likely very far from the truth for arbitrary  graphs. Specific to $2$-limited broadcast domination, Yang proved that if $G$ is a cubic graph on $n$ vertices with no 4-cycles or 6-cycles, then $\gamma_{b,2}(G) \leq n/3$ \cite{thesis}.

For each fixed positive integer $k$, the problem of deciding whether there exists a $k$-limited dominating broadcast of cost at most a given integer $B$ is NP-complete \cite{thing2}.  There are polynomial time algorithms to compute $\gamma_{b,k}(G)$ for strongly chordal graphs, interval graphs, circular arc graphs, and proper interval bigraphs \cite{thesis}. 

The $k$-limited broadcast domination problem is a restriction of the broadcast domination problem (which was also introduced in \cite{ogerwin}). The \textit{broadcast domination number} of a graph can be obtained from \ref{ILP} by setting $k=rad(G)$. The broadcast domination number of a graph with $n$ vertices can be found in $O(n^6)$ time \cite{HEGGERNES20063267}. There exists improved algorithms for trees \cite{lineartreesmaster,treesnotmasters}, interval graphs \cite{broadinterval}, and strongly chordal graphs \cite{frankmasters,frankrickandgarry}. A survey of results on broadcast domination can be found in \cite{Henning2021}.

\subsection{$k$-Limited Multipacking} 
The dual of the linear programming (LP) relaxation of \ref{ILP} leads to \textit{fractional $k$-limited multipacking}.  For each vertex $i \in V(G)$, define a variable $y_i$. The \textit{fractional $k$-limited multipacking packing number} $mp_{f,k}(G)$ of a graph $G$ is the cost of an optimal solution to \ref{LP:multipackc}.
\renewcommand{\theequation}{LP 2}
\begin{equation}\label{LP:multipackc}
\begin{aligned}
\textnormal{Maximize:} \quad &\sum_{i \in V(G)} y_i\\
\textnormal{Subject to:} \quad \hspace{-0.1cm} &(1) \quad \sum_{i \in V(G) \textnormal{ s.t. } d(i,j) \leq \ell} y_i \leq \ell, \quad \textnormal{for each vertex } j \in V(G) \textnormal{ and } \ell \in \left\{ 1,2, \dots, k \right\} , \textnormal{ and}   \\
                                               &(2) \quad  y_i\geq 0, \quad \textnormal{for each vertex } i \in V(G). \\
\end{aligned}
\end{equation}
\renewcommand{\theequation}{\arabic{equation}}

The \textit{multipacking number} $mp(G)$ of a graph $G$ is the cost of an optimal solution to \ref{LP:multipackc} when interpreted as a 0-1 ILP and setting $k=rad(G)$. This problem  was introduced in Teshima's Master's thesis \cite{lteshima} in 2012 (also see \cite{laurakeikarick}).   The complexity of deciding whether a given graph $G$ has a multipacking number at least a given integer $B$ remains an open problem. There is a $(2 +O(1)) $ approximation algorithm for this decision problem on undirected graphs \cite{rickbeaudoufocaud}, and there are polynomial time algorithms for trees  \cite{lteshima,multipackkiekaandlaura,frankrickandgarry} and strongly chordal graphs \cite{frankmasters, frankrickandgarry}. For further information, see \cite{Henning2021}.

The cost of an optimal solution of \ref{LP:multipackc}, when interpreted as a 0-1 ILP, is the \textit{$k$-limited multipacking number} of a graph.  The $k$-limited multipacking problem is a restriction of the multipacking problem. Relatively little is known about $k$-limited multipackings on graphs. For each fixed integer $k$, the problem of deciding whether a given graph $G$ has a $k$-limited multipacking number at least a given integer $B$ is NP-complete \cite{thesis}.  There is a polynomial time algorithm for  the $k$-limited multipacking number for strongly chordal graphs \cite{frankrickandgarry}.

\subsection{Outline}
After a sequence of many papers, the $1$-limited broadcast domination number  (i.e. the domination number) of the Cartesian product of two paths was finally determined (in \cite{domgrids}). This paper focuses on extending this work to determine bounds for the $2$-limited broadcast domination number of the Cartesian product of two paths, a path and a cycle, and two cycles. Section \ref{sec:upperbound} introduces the method of constructing $2$-limited dominating broadcasts by ``tiling'' the graph with ``broadcast tiles.'' Sections \ref{sec:paths}, \ref{sec:pathsandcycle}, and \ref{sec:cycles} summarize the tiling results from Chapters 2, 3 and 4 of \cite{slobodin}  which give $2$-limited dominating broadcasts on the Cartesian products of two paths, a  path and a cycle, and two cycles. Section \ref{sec:multipack} utilizes  \textit{fractional $2$-limited multipacking} to establish lower bounds for the $2$-limited broadcast domination numbers of the Cartesian product of two paths, a path and a cycle, and two cycles. Section \ref{sec:conclusion} includes suggested questions for future research.

\section{The Tiling Method for establishing Upper Bounds}\label{sec:upperbound}
While we are able to exactly compute the $2$-limited broadcast domination number of $P_m \square P_n$, $ P_m \square C_n$, or $ C_m \square C_n$ for many values of $m$ and $n$, we have yet to find a general formula. Thus, we establish upper bounds.  This section describes the method of ``tiling'' which we use in \cite{slobodin} to construct $2$-limited dominating broadcasts on these graphs. Our description uses  examples specific to $P_4 \square P_n$. Results for all $m \geq 2$ and $n \geq m$ are stated subsequently in Sections \ref{sec:paths}, \ref{sec:pathsandcycle}, and \ref{sec:cycles}.

\begin{figure}[htbp]
\centering
\begin{tabular}{cccccc}
 $B$ & $B_1$ & $B_2$ & $B_3$ & $B_4$ & $B_5$\\
\begin{tikzpicture}[scale = 0.45,
baseline={([yshift=-.5ex]current bounding box.center)},
vertex/.style = {circle, fill, inner sep=1.4pt, outer sep=0pt},
every edge quotes/.style = {auto=left, sloped, font=\scriptsize, inner sep=1pt}
]
% LATTICE
\path
(0,0) edge  (0,3)
(1,0) edge  (1,3)
(2,0) edge  (2,3)
(3,0) edge  (3,3)
(4,0) edge  (4,3)
(-1,0) edge   (5,0)
(-1,1) edge   (5,1)
(-1,2) edge   (5,2)
(-1,3) edge   (5,3)
;
\draw[ultra thick, blue] (0,0) rectangle (4,3);

% NODES
\node[vertex] at (2,1) {} ;
\node[vertex] at (0,3) {} ;
\node[vertex] at (4,3) {} ;

%UNDOMINATED

\draw[ForestGreen,ultra thick, fill = white] (0,0) circle  (10pt);
\draw[ForestGreen,ultra thick, fill = white] (4,0) circle  (10pt);

%DOMINATED

\draw[BurntOrange,ultra thick, fill = BurntOrange, fill opacity = 0.5] (-1,3) circle  (10pt);
\draw[BurntOrange,ultra thick, fill = BurntOrange, fill opacity = 0.5] (5,3) circle  (10pt);

\clip (-1,0) rectangle (5,3);

% NEIGHBORHOODS EDGES
\path
%Dominating vertex at (2,1) of strength 2
(4, 1) edge [ultra thick, dotted, red] (2, -1)
%Dominating vertex at (2,1) of strength 2
(2, -1) edge [ultra thick, dotted, red] (0, 1)
%Dominating vertex at (2,1) of strength 2
(0, 1) edge [ultra thick, dotted, red] (2, 3)
%Dominating vertex at (2,1) of strength 2
(2, 3) edge [ultra thick, dotted, red] (4, 1)
%Dominating vertex at (0,3) of strength 1
(1, 3) edge [ultra thick, dotted, red] (0, 2)
%Dominating vertex at (0,3) of strength 1
(0, 2) edge [ultra thick, dotted, red] (-1, 3)
%Dominating vertex at (0,3) of strength 1
(-1, 3) edge [ultra thick, dotted, red] (0, 4)
%Dominating vertex at (0,3) of strength 1
(0, 4) edge [ultra thick, dotted, red] (1, 3)
%Dominating vertex at (4,3) of strength 1
(5, 3) edge [ultra thick, dotted, red] (4, 2)
%Dominating vertex at (4,3) of strength 1
(4, 2) edge [ultra thick, dotted, red] (3, 3)
%Dominating vertex at (4,3) of strength 1
(3, 3) edge [ultra thick, dotted, red] (4, 4)
%Dominating vertex at (4,3) of strength 1
(4, 4) edge [ultra thick, dotted, red] (5, 3)
;
\end{tikzpicture}
&
\begin{tikzpicture}[scale = 0.45,
baseline={([yshift=-.5ex]current bounding box.center)},
vertex/.style = {circle, fill, inner sep=1.4pt, outer sep=0pt},
every edge quotes/.style = {auto=left, sloped, font=\scriptsize, inner sep=1pt}
]
% LATTICE
\path
(0,0) edge  (0,3)
(1,0) edge  (1,3)
(0,0) edge   (2,0)
(0,1) edge   (2,1)
(0,2) edge   (2,2)
(0,3) edge   (2,3)

;
\draw[ultra thick, blue] (0,0) rectangle (1,3);

% NODES
\node[vertex] at (1,0) {} ;
\node[vertex] at (0,2) {} ;

%UNDOMINATED

\draw[ForestGreen,ultra thick, fill = white] (1,3) circle  (10pt);

%DOMINATED

\draw[BurntOrange,ultra thick, fill = BurntOrange, fill opacity = 0.5] (2,0) circle  (10pt);

\clip (0,0) rectangle (2,3);

% NEIGHBORHOODS EDGES
\path
%Dominating vertex at (1,0) of strength 1
(2, 0) edge [ultra thick, dotted, red] (1, -1)
%Dominating vertex at (1,0) of strength 1
(1, -1) edge [ultra thick, dotted, red] (0, 0)
%Dominating vertex at (1,0) of strength 1
(0, 0) edge [ultra thick, dotted, red] (1, 1)
%Dominating vertex at (1,0) of strength 1
(1, 1) edge [ultra thick, dotted, red] (2, 0)
%Dominating vertex at (0,2) of strength 1
(1, 2) edge [ultra thick, dotted, red] (0, 1)
%Dominating vertex at (0,2) of strength 1
(0, 1) edge [ultra thick, dotted, red] (-1, 2)
%Dominating vertex at (0,2) of strength 1
(-1, 2) edge [ultra thick, dotted, red] (0, 3)
%Dominating vertex at (0,2) of strength 1
(0, 3) edge [ultra thick, dotted, red] (1, 2)
;
\end{tikzpicture}

&
\begin{tikzpicture}[scale = 0.45,
baseline={([yshift=-.5ex]current bounding box.center)},
vertex/.style = {circle, fill, inner sep=1.4pt, outer sep=0pt},
every edge quotes/.style = {auto=left, sloped, font=\scriptsize, inner sep=1pt}
]
% LATTICE
\path
(0,0) edge  (0,3)
(-1,0) edge   (0,0)
(-1,1) edge   (0,1)
(-1,2) edge   (0,2)
(-1,3) edge   (0,3)

;
\draw[ultra thick, blue] (0,0) rectangle (0,3);

% NODES
\node[vertex] at (0,1) {} ;

%DOMINATED

\draw[BurntOrange,ultra thick, fill = BurntOrange, fill opacity = 0.5] (-1,0) circle  (10pt);
\draw[BurntOrange,ultra thick, fill = BurntOrange, fill opacity = 0.5] (-2,1) circle  (10pt);
\draw[BurntOrange,ultra thick, fill = BurntOrange, fill opacity = 0.5] (-1,1) circle  (10pt);
\draw[BurntOrange,ultra thick, fill = BurntOrange, fill opacity = 0.5] (-1,2) circle  (10pt);

\clip (-2,0) rectangle (0,3);

% NEIGHBORHOODS EDGES
\path
%Dominating vertex at (0,1) of strength 2
(2, 1) edge [ultra thick, dotted, red] (0, -1)
%Dominating vertex at (0,1) of strength 2
(0, -1) edge [ultra thick, dotted, red] (-2, 1)
%Dominating vertex at (0,1) of strength 2
(-2, 1) edge [ultra thick, dotted, red] (0, 3)
%Dominating vertex at (0,1) of strength 2
(0, 3) edge [ultra thick, dotted, red] (2, 1)
;
\end{tikzpicture}
&
\begin{tikzpicture}[scale = 0.45,
baseline={([yshift=-.5ex]current bounding box.center)},
vertex/.style = {circle, fill, inner sep=1.4pt, outer sep=0pt},
every edge quotes/.style = {auto=left, sloped, font=\scriptsize, inner sep=1pt}
]
% LATTICE
\path
(0,0) edge  (0,3)
(1,0) edge  (1,3)
(2,0) edge  (2,3)
(3,0) edge  (3,3)
(-1,0) edge   (3,0)
(-1,1) edge   (3,1)
(-1,2) edge   (3,2)
(-1,3) edge   (3,3)

;
\draw[ultra thick, blue] (0,0) rectangle (3,3);

% NODES
\node[vertex] at (1,0) {} ;
\node[vertex] at (2,2) {} ;

%UNDOMINATED

\draw[ForestGreen,ultra thick, fill = white] (0,3) circle  (10pt);

%DOMINATED

\draw[BurntOrange,ultra thick, fill = BurntOrange, fill opacity = 0.5] (-1,0) circle  (10pt);

\clip (-1,0) rectangle (3,3);

% NEIGHBORHOODS EDGES
\path
%Dominating vertex at (1,0) of strength 2
(3, 0) edge [ultra thick, dotted, red] (1, -2)
%Dominating vertex at (1,0) of strength 2
(1, -2) edge [ultra thick, dotted, red] (-1, 0)
%Dominating vertex at (1,0) of strength 2
(-1, 0) edge [ultra thick, dotted, red] (1, 2)
%Dominating vertex at (1,0) of strength 2
(1, 2) edge [ultra thick, dotted, red] (3, 0)
%Dominating vertex at (2,2) of strength 2
(4, 2) edge [ultra thick, dotted, red] (2, 0)
%Dominating vertex at (2,2) of strength 2
(2, 0) edge [ultra thick, dotted, red] (0, 2)
%Dominating vertex at (2,2) of strength 2
(0, 2) edge [ultra thick, dotted, red] (2, 4)
%Dominating vertex at (2,2) of strength 2
(2, 4) edge [ultra thick, dotted, red] (4, 2)
;
\end{tikzpicture}
&
\begin{tikzpicture}[scale = 0.45,
baseline={([yshift=-.5ex]current bounding box.center)},
vertex/.style = {circle, fill, inner sep=1.4pt, outer sep=0pt},
every edge quotes/.style = {auto=left, sloped, font=\scriptsize, inner sep=1pt}
]
% LATTICE
\path
(0,0) edge  (0,3)
(1,0) edge  (1,3)
(2,0) edge  (2,3)
(-1,0) edge   (2,0)
(-1,1) edge   (2,1)
(-1,2) edge   (2,2)
(-1,3) edge   (2,3)
;
\draw[ultra thick, blue] (0,0) rectangle (2,3);

% NODES
\node[vertex] at (2,2) {} ;
\node[vertex] at (0,0) {} ;

%UNDOMINATED

\draw[ForestGreen,ultra thick, fill = white] (0,3) circle  (10pt);

%DOMINATED

\draw[BurntOrange,ultra thick, fill = BurntOrange, fill opacity = 0.5] (-1,0) circle  (10pt);

\clip (-1,0) rectangle (2,3);

% NEIGHBORHOODS EDGES
\path
%Dominating vertex at (2,2) of strength 2
(4, 2) edge [ultra thick, dotted, red] (2, 0)
%Dominating vertex at (2,2) of strength 2
(2, 0) edge [ultra thick, dotted, red] (0, 2)
%Dominating vertex at (2,2) of strength 2
(0, 2) edge [ultra thick, dotted, red] (2, 4)
%Dominating vertex at (2,2) of strength 2
(2, 4) edge [ultra thick, dotted, red] (4, 2)
%Dominating vertex at (0,0) of strength 1
(1, 0) edge [ultra thick, dotted, red] (0, -1)
%Dominating vertex at (0,0) of strength 1
(0, -1) edge [ultra thick, dotted, red] (-1, 0)
%Dominating vertex at (0,0) of strength 1
(-1, 0) edge [ultra thick, dotted, red] (0, 1)
%Dominating vertex at (0,0) of strength 1
(0, 1) edge [ultra thick, dotted, red] (1, 0)
;
\end{tikzpicture}
&
\begin{tikzpicture}[scale = 0.45,
baseline={([yshift=-.5ex]current bounding box.center)},
vertex/.style = {circle, fill, inner sep=1.4pt, outer sep=0pt},
every edge quotes/.style = {auto=left, sloped, font=\scriptsize, inner sep=1pt}
]
% LATTICE
\path
(0,0) edge  (0,3)
(1,0) edge  (1,3)
(2,0) edge  (2,3)
(3,0) edge  (3,3)
(4,0) edge  (4,3)
(-1,0) edge   (4,0)
(-1,1) edge   (4,1)
(-1,2) edge   (4,2)
(-1,3) edge   (4,3)
;
\draw[ultra thick, blue] (0,0) rectangle (4,3);

% NODES
\node[vertex] at (2,2) {} ;
\node[vertex] at (0,0) {} ;
\node[vertex] at (4,0) {} ;
\node[vertex] at (4,3) {} ;

%UNDOMINATED

\draw[ForestGreen,ultra thick, fill = white] (0,3) circle  (10pt);

%DOMINATED

\draw[BurntOrange,ultra thick, fill = BurntOrange, fill opacity = 0.5] (-1,0) circle  (10pt);

\clip (-1,0) rectangle (4,3);

% NEIGHBORHOODS EDGES
\path
%Dominating vertex at (2,2) of strength 2
(4, 2) edge [ultra thick, dotted, red] (2, 0)
%Dominating vertex at (2,2) of strength 2
(2, 0) edge [ultra thick, dotted, red] (0, 2)
%Dominating vertex at (2,2) of strength 2
(0, 2) edge [ultra thick, dotted, red] (2, 4)
%Dominating vertex at (2,2) of strength 2
(2, 4) edge [ultra thick, dotted, red] (4, 2)
%Dominating vertex at (0,0) of strength 1
(1, 0) edge [ultra thick, dotted, red] (0, -1)
%Dominating vertex at (0,0) of strength 1
(0, -1) edge [ultra thick, dotted, red] (-1, 0)
%Dominating vertex at (0,0) of strength 1
(-1, 0) edge [ultra thick, dotted, red] (0, 1)
%Dominating vertex at (0,0) of strength 1
(0, 1) edge [ultra thick, dotted, red] (1, 0)
%Dominating vertex at (4,0) of strength 1
(5, 0) edge [ultra thick, dotted, red] (4, -1)
%Dominating vertex at (4,0) of strength 1
(4, -1) edge [ultra thick, dotted, red] (3, 0)
%Dominating vertex at (4,0) of strength 1
(3, 0) edge [ultra thick, dotted, red] (4, 1)
%Dominating vertex at (4,0) of strength 1
(4, 1) edge [ultra thick, dotted, red] (5, 0)
%Dominating vertex at (4,3) of strength 1
(5, 3) edge [ultra thick, dotted, red] (4, 2)
%Dominating vertex at (4,3) of strength 1
(4, 2) edge [ultra thick, dotted, red] (3, 3)
%Dominating vertex at (4,3) of strength 1
(3, 3) edge [ultra thick, dotted, red] (4, 4)
%Dominating vertex at (4,3) of strength 1
(4, 4) edge [ultra thick, dotted, red] (5, 3)
;
\end{tikzpicture}
\end{tabular}
\caption{Tiles used in constructions of $2$-limited dominating broadcasts on $P_4 \square P_n$.}
\label{fig:P4}
\end{figure}
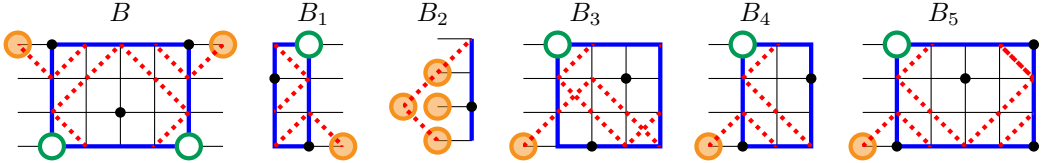

Figure \ref{fig:P4} depicts broadcast tiles $B, B_1, B_2, B_3, B_4,$ and $ B_5$ which can be used to construct $2$-limited dominating broadcasts on $P_4 \square P_n$ for all $n \geq 4$.   
% !!!
% Colour required for figure
% !!!
% !!!
% Colour required for figure
% !!!
The thick blue border indicates the border of each tile. The black circles with a black inner fill indicate vertices on each tile broadcasting at a non-zero strength. In this example, $B$ has two vertices broadcasting at strength one (as found in the top most corners) and one vertex broadcasting at strength two (as found in the lower middle row of the center column). The thick red dotted lines indicate the broadcast ranges of the broadcasting vertices at their centers. Thick green circles with a white inner fill indicate vertices on each tile which do not hear a broadcast within the tile. Thick orange circles with an opaque orange inner fill indicate vertices exterior to the tile which can hear broadcasts from vertices on the tile. 

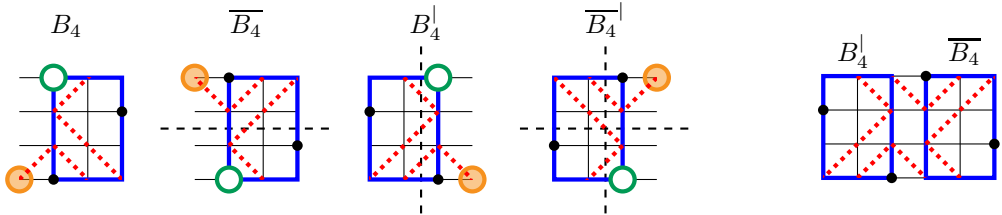
\begin{figure}[htbp]
\centering
\begin{tabular}{cccc}
$B_4$ & $\overline{B_4}$ & $B_4^{|}$ & $\overline{B_4}^{|}$ \\
\begin{tikzpicture}[scale = 0.45,
baseline={([yshift=-.5ex]current bounding box.center)},
vertex/.style = {circle, fill, inner sep=1.4pt, outer sep=0pt},
every edge quotes/.style = {auto=left, sloped, font=\scriptsize, inner sep=1pt}
]
% LATTICE
\path
(0,0) edge  (0,3)
(1,0) edge  (1,3)
(2,0) edge  (2,3)
(-1,0) edge   (2,0)
(-1,1) edge   (2,1)
(-1,2) edge   (2,2)
(-1,3) edge   (2,3)
;
\draw[ultra thick, blue] (0,0) rectangle (2,3);

% NODES
\node[vertex] at (0,0) {} ;
\node[vertex] at (2,2) {} ;

%UNDOMINATED

\draw[ForestGreen,ultra thick, fill = white] (0,3) circle  (10pt);

%DOMINATED

\draw[BurntOrange,ultra thick, fill = BurntOrange, fill opacity = 0.5] (-1,0) circle  (10pt);

\clip (-1,0) rectangle (2,3);

% NEIGHBORHOODS EDGES
\path
%Dominating vertex at (0,0) of strength 1
(1, 0) edge [ultra thick, dotted, red] (0, -1)
%Dominating vertex at (0,0) of strength 1
(0, -1) edge [ultra thick, dotted, red] (-1, 0)
%Dominating vertex at (0,0) of strength 1
(-1, 0) edge [ultra thick, dotted, red] (0, 1)
%Dominating vertex at (0,0) of strength 1
(0, 1) edge [ultra thick, dotted, red] (1, 0)
%Dominating vertex at (2,2) of strength 2
(4, 2) edge [ultra thick, dotted, red] (2, 0)
%Dominating vertex at (2,2) of strength 2
(2, 0) edge [ultra thick, dotted, red] (0, 2)
%Dominating vertex at (2,2) of strength 2
(0, 2) edge [ultra thick, dotted, red] (2, 4)
%Dominating vertex at (2,2) of strength 2
(2, 4) edge [ultra thick, dotted, red] (4, 2)
;
\end{tikzpicture}
        &
\begin{tikzpicture}[scale = 0.45,
baseline={([yshift=-.5ex]current bounding box.center)},
vertex/.style = {circle, fill, inner sep=1.4pt, outer sep=0pt},
every edge quotes/.style = {auto=left, sloped, font=\scriptsize, inner sep=1pt}
]
% LATTICE
\path
(0,0) edge  (0,3)
(1,0) edge  (1,3)
(2,0) edge  (2,3)
(-1,0) edge   (2,0)
(-1,1) edge   (2,1)
(-1,2) edge   (2,2)
(-1,3) edge   (2,3)
(-2, 1.5) edge [thick, dashed] (3, 1.5)
;
\draw[ultra thick, blue] (0,0) rectangle (2,3);

% NODES
\node[vertex] at (2,1) {} ;
\node[vertex] at (0,3) {} ;

%UNDOMINATED

\draw[ForestGreen,ultra thick, fill = white] (0,0) circle  (10pt);

%DOMINATED

\draw[BurntOrange,ultra thick, fill = BurntOrange, fill opacity = 0.5] (-1,3) circle  (10pt);

\clip (-1,0) rectangle (2,3);

% NEIGHBORHOODS EDGES
\path
%Dominating vertex at (2,1) of strength 2
(4, 1) edge [ultra thick, dotted, red] (2, -1)
%Dominating vertex at (2,1) of strength 2
(2, -1) edge [ultra thick, dotted, red] (0, 1)
%Dominating vertex at (2,1) of strength 2
(0, 1) edge [ultra thick, dotted, red] (2, 3)
%Dominating vertex at (2,1) of strength 2
(2, 3) edge [ultra thick, dotted, red] (4, 1)
%Dominating vertex at (0,3) of strength 1
(1, 3) edge [ultra thick, dotted, red] (0, 2)
%Dominating vertex at (0,3) of strength 1
(0, 2) edge [ultra thick, dotted, red] (-1, 3)
%Dominating vertex at (0,3) of strength 1
(-1, 3) edge [ultra thick, dotted, red] (0, 4)
%Dominating vertex at (0,3) of strength 1
(0, 4) edge [ultra thick, dotted, red] (1, 3)
;
\end{tikzpicture}

&

\begin{tikzpicture}[scale = 0.45,
baseline={([yshift=-.5ex]current bounding box.center)},
vertex/.style = {circle, fill, inner sep=1.4pt, outer sep=0pt},
every edge quotes/.style = {auto=left, sloped, font=\scriptsize, inner sep=1pt}
]
% LATTICE
\path
(1,0) edge  (1,3)
(2,0) edge  (2,3)
(3,0) edge  (3,3)
(1,0) edge   (4,0)
(1,1) edge   (4,1)
(1,2) edge   (4,2)
(1,3) edge   (4,3)

(2.5, -1) edge [thick, dashed] (2.5, 4)
;
\draw[ultra thick, blue] (1,0) rectangle (3,3);

% NODES
\node[vertex] at (3,0) {} ;
\node[vertex] at (1,2) {} ;

%UNDOMINATED

\draw[ForestGreen,ultra thick, fill = white] (3,3) circle  (10pt);

%DOMINATED

\draw[BurntOrange,ultra thick, fill = BurntOrange, fill opacity = 0.5] (4,0) circle  (10pt);

\clip (1,0) rectangle (4,3);

% NEIGHBORHOODS EDGES
\path
%Dominating vertex at (3,0) of strength 1
(4, 0) edge [ultra thick, dotted, red] (3, -1)
%Dominating vertex at (3,0) of strength 1
(3, -1) edge [ultra thick, dotted, red] (2, 0)
%Dominating vertex at (3,0) of strength 1
(2, 0) edge [ultra thick, dotted, red] (3, 1)
%Dominating vertex at (3,0) of strength 1
(3, 1) edge [ultra thick, dotted, red] (4, 0)
%Dominating vertex at (1,2) of strength 2
(3, 2) edge [ultra thick, dotted, red] (1, 0)
%Dominating vertex at (1,2) of strength 2
(1, 0) edge [ultra thick, dotted, red] (-1, 2)
%Dominating vertex at (1,2) of strength 2
(-1, 2) edge [ultra thick, dotted, red] (1, 4)
%Dominating vertex at (1,2) of strength 2
(1, 4) edge [ultra thick, dotted, red] (3, 2)

;
\end{tikzpicture}
&
\begin{tikzpicture}[scale = 0.45,
baseline={([yshift=-.5ex]current bounding box.center)},
vertex/.style = {circle, fill, inner sep=1.4pt, outer sep=0pt},
every edge quotes/.style = {auto=left, sloped, font=\scriptsize, inner sep=1pt}
]
% LATTICE
\path
(1,0) edge  (1,3)
(2,0) edge  (2,3)
(3,0) edge  (3,3)
(1,0) edge   (4,0)
(1,1) edge   (4,1)
(1,2) edge   (4,2)
(1,3) edge   (4,3)

(2.5, -1) edge [thick, dashed] (2.5, 4)

(0, 1.5) edge [thick, dashed] (5, 1.5)
;
\draw[ultra thick, blue] (1,0) rectangle (3,3);

% NODES
\node[vertex] at (3,3) {} ;
\node[vertex] at (1,1) {} ;

%UNDOMINATED

\draw[ForestGreen,ultra thick, fill = white] (3,0) circle  (10pt);

%DOMINATED

\draw[BurntOrange,ultra thick, fill = BurntOrange, fill opacity = 0.5] (4,3) circle  (10pt);

\clip (1,0) rectangle (4,3);

% NEIGHBORHOODS EDGES
\path
%Dominating vertex at (3,3) of strength 1
(4, 3) edge [ultra thick, dotted, red] (3, 2)
%Dominating vertex at (3,3) of strength 1
(3, 2) edge [ultra thick, dotted, red] (2, 3)
%Dominating vertex at (3,3) of strength 1
(2, 3) edge [ultra thick, dotted, red] (3, 4)
%Dominating vertex at (3,3) of strength 1
(3, 4) edge [ultra thick, dotted, red] (4, 3)
%Dominating vertex at (1,1) of strength 2
(3, 1) edge [ultra thick, dotted, red] (1, -1)
%Dominating vertex at (1,1) of strength 2
(1, -1) edge [ultra thick, dotted, red] (-1, 1)
%Dominating vertex at (1,1) of strength 2
(-1, 1) edge [ultra thick, dotted, red] (1, 3)
%Dominating vertex at (1,1) of strength 2
(1, 3) edge [ultra thick, dotted, red] (3, 1)

;
\end{tikzpicture}
\end{tabular}
\hspace{1cm}
\begin{tabular}{c}
$B_4^{|}$ \hspace{0.8cm} $\overline{B_4}$\\
\begin{tikzpicture}[scale = 0.45,
baseline={([yshift=-.5ex]current bounding box.center)},
vertex/.style = {circle, fill, inner sep=1.4pt, outer sep=0pt},
every edge quotes/.style = {auto=left, sloped, font=\scriptsize, inner sep=1pt}
]
% LATTICE
\path
(1,0) edge  (1,3)
(2,0) edge  (2,3)
(3,0) edge  (3,3)
(4,0) edge  (4,3)
(5,0) edge  (5,3)
(6,0) edge  (6,3)
(1,0) edge   (6,0)
(1,1) edge   (6,1)
(1,2) edge   (6,2)
(1,3) edge   (6,3)

;
\draw[ultra thick, blue] (1,0) rectangle (3,3);
\draw[ultra thick, blue] (4,0) rectangle (6,3);

% NODES
\node[vertex] at (3,0) {} ;
\node[vertex] at (1,2) {} ;
\node[vertex] at (4,3) {} ;
\node[vertex] at (6,1) {} ;

%UNDOMINATED

\clip (1,0) rectangle (6,3);

% NEIGHBORHOODS EDGES
\path
%Dominating vertex at (3,0) of strength 1
(4, 0) edge [ultra thick, dotted, red] (3, -1)
%Dominating vertex at (3,0) of strength 1
(3, -1) edge [ultra thick, dotted, red] (2, 0)
%Dominating vertex at (3,0) of strength 1
(2, 0) edge [ultra thick, dotted, red] (3, 1)
%Dominating vertex at (3,0) of strength 1
(3, 1) edge [ultra thick, dotted, red] (4, 0)
%Dominating vertex at (1,2) of strength 2
(3, 2) edge [ultra thick, dotted, red] (1, 0)
%Dominating vertex at (1,2) of strength 2
(1, 0) edge [ultra thick, dotted, red] (-1, 2)
%Dominating vertex at (1,2) of strength 2
(-1, 2) edge [ultra thick, dotted, red] (1, 4)
%Dominating vertex at (1,2) of strength 2
(1, 4) edge [ultra thick, dotted, red] (3, 2)
%Dominating vertex at (4,3) of strength 1
(5, 3) edge [ultra thick, dotted, red] (4, 2)
%Dominating vertex at (4,3) of strength 1
(4, 2) edge [ultra thick, dotted, red] (3, 3)
%Dominating vertex at (4,3) of strength 1
(3, 3) edge [ultra thick, dotted, red] (4, 4)
%Dominating vertex at (4,3) of strength 1
(4, 4) edge [ultra thick, dotted, red] (5, 3)
%Dominating vertex at (6,1) of strength 2
(8, 1) edge [ultra thick, dotted, red] (6, -1)
%Dominating vertex at (6,1) of strength 2
(6, -1) edge [ultra thick, dotted, red] (4, 1)
%Dominating vertex at (6,1) of strength 2
(4, 1) edge [ultra thick, dotted, red] (6, 3)
%Dominating vertex at (6,1) of strength 2
(6, 3) edge [ultra thick, dotted, red] (8, 1)

;
\end{tikzpicture}
\end{tabular}
\caption{(Four left most tiles) Possible states of tile $B_4$. (Right) $B_4^{|}\overline{B}$ given the tile $B_4$.}
\label{fig:p4states}
\end{figure}

Given the tile $B_4$, $\overline{B_4}$ denotes flipping $B_4$ about the horizontal axis, $B_4^{|}$ indicates flipping $B_4$ about the vertical axis, and $\overline{B_4}^{|}$ indicates flipping $B_4$ about the horizontal and the vertical axes. The four possible states of $B_4$ are shown in Figure \ref{fig:p4states} (Four left most tiles).
% !!!
% Colour required for figure
% !!!
% !!!
% Colour required for figure
% !!!

Given two tiles $T_1$ and $T_2$, the sequence $T_1T_2$ denotes placing the right side of $T_1$ beside the left side of $T_2$ such that the vertices along the right border of $T_1$ and the left border of $T_2$ are at distance one from one another. Figure  \ref{fig:p4states} (Right) depicts $B_4^{|}\overline{B_4}$. Observe that, in Figure \ref{fig:p4states} (Right), the undominated vertex in the upper right corner of $B_4^{|}$ hears the broadcast from the vertex in the upper left corner of $\overline{B_4}$. Similarly, in Figure \ref{fig:p4states} (Right), the undominated vertex in the bottom left corner of $\overline{B_4}$ hears the broadcast from the vertex in the bottom right corner of $B_4^{|}$. As all vertices of $ P_4 \square P_6$ in Figure \ref{fig:p4states} (Right) hear a broadcast under the tiling $B_4^{|}\overline{B_4}$,  the sequence $B_4^{|}\overline{B_4}$ gives a $2$-limited dominating broadcast scheme on $P_4 \square P_6$ of cost six. Therefore, $\gamma_{b,2} \left(P_4 \square P_6 \right) \leq 6$. 

As an example of our method in Section \ref{sec:multipack}, we construct a valid $2$-limited fractional multipacking on $P_4 \square P_6$. For each vertex $i \in V(P_4 \square P_6)$, define the variable $y_i$. Define the fractional $2$-limited multipacking on $P_4 \square P_6$ by $y_i = 1/3 $ for all $i \in V(P_4 \square P_6)$ such that $i$ is on the boundary of $P_4 \square P_6$ and $y_i = 0$ otherwise. This assignment satisfies the constraints of \ref{LP:multipackc} with $k=2$. Thus, $16/3 \leq mp_{f,2}(P_4 \square P_6)$. By the duality theorem of linear programming, $mp_{f,2}(P_4 \square P_6) \leq \gamma_{b,2}(P_4 \square P_6)$. As $ \gamma_{b,2}(P_4 \square P_6)$ is an integer, it follows that $6 \leq  \gamma_{b,2}(P_4 \square P_6)$. Thus, $\gamma_{b,2}(P_4 \square P_6) = 6$. 

For each least residue of $n$ modulo 10, Table \ref{table:P4} contains the associated sequences of broadcast tiles from Figure \ref{fig:P4} which give a  $2$-limited dominating broadcast on $P_4 \square P_n$. The costs of the constructions in Table \ref{table:P4} establish Proposition \ref{prop:P4}. The costs of these constructions are easily calculated from the costs of the individual broadcast tiles in Figure  \ref{table:P4}. The values of Proposition \ref{prop:P4} are optimal for all $n \leq 99$ (as found by computation). 

\begin{table}[htbp]
\renewcommand\arraystretch{1.6}
\centering
\begin{tabular}{r|c|l}
 & Construction & \multicolumn{1}{c}{Cost} \\ \hline
$n \equiv 0 \pmod{10} $: & $B_1 \left[ B \overline{B} \right]^{ \frac{ n-10 }{ 10 } }BB_4$  & $ 8 \left( \frac{ n-10 }{ 10 }  \right) +9 = 8 \left( \frac{ n }{ 10 }  \right) + 1$\\
$n \equiv 1 \pmod{10} $: & $B_1 \left[ B \overline{B} \right]^{ \frac{ n-11 }{ 10 } }BB_3$ & $ 8 \left( \frac{ n-11 }{ 10 }  \right) + 10 = 8 \left( \frac{ n-1 }{ 10 }  \right) + 2  $\\
$n \equiv 2 \pmod{10}$: & $B_1 \left[ B \overline{B} \right]^{ \frac{ n-12 }{ 10 } }BB_5  $ & $ 8 \left( \frac{ n-12 }{ 10 }  \right) + 11 = 8 \left( \frac{ n-2 }{ 10 }  \right) + 3 $\\
$n \equiv 3 \pmod{10}$: & $B_1 \left[ B\overline{B} \right]^{ \frac{ n-3 }{ 10 } } \overline{B_2}  $ & $ 8 \left( \frac{ n-3 }{ 10 }  \right) + 4$\\
$n \equiv 4 \pmod{10} $: & $B_1 \left[ B\overline{B} \right]^{ \frac{ n-4 }{ 10 } }\overline{B_1}^{|} $ & $ 8 \left( \frac{ n-4 }{ 10 }  \right) +4 $\\
$n \equiv 5 \pmod{10} $: & $B_1 \left[ B \overline{B} \right]^{ \frac{ n-5 }{ 10 } }\overline{B_4}$  & $ 8 \left( \frac{ n-5 }{ 10 }  \right) +5$\\
$n \equiv 6 \pmod{10} $: & $B_1 \left[ B \overline{B} \right]^{ \frac{ n-6 }{ 10 } }\overline{B_3}$ & $ 8 \left( \frac{ n-6 }{ 10 }  \right) + 6 $\\
$n \equiv 7 \pmod{10} $: & $B_1 \left[ B \overline{B} \right]^{ \frac{ n-7 }{ 10 } }\overline{B_5}  $ & $ 8 \left( \frac{ n-7 }{ 10 }  \right) + 7  $\\
$n \equiv 8 \pmod{10} $: & $B_1 \left[ B\overline{B} \right]^{ \frac{ n-8 }{ 10 } }B B_2  $ & $ 8 \left( \frac{ n-8 }{ 10 }  \right) + 8$ \\
$n \equiv 9 \pmod{10} $: & $B_1 \left[ B\overline{B} \right]^{ \frac{ n-9 }{ 10 } }BB_1^{|} $ & $ 8 \left( \frac{ n-9 }{ 10 }  \right) +8 $\\
\end{tabular}
\caption{Constructions of $2$-limited dominating broadcasts on $P_4 \square P_n$ using the tiles from Figure \ref{fig:P4}.}
\label{table:P4}
\end{table}
% default
\renewcommand\arraystretch{1}

\begin{prop}\label{prop:P4}
For $n \geq 4$, $\gamma_{b,2} \left( P_4 \square P_n \right) \leq 8 \left\lfloor \frac{ n }{ 10 }  \right\rfloor + d(n_{10})$ where $n_{10}$ is  the least residue  of $n$ modulo 10 and $d(n_{10})$ is given in Table \ref{table:pathxpathExample}.
\begin{table}[htbp]
\centering
\begin{tabular}{c|ccccccccccc}
$n_{10}$ & $0$ & $1$ & $2$ & $3$ & $4$ & $5$ & $6$ & $7$ & $8$ & $9$ \\ \hline
$d(n_{10}) $& $1$ & $2$ & $3$ & $4$ & $4$ & $5$ & $6$ & $7$ & $8$ & $8$       \\
\end{tabular}
\caption{Values of $d(n_{10})$ in the upper bound of $ \gamma_{b,2} \left(P_4 \square P_n\right)$ for $n \geq 4$.}
\label{table:pathxpathExample}
\end{table}
\end{prop}

Tilings are used in Sections \ref{sec:paths}, \ref{sec:pathsandcycle}, and \ref{sec:cycles} to establish upper bounds for the $2$-limited broadcast domination number on all $ P_m \square P_n$, $P_m \square C_n$, and $C_m \square C_n$ grid graphs, respectively. This paper does not include the broadcast tiles to prove these bounds as, when combined, they span over 100 pages. For each stated theorem, we include the section(s) of \cite{slobodin} which contain the appropriate broadcast tiles. Most broadcast tiles of \cite{slobodin} were identified by hand by examining known optimal solutions, as found by computation, and looking for recurring patterns. This manual process was greatly facilitated by automatically creating PDF visuals of known optimal solutions using a Python script (written by Slobodin).  

\section{Upper Bounds for $\gamma_{b,2} \left(P_m \square P_n\right)$}\label{sec:paths}
In this section, we establish upper bounds for the $2$-limited broadcast domination number of the Cartesian product of two paths. Section \ref{sec:pathfirst} states upper bounds for $\gamma_{b,2}(P_m \square P_n)$ for $2 \leq m \leq 12$ and all $n \geq m$. These results  are established using the tiling method described in Section \ref{sec:upperbound}.  For $m,n \geq 13$, Section \ref{path:sec:p13} describes a general construction  on $ P_m \square P_n$ derived from a $2$-limited dominating broadcast on the Cartesian plane.

\subsection{$P_{2 \leq m \leq 12} \square P_{n\geq m}$}\label{sec:pathfirst}
The tilings used to establish Theorem \ref{thm:path} are in \cite[Sections 2.2 through 2.12]{slobodin}.

\begin{thm}\label{thm:path}
Fix $2 \leq m \leq 12$ and $n \geq m$. Let $x$ be the value in Table \ref{table:pathxpathx} dependant upon $m$ and define $n_x$ as the least residue of $n$ modulo $x$. By construction,  
\begin{equation*}
\begin{aligned}
\gamma_{b,2} \left(P_m \square P_n\right) \leq b(m) + c(m,n,n_x)
\end{aligned}
\end{equation*}
where $b(m)$ corresponds with the $O(n)$ terms in Table \ref{table:pathxpathlinearterm} and $c(m,n,n_x)$ corresponds with the $O(1)$ terms in Table \ref{table:pathxpathconstanterm} (see Section \ref{appendix}).
\end{thm}

\begin{table}[htbp]
\centering
\begin{tabular}{c|cccccccccccccc}
    $m$ &  $2,3$ &$4$  & $5$ & $6$ & $7$ & $8$ & $9$ & $10$ & $11$ & $12$\\ \hline
    $x $ & $1$ & $10$ & $1$ & $16$ & $14$ & $22$ & $10$ & $18$ & $26$ & $24$ \\
\end{tabular}
\caption{Value of $x$ in the upper bound of $ \gamma_{b,2} \left(P_m \square P_n\right)$ for $2 \leq m \leq 12$ and $n \geq m$.}
\label{table:pathxpathx}
\end{table}

\renewcommand\arraystretch{1.25}
\begin{table}[htbp]
\centering
    \begin{tabular}{c|ccccccccccc}
    $m$ & $2$ & $3$ & $4$ & $5$ & $6$ & $7$ & $8$ & $9$ & $10$ & $11$ & $12$\\ \hline
    $b(m) $&        $ \left\lfloor \frac{n}{2}  \right\rfloor  $ & $ \left\lceil \frac{2n}{3}  \right\rceil $ & $ 8 \left\lfloor \frac{n}{10}  \right\rfloor $ & $n$ & $ 18 \left\lfloor \frac{n}{16}  \right\rfloor $ & $ 18 \left\lfloor \frac{n}{14}  \right\rfloor $ & $32 \left\lfloor \frac{n}{22}  \right\rfloor  $ & $ 16 \left\lfloor \frac{n}{10}  \right\rfloor  $ & $ 32 \left\lfloor  \frac{n}{18}  \right\rfloor $ & $ 50 \left\lfloor \frac{n}{26}  \right\rfloor  $ & $ 50 \left\lfloor \frac{n}{ 24}  \right\rfloor $ \\
\end{tabular}
\caption{Value of $b(m)$ in the upper bound of $ \gamma_{b,2} \left(P_m \square P_n\right)$ for $2 \leq m \leq 12$ and $n \geq m$.}
\label{table:pathxpathlinearterm}
\end{table}
% default
\renewcommand\arraystretch{1}

\subsection{$P_{m\geq13} \square P_{n\geq13}$}\label{path:sec:p13}
This section provides upper bounds for $\gamma_{b,2} \left( P_m \square P_n \right) $ for $m,n \geq 13$ by modifying  $2$-limited dominating broadcasts on the Cartesian plane. Our expectation is that this method will achieve good upper bounds  as, for very large $m$ and $n$, there likely exists a $2$-limited dominating broadcast on $P_m \square P_n$ which, for much of the graph, resembles this $2$-limited dominating broadcast on the Cartesian plane. This method was inspired by  the work on the distance domination number of grids in \cite{undergrad}.

\subsubsection{$2$-limited dominating broadcasts on $\mathbb{Z}^2$}\label{subsec:$2$-limitedd}
Let $ \mathbb{Z}\times \mathbb{Z} = \mathbb{Z}^2$ denote the integer lattice. Let $G_{m,n}$ and $Y_{m+4,n+4}$ denote $P_m \square P_n$ and $P_{m+4} \square P_{n+4}$, respectively, embedded in $ \mathbb{Z}^2$ with their bottom left corners at $(0,0)$ and $(-2,-2)$, respectively. That is, let  
\begin{equation*}
\begin{aligned}
V \left( G_{m,n} \right)  &= \left\{ (i,j) \in \mathbb{Z}^2 : 0 \leq i \leq n-1 \textnormal{ and } 0 \leq j \leq m-1 \right\} \textnormal{ and}  \\
V \left( Y_{m+4,n+4} \right)  &= \left\{ (i,j) \in \mathbb{Z}^2 : -2 \leq i \leq n+1 \textnormal{ and } -2 \leq j \leq m+1 \right\} .\\
\end{aligned}
\end{equation*}
Let $ \mathbb{Z}_{13} = \left\{ 0, 1, \dots,12 \right\} $.  Define the function $ \phi: \mathbb{Z}\times \mathbb{Z}\mapsto \mathbb{Z}_{13}$ by 
\begin{equation*}
\begin{aligned}
\left( i,j \right) \mapsto 3i + 2j \pmod{ 13 } .
\end{aligned}
\end{equation*}
Fix $ \ell \in \mathbb{Z}_{13}$. The set given by
\begin{equation*}
\begin{aligned}
\phi^{-1} \left( \ell \right)  &=  \left\{ (i,j) \in \mathbb{Z}^2: 3i+2j \equiv \ell \pmod{ 13 } \right\} 
\end{aligned}
\end{equation*}
defines the vertices of $ \mathbb{Z}^2$ which, if broadcasting at strength two, give a $2$-limited dominating broadcast on $ \mathbb{Z}^2$ \cite[Lemma V.7]{undergrad2}. Figure \ref{fig:optimalp13} depicts $G_{13,13}$, $Y_{17,17}$, and the broadcast vertices and their neighbourhoods given by $\phi^{-1} \left( 8 \right) $. The thick blue and green lines indicate the borders of $G_{13,13}$ and $Y_{17,17}$, respectively. The light blue fill and green inner fill indicate the interiors of $G_{13,13}$ and $Y_{17,17}$, respectively. The thick black lines indicate the axes of $ \mathbb{Z}^2$. To see that the broadcast set in Figure \ref{fig:optimalp13} is given by $ \phi^{-1}(8)$ observe, for instance, that $(2,1)$ is one such broadcaster and 
\begin{equation*}
\begin{aligned}
\phi(2,1)=   3(2)+2(1)=8 &  \equiv 8 \pmod{13}  \Rightarrow (5,1) \in \phi^{-1}(8).
\end{aligned}
\end{equation*}

Prior to constructing $2$-limited dominating broadcasts on $P_m \square P_n$ for all $m,n \geq 13$ from  $2$-limited dominating broadcasts on $ \mathbb{Z}^2$, we include Example \ref{example:pxpinf} to provide an introduction to our methodology.

\begin{eg}\label{example:pxpinf}
Suppose we wish to find a $2$-limited dominating broadcast on $G_{13,13}$. Overlay one of the 13 possible $2$-limited dominating broadcasts of $ \mathbb{Z}^2$  given by $ \phi^{-1} \left( \ell \right) $ for some $\ell \in \mathbb{Z}_{13}$. The broadcast vertices on $Y_{17,17}$ dominate the vertices of $G_{13,13}$. Figure \ref{fig:optimalp13} (Left) depicts $G_{13,13}$, $Y_{17,17}$, and the overlay of $ \phi^{-1} \left( 8 \right) $ where the thick purple circles indicate the broadcast vertices on $Y_{17,17}$ but not $G_{13,13}$. 

% !!!
% Colour required for figure
% !!!
\begin{figure}[htbp]
\centering
\begin{tikzpicture}[scale = 0.35,
baseline={([yshift=-.5ex]current bounding box.center)},
vertex/.style = {circle, fill, inner sep=1.4pt, outer sep=0pt},
every edge quotes/.style = {auto=left, sloped, font=\scriptsize, inner sep=1pt}
]

\node[vertex] at (0,0) {} ;
\node[vertex] at (13,0) {} ;
\node[vertex] at (8,1) {} ;
\node[vertex] at (3,2) {} ;
\node[vertex] at (11,3) {} ;
\node[vertex] at (6,4) {} ;
\node[vertex] at (1,5) {} ;
\node[vertex] at (14,5) {} ;
\node[vertex] at (9,6) {} ;
\node[vertex] at (4,7) {} ;
\node[vertex] at (12,8) {} ;
\node[vertex] at (7,9) {} ;
\node[vertex] at (2,10) {} ;
\node[vertex] at (15,10) {} ;
\node[vertex] at (10,11) {} ;
\node[vertex] at (5,12) {} ;
\node[vertex] at (0,13) {} ;
\node[vertex] at (13,13) {} ;
\node[vertex] at (3,15) {}; 
\node[vertex]  at (5,-1) {}; 
\node[vertex]  at (8,14)  {}; 
\node[vertex]  at (-1,8)  {}; 

\draw [DarkOrchid, ultra thick] (0,0) circle (10pt);
\draw [DarkOrchid, ultra thick] (5,-1) circle (10pt);
\draw [DarkOrchid, ultra thick] (13,0) circle (10pt);
\draw [DarkOrchid, ultra thick] (14,5) circle (10pt);
\draw [DarkOrchid, ultra thick] (15,10) circle (10pt);
\draw [DarkOrchid, ultra thick] (3,15) circle (10pt);
\draw [DarkOrchid, ultra thick] (8,14) circle (10pt);
\draw [DarkOrchid, ultra thick] (0,13) circle (10pt);
\draw [DarkOrchid, ultra thick] (-1,8) circle (10pt);

\path
(5,1) edge [red, ultra thick, dotted] (7,-1)
(7,-1) edge [red, ultra thick, dotted] (6.5,-1.5)
(3.5,-1.5) edge [red, ultra thick, dotted] (3,-1)
(3,-1) edge [red, ultra thick, dotted] (5,1) 
(1.5,-1.5) edge [red, ultra thick, dotted] (2,-1)
(2,-1) edge [red, ultra thick, dotted] (2.5,-1.5)
(8.5,-1.5) edge [red, ultra thick, dotted] (10,0)
(10,0) edge [red, ultra thick, dotted] (11.5,-1.5)
(14.5,-1.5) edge [red, ultra thick, dotted] (15,-1)
(15,-1) edge [red, ultra thick, dotted] (15.5,-1.5)
(11,14) edge [red, ultra thick, dotted] (9.5,15.5)
(11,14) edge [red, ultra thick, dotted] (12.5,15.5)
(3,13) edge [red, ultra thick, dotted] (1,15)
(1,15) edge [red, ultra thick, dotted] (1.5,15.5)
(5,15) edge [red, ultra thick, dotted] (4.5,15.5)
(3,13) edge [red, ultra thick, dotted] (5,15)
(8,12) edge [red, ultra thick, dotted] (6,14)
(8,12) edge [red, ultra thick, dotted] (10,14)
(6,14) edge [red, ultra thick, dotted] (7.5,15.5)
(10,14) edge [red, ultra thick, dotted] (8.5,15.5)
(6,15) edge [red, ultra thick, dotted] (6.5,15.5)
(6,15) edge [red, ultra thick, dotted] (5.5,15.5)
(-1.5,1.5) edge [red, ultra thick, dotted] (0,3)
(0,3) edge [red, ultra thick, dotted] (-1.5,4.5)
(-1,6) edge [red, ultra thick, dotted] (1,8)
(1,8) edge [red, ultra thick, dotted] (-1,10)
(-1,10) edge [red, ultra thick, dotted] (-1.5,9.5)
(-1,6) edge [red, ultra thick, dotted] (-1.5,6.5)
(-1.5,10.5) edge [red, ultra thick, dotted] (-1,11)
(-1,11) edge [red, ultra thick, dotted] (-1.5,11.5)
(14,15) edge [red, ultra thick, dotted] (15.5,13.5)
(14,15) edge [red, ultra thick, dotted] (14.5,15.5)
(-0.5,15.5) edge [red, ultra thick, dotted] (-1.5,14.5)
(0,2) edge [red, ultra thick, dotted] (2,0)
(2,0) edge [red, ultra thick, dotted] (0.5,-1.5)
(-0.5,-1.5) edge [red, ultra thick, dotted] (-1.5,-0.5)
(-1.5,0.5) edge [red, ultra thick, dotted] (0,2) 
(13,2) edge [red, ultra thick, dotted] (15,0)
(15,0) edge [red, ultra thick, dotted] (13.5,-1.5)
(12.5,-1.5) edge [red, ultra thick, dotted] (11,0)
(11,0) edge [red, ultra thick, dotted] (13,2) 
(8,3) edge [red, ultra thick, dotted] (10,1)
(10,1) edge [red, ultra thick, dotted] (8,-1)
(8,-1) edge [red, ultra thick, dotted] (6,1)
(6,1) edge [red, ultra thick, dotted] (8,3) 
(3,4) edge [red, ultra thick, dotted] (5,2)
(5,2) edge [red, ultra thick, dotted] (3,0)
(3,0) edge [red, ultra thick, dotted] (1,2)
(1,2) edge [red, ultra thick, dotted] (3,4) 
(15.5,0.5) edge [red, ultra thick, dotted] (14,2)
(14,2) edge [red, ultra thick, dotted] (15.5,3.5) 
(11,5) edge [red, ultra thick, dotted] (13,3)
(13,3) edge [red, ultra thick, dotted] (11,1)
(11,1) edge [red, ultra thick, dotted] (9,3)
(9,3) edge [red, ultra thick, dotted] (11,5) 
(6,6) edge [red, ultra thick, dotted] (8,4)
(8,4) edge [red, ultra thick, dotted] (6,2)
(6,2) edge [red, ultra thick, dotted] (4,4)
(4,4) edge [red, ultra thick, dotted] (6,6) 
(1,7) edge [red, ultra thick, dotted] (3,5)
(3,5) edge [red, ultra thick, dotted] (1,3)
(1,3) edge [red, ultra thick, dotted] (-1,5)
(-1,5) edge [red, ultra thick, dotted] (1,7) 
(14,7) edge [red, ultra thick, dotted] (15.5,5.5)
(15.5,4.5) edge [red, ultra thick, dotted] (14,3)
(14,3) edge [red, ultra thick, dotted] (12,5)
(12,5) edge [red, ultra thick, dotted] (14,7) 
(9,8) edge [red, ultra thick, dotted] (11,6)
(11,6) edge [red, ultra thick, dotted] (9,4)
(9,4) edge [red, ultra thick, dotted] (7,6)
(7,6) edge [red, ultra thick, dotted] (9,8) 
(4,9) edge [red, ultra thick, dotted] (6,7)
(6,7) edge [red, ultra thick, dotted] (4,5)
(4,5) edge [red, ultra thick, dotted] (2,7)
(2,7) edge [red, ultra thick, dotted] (4,9) 
(15.5,6.5) edge [red, ultra thick, dotted] (15,7)
(15,7) edge [red, ultra thick, dotted] (15.5,7.5) 
(12,10) edge [red, ultra thick, dotted] (14,8)
(14,8) edge [red, ultra thick, dotted] (12,6)
(12,6) edge [red, ultra thick, dotted] (10,8)
(10,8) edge [red, ultra thick, dotted] (12,10) 
(7,11) edge [red, ultra thick, dotted] (9,9)
(9,9) edge [red, ultra thick, dotted] (7,7)
(7,7) edge [red, ultra thick, dotted] (5,9)
(5,9) edge [red, ultra thick, dotted] (7,11) 
(2,12) edge [red, ultra thick, dotted] (4,10)
(4,10) edge [red, ultra thick, dotted] (2,8)
(2,8) edge [red, ultra thick, dotted] (0,10)
(0,10) edge [red, ultra thick, dotted] (2,12) 
(15,12) edge [red, ultra thick, dotted] (15.5,11.5)
(15.5,8.5) edge [red, ultra thick, dotted] (15,8)
(15,8) edge [red, ultra thick, dotted] (13,10)
(13,10) edge [red, ultra thick, dotted] (15,12) 
(10,13) edge [red, ultra thick, dotted] (12,11)
(12,11) edge [red, ultra thick, dotted] (10,9)
(10,9) edge [red, ultra thick, dotted] (8,11)
(8,11) edge [red, ultra thick, dotted] (10,13) 
(5,14) edge [red, ultra thick, dotted] (7,12)
(7,12) edge [red, ultra thick, dotted] (5,10)
(5,10) edge [red, ultra thick, dotted] (3,12)
(3,12) edge [red, ultra thick, dotted] (5,14) 
(0,15) edge [red, ultra thick, dotted] (2,13)
(2,13) edge [red, ultra thick, dotted] (0,11)
(0,11) edge [red, ultra thick, dotted] (-1.5,12.5)
(-1.5,13.5) edge [red, ultra thick, dotted] (0,15) 
(13,15) edge [red, ultra thick, dotted] (15,13)
(15,13) edge [red, ultra thick, dotted] (13,11)
(13,11) edge [red, ultra thick, dotted] (11,13)
(11,13) edge [red, ultra thick, dotted] (13,15)
(-1,-1.5) edge[thick]  (-1,15.5)
(0,-1.5) edge [thick] (0,15.5)
(1,-1.5) edge [thick] (1,15.5)
(2,-1.5) edge [thick] (2,15.5)
(3,-1.5) edge [thick] (3,15.5)
(4,-1.5) edge [thick] (4,15.5)
(5,-1.5) edge [thick] (5,15.5)
(6,-1.5) edge [thick] (6,15.5)
(7,-1.5) edge [thick] (7,15.5)
(8,-1.5) edge [thick] (8,15.5)
(9,-1.5) edge [thick] (9,15.5)
(10,-1.5) edge[thick]  (10,15.5)
(11,-1.5) edge[thick]  (11,15.5)
(12,-1.5) edge[thick]  (12,15.5)
(13,-1.5) edge[thick]  (13,15.5)
(14,-1.5) edge[thick]  (14,15.5)
(15,-1.5) edge[thick]  (15,15.5)
(-1.5,-1) edge[thick]  (15.5,-1)
(-1.5,0) edge [thick] (15.5,0)
(-1.5,1) edge [thick] (15.5,1)
(-1.5,2) edge [thick] (15.5,2)
(-1.5,3) edge [thick] (15.5,3)
(-1.5,4) edge [thick] (15.5,4)
(-1.5,5) edge [thick] (15.5,5)
(-1.5,6) edge [thick] (15.5,6)
(-1.5,7) edge [thick] (15.5,7)
(-1.5,8) edge [thick] (15.5,8)
(-1.5,9) edge [thick] (15.5,9)
(-1.5,10) edge[thick]  (15.5,10)
(-1.5,11) edge[thick]  (15.5,11)
(-1.5,12) edge[thick]  (15.5,12)
(-1.5,13) edge[thick]  (15.5,13)
(-1.5,14) edge[thick]  (15.5,14)
(-1.5,15) edge[thick]  (15.5,15)
(-1.5, 1) edge [line width = 2.5pt, black] (15.5,1)
(1, -1.5) edge [line width = 2.5pt, black] (1, 15.5)
;
\draw[fill = ForestGreen, opacity=0.08] (-1,-1) rectangle (15,15);
\draw[ultra thick,ForestGreen] (-1,-1) rectangle (15,15);
\draw[fill = blue, opacity=0.08] (1,1) rectangle (13,13);
\draw[ultra thick,blue] (1,1) rectangle (13,13);
\draw [DarkOrchid, ultra thick] (0,0) circle (10pt);
\draw [DarkOrchid, ultra thick] (5,-1) circle (10pt);
\draw [DarkOrchid, ultra thick] (13,0) circle (10pt);
\draw [DarkOrchid, ultra thick] (14,5) circle (10pt);
\draw [DarkOrchid, ultra thick] (15,10) circle (10pt);
\draw [DarkOrchid, ultra thick] (3,15) circle (10pt);
\draw [DarkOrchid, ultra thick] (8,14) circle (10pt);
\draw [DarkOrchid, ultra thick] (0,13) circle (10pt);
\draw [DarkOrchid, ultra thick] (-1,8) circle (10pt);
\end{tikzpicture}
\hspace{1cm}
\begin{tikzpicture}[scale = 0.35,
baseline={([yshift=-.5ex]current bounding box.center)},
vertex/.style = {circle, fill, inner sep=1.4pt, outer sep=0pt},
every edge quotes/.style = {auto=left, sloped, font=\scriptsize, inner sep=1pt}
]
\node[vertex] at (1,1) {} ;
\node[vertex] at (13,1) {} ;
\node[vertex] at (8,1) {} ;
\node[vertex] at (3,2) {} ;
\node[vertex] at (11,3) {} ;
\node[vertex] at (6,4) {} ;
\node[vertex] at (1,5) {} ;
\node[vertex] at (13,5) {} ;
\node[vertex] at (9,6) {} ;
\node[vertex] at (4,7) {} ;
\node[vertex] at (12,8) {} ;
\node[vertex] at (7,9) {} ;
\node[vertex] at (2,10) {} ;
\node[vertex] at (13,10) {} ;
\node[vertex] at (10,11) {} ;
\node[vertex] at (5,12) {} ;
\node[vertex] at (1,13) {} ;
\node[vertex] at (13,13) {} ;
\node[vertex] at (3,13) {}; 
\node[vertex]  at (5,1) {}; 
\node[vertex]  at (8,13)  {}; 
\node[vertex]  at (1,8)  {}; 
\node[vertex] at (0,0) {} ;
\node[vertex] at (13,0) {} ;
\node[vertex] at (8,1) {} ;
\node[vertex] at (3,2) {} ;
\node[vertex] at (11,3) {} ;
\node[vertex] at (6,4) {} ;
\node[vertex] at (1,5) {} ;
\node[vertex] at (14,5) {} ;
\node[vertex] at (9,6) {} ;
\node[vertex] at (4,7) {} ;
\node[vertex] at (12,8) {} ;
\node[vertex] at (7,9) {} ;
\node[vertex] at (2,10) {} ;
\node[vertex] at (15,10) {} ;
\node[vertex] at (10,11) {} ;
\node[vertex] at (5,12) {} ;
\node[vertex] at (0,13) {} ;
\node[vertex] at (13,13) {} ;
\node[vertex] at (3,15) {}; 
\node[vertex]  at (5,-1) {}; 
\node[vertex]  at (8,14)  {}; 
\node[vertex]  at (-1,8)  {}; 

\draw [BurntOrange, ultra thick] (1,1) circle (10pt);
\draw [BurntOrange, ultra thick] (5,1) circle (10pt);
\draw [BurntOrange, ultra thick] (13,1) circle (10pt);

\draw [BurntOrange, ultra thick] (13,5) circle (10pt);
\draw [BurntOrange, ultra thick] (13,10) circle (10pt);
\draw [BurntOrange, ultra thick] (3,13) circle (10pt);

\draw [BurntOrange, ultra thick] (8,13) circle (10pt);
\draw [BurntOrange, ultra thick] (1,13) circle (10pt);
\draw [BurntOrange, ultra thick] (1,8) circle (10pt);

\path
    \foreach \x / \y in {1/1, 5/1, 13/1, 13/5, 13/10, 3/13, 8/13,1/13,1/8}{
            (\x,\y+1) edge [red, ultra thick, dotted] (\x+1,\y)
            (\x+1,\y) edge [red, ultra thick, dotted] (\x,\y-1)
            (\x,\y-1) edge [red, ultra thick, dotted] (\x-1,\y)
            (\x-1,\y) edge [red, ultra thick, dotted] (\x,\y+1)};
\path
(1.5,-1.5) edge [red, ultra thick, dotted] (2,-1)
(2,-1) edge [red, ultra thick, dotted] (2.5,-1.5)
(8.5,-1.5) edge [red, ultra thick, dotted] (10,0)
(10,0) edge [red, ultra thick, dotted] (11.5,-1.5)
(14.5,-1.5) edge [red, ultra thick, dotted] (15,-1)
(15,-1) edge [red, ultra thick, dotted] (15.5,-1.5)
(11,14) edge [red, ultra thick, dotted] (9.5,15.5)
(11,14) edge [red, ultra thick, dotted] (12.5,15.5)
(6,15) edge [red, ultra thick, dotted] (6.5,15.5)
(6,15) edge [red, ultra thick, dotted] (5.5,15.5)
(-1.5,1.5) edge [red, ultra thick, dotted] (0,3)
(0,3) edge [red, ultra thick, dotted] (-1.5,4.5)
(-1.5,10.5) edge [red, ultra thick, dotted] (-1,11)
(-1,11) edge [red, ultra thick, dotted] (-1.5,11.5)
(14,15) edge [red, ultra thick, dotted] (15.5,13.5)
(14,15) edge [red, ultra thick, dotted] (14.5,15.5)
(-0.5,15.5) edge [red, ultra thick, dotted] (-1.5,14.5)
(8,3) edge [red, ultra thick, dotted] (10,1)
(10,1) edge [red, ultra thick, dotted] (8,-1)
(8,-1) edge [red, ultra thick, dotted] (6,1)
(6,1) edge [red, ultra thick, dotted] (8,3) 
(3,4) edge [red, ultra thick, dotted] (5,2)
(5,2) edge [red, ultra thick, dotted] (3,0)
(3,0) edge [red, ultra thick, dotted] (1,2)
(1,2) edge [red, ultra thick, dotted] (3,4) 
(15.5,0.5) edge [red, ultra thick, dotted] (14,2)
(14,2) edge [red, ultra thick, dotted] (15.5,3.5) 
(11,5) edge [red, ultra thick, dotted] (13,3)
(13,3) edge [red, ultra thick, dotted] (11,1)
(11,1) edge [red, ultra thick, dotted] (9,3)
(9,3) edge [red, ultra thick, dotted] (11,5) 
(6,6) edge [red, ultra thick, dotted] (8,4)
(8,4) edge [red, ultra thick, dotted] (6,2)
(6,2) edge [red, ultra thick, dotted] (4,4)
(4,4) edge [red, ultra thick, dotted] (6,6) 
(1,7) edge [red, ultra thick, dotted] (3,5)
(3,5) edge [red, ultra thick, dotted] (1,3)
(1,3) edge [red, ultra thick, dotted] (-1,5)
(-1,5) edge [red, ultra thick, dotted] (1,7) 
(9,8) edge [red, ultra thick, dotted] (11,6)
(11,6) edge [red, ultra thick, dotted] (9,4)
(9,4) edge [red, ultra thick, dotted] (7,6)
(7,6) edge [red, ultra thick, dotted] (9,8) 
(4,9) edge [red, ultra thick, dotted] (6,7)
(6,7) edge [red, ultra thick, dotted] (4,5)
(4,5) edge [red, ultra thick, dotted] (2,7)
(2,7) edge [red, ultra thick, dotted] (4,9) 
(15.5,6.5) edge [red, ultra thick, dotted] (15,7)
(15,7) edge [red, ultra thick, dotted] (15.5,7.5) 
(12,10) edge [red, ultra thick, dotted] (14,8)
(14,8) edge [red, ultra thick, dotted] (12,6)
(12,6) edge [red, ultra thick, dotted] (10,8)
(10,8) edge [red, ultra thick, dotted] (12,10) 
(7,11) edge [red, ultra thick, dotted] (9,9)
(9,9) edge [red, ultra thick, dotted] (7,7)
(7,7) edge [red, ultra thick, dotted] (5,9)
(5,9) edge [red, ultra thick, dotted] (7,11) 
(2,12) edge [red, ultra thick, dotted] (4,10)
(4,10) edge [red, ultra thick, dotted] (2,8)
(2,8) edge [red, ultra thick, dotted] (0,10)
(0,10) edge [red, ultra thick, dotted] (2,12) 
(10,13) edge [red, ultra thick, dotted] (12,11)
(12,11) edge [red, ultra thick, dotted] (10,9)
(10,9) edge [red, ultra thick, dotted] (8,11)
(8,11) edge [red, ultra thick, dotted] (10,13) 
(5,14) edge [red, ultra thick, dotted] (7,12)
(7,12) edge [red, ultra thick, dotted] (5,10)
(5,10) edge [red, ultra thick, dotted] (3,12)
(3,12) edge [red, ultra thick, dotted] (5,14) 
(13,15) edge [red, ultra thick, dotted] (15,13)
(15,13) edge [red, ultra thick, dotted] (13,11)
(13,11) edge [red, ultra thick, dotted] (11,13)
(11,13) edge [red, ultra thick, dotted] (13,15)
(1,1) edge [thick] (1,13)
(2,1) edge [thick] (2,13)
(3,1) edge [thick] (3,13)
(4,1) edge [thick] (4,13)
(5,1) edge [thick] (5,13)
(6,1) edge [thick] (6,13)
(7,1) edge [thick] (7,13)
(8,1) edge [thick] (8,13)
(9,1) edge [thick] (9,13)
(10,1) edge[thick]  (10,13)
(11,1) edge[thick]  (11,13)
(12,1) edge[thick]  (12,13)
(13,1) edge[thick]  (13,13)
(1,1) edge [thick] (13,1)
(1,2) edge [thick] (13,2)
(1,3) edge [thick] (13,3)
(1,4) edge [thick] (13,4)
(1,5) edge [thick] (13,5)
(1,6) edge [thick] (13,6)
(1,7) edge [thick] (13,7)
(1,8) edge [thick] (13,8)
(1,9) edge [thick] (13,9)
(1,10) edge[thick]  (13,10)
(1,11) edge[thick]  (13,11)
(1,12) edge[thick]  (13,12)
(1,13) edge[thick]  (13,13)
;
\draw [line width = 0.4mm, ->] (0,0) to (0.75,0.7);
\draw [line width = 0.4mm, ->] (5,-1) to (5,0.6);
\draw [line width = 0.4mm, ->] (13,0) to (13,0.6);
\draw [line width = 0.4mm, ->] (15,10) to (13.4,10);
\draw [line width = 0.4mm, ->] (8,14) to (8,13.4);
\draw [line width = 0.4mm, ->] (0,13) to (0.6,13);
\draw [line width = 0.4mm, ->] (-1,8) to (0.6,8);
\draw [line width = 0.4mm, ->] (3,15) to (3,13.4);    
\draw [line width = 0.4mm, ->] (14,5) to (13.4,5);
\end{tikzpicture}
\caption{(Left) $G_{13,13}, Y_{17,17}$ and $2$-limited dominating broadcast on $ \mathbb{Z}^2$ defined by  $\phi^{-1} \left( 8 \right) $. (Right) Resulting $2$-limited dominating broadcast on $G_{13,13}$ given $\phi^{-1} \left( 8 \right) $.}
\label{fig:optimalp13}
\end{figure}
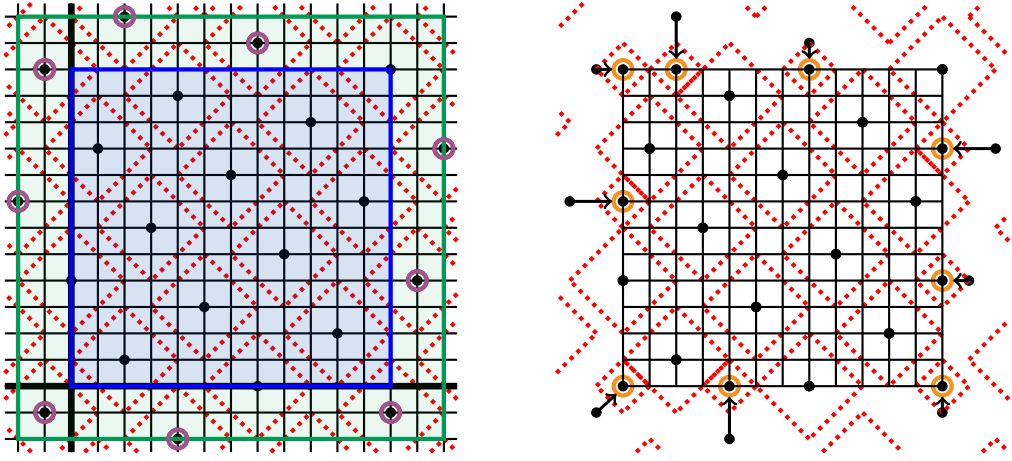   
% !!!
% Colour required for figure
% !!!

Moving each of the broadcast vertices exterior to $G_{13,13}$ and on $Y_{17,17}$ to the nearest point on the border of $G_{13,13}$ and reducing their broadcast strength to one produces a $2$-limited dominating broadcast on $G_{13,13}$. Figure \ref{fig:optimalp13} (Right) depicts the resulting $2$-limited dominating broadcast on $G_{13,13}$ given $ \phi^{-1} \left( 8 \right) $ which establishes $ \gamma_{b,2} \left( P_{13} \square P_{13} \right) \leq 35$. Note that, although this bound is not tight as $ \gamma_{b,2} \left( P_{13} \square P_{13} \right) = 32$ (as determined by computation), it is also not too far from the truth. As $m$ and $n$ get large, we show that the ratio of this upper bound and the lower bound of  $ \gamma_{b,2} \left( P_{m} \square P_{n} \right) $ given by Theorem \ref{thm:mutipack:pathxpathgeneral} approaches one.
\end{eg}

Example \ref{example:pxpinf} illustrates that, given some fixed $\ell \in \mathbb{Z}_{13}$ and $m,n \geq 13$, if there are $x$ and $y$ vertices broadcasting at non-zero strength under $\phi^{-1} \left(\ell  \right) $ on $G_{m,n}$ and $Y_{m+4,n+4}$, respectively, then there is a $2$-limited dominating broadcast on $P_m \square P_n$ of cost $2x + \left( y-x \right) $. Our general construction, as described in the subsequent section, uses the same approach as Example \ref{example:pxpinf} but uses a choice of $\ell \in \mathbb{Z}_{13}$ which gives a best possible construction (in terms of cost).

\subsubsection{General Constructions}\label{subsec:genPxP}
This section proves Theorem \ref{thm:pathxpathgeneralized} which, with Theorem \ref{thm:path}, establishes upper bounds for $\gamma_{b,2} \left( P_m \square P_n \right)  $ for all $m,n \geq 2$. The proof of Theorem \ref{thm:pathxpathgeneralized} uses Lemma \ref{lem:undergrad2} from \cite {undergrad}, stated as follows.
\begin{lem}\label{lem:undergrad2}\rm{\cite[Lemma 3.2]{undergrad}}
Let $\ell \in \mathbb{Z}_{13}$. If either $m$ or $n$ is a multiple of 13, then for any $\ell \in \mathbb{Z}_{13}$ then 
\begin{equation*}
\begin{aligned}
\left| \phi^{-1} \left( \ell \right) \cap V(G_{m,n}) \right|  &=  \frac{ mn }{ 13 } .
\end{aligned}
\end{equation*}
\end{lem}

\begin{thm}\label{thm:pathxpathgeneralized}
If $m,n\geq 13$ then
\begin{equation*}
\begin{aligned}
\gamma_{b,2} \left( P_m \square P_n \right) \leq 2 \left( \frac{ mn }{ 13 }  \right) + 4\left(\frac{   m+n }{ 13 } \right)  + \frac{c \left( m_{13},n_{13} \right)}{13}
\end{aligned}
\end{equation*}
where $ c \left( m_{13},n_{13} \right) $ corresponds with the $O(1)$ terms in Table \ref{table:pathconstant} where $m_{13}$ and $n_{13}$ are the least residues of $m$ and $n$ modulo 13, respectively.
\end{thm}

\begin{table}[htbp]
\centering
\begin{tabular}{cc|ccccccccccccc}
& & \multicolumn{13}{c}{Least residue $n_{13}$ of $n$ modulo $13$} \tabularnewline
& & $0$ & $1$ & $2$ & $3$ & $4$ & $5$ & $6$ &$ 7$ &$ 8$ &$ 9$ & $10$ & $11$ &$ 12$\\ \hline         
\parbox[t]{2mm}{\multirow{13}{*}{\rotatebox[origin=c]{90}{Least residue $m_{13}$ of $m$ modulo $13$}}}
& 0& $0$ & $9$ & $5$ & $1$ & $10$ & $6$ & $2$ & $11$ & $7$ & $16$ & $12$ & $8$ & $4$\tabularnewline 
& 1 & $9$ & $16$ & $10$ & $4$ & $11$ & $5$ & $12$ & $6$ & $0$ & $7$ & $1$ & $8$ & $2$\tabularnewline 
&2 & $5$ & $10$ & $2$ & -$6$ & -$1$ & -$9$ & -$4$ & -$12$ & $6$ & $11$ & $3$ & $8$ & $0$\tabularnewline 
&3 & $1$ & $4$ & -$6$ & $10$ & $0$ & $3$ & $6$ & -$4$ & $12$ & $2$ & $5$ & $8$ & -$2$\tabularnewline 
&4 & $10$ & $11$ & -$1$ & $0$ & -$12$ & $2$ & $3$ & -$9$ & $5$ & $6$ & -$6$ & $8$ & -$4$\tabularnewline 
&5 & $6$ & $5$ & -$9$ & $3$ & $2$ & -$12$ & $0$ & -$1$ & $11$ & $10$ & -$4$ & $8$ & -$6$\tabularnewline 
&6 & $2$ & $12$ & -$4$ & $6$ & $3$ & $0$ & $10$ & -$6$ & $4$ & $1$ & -$2$ & $8$ & $5$\tabularnewline 
&7 & $11$ & $6$ & -$12$ & -$4$ & -$9$ & -$1$ & -$6$ & $2$ & $10$ & $5$ & $0$ & $8$ & $3$\tabularnewline 
&8 & $7$ & $0$ & $6$ & $12$ & $5$ & $11$ & $4$ & $10$ & $16$ & $9$ & $2$ & $8$ & $1$\tabularnewline 
&9 & $16$ & $7$ & $11$ & $2$ & $6$ & $10$ & $1$ & $5$ & $9$ & $0$ & $4$ & $8$ & $12$\tabularnewline 
&10 & $12$ & $1$ & $3$ & $5$ & -$6$ & -$4$ & -$2$ & $0$ & $2$ & $4$ & $6$ & $8$ & $10$\tabularnewline 
&11 & $8$ & $8$ & $8$ & $8$ & $8$ & $8$ & $8$ & $8$ & $8$ & $8$ & $8$ & $8$ & $8$\tabularnewline 
&12 & $4$ & $2$ & $0$ & -$2$ & -$4$ & -$6$ & $5$ & $3$ & $1$ & $12$ & $10$ & $8$ & $6$\tabularnewline 
\end{tabular}
\caption{Value of $c(m_{13},n_{13})$ for the upper bound for $ \gamma_{b,2} \left( P_{m} \square P_{n} \right) $ stated in Theorem \ref{thm:pathxpathgeneralized}.}
\label{table:pathconstant}
\end{table}

\begin{pf}
Let $G_{m,n}$, $Y_{m+4,n+4}$, $ \mathbb{Z}^2$, $ \mathbb{Z}_{13}$, and $\phi$ be defined as in Section \ref{subsec:$2$-limitedd}. Fix some $\ell \in \mathbb{Z}^2$. As $ Y_{m+4,n+4}$ is the distance two neighbourhood of $G_{m,n}$, the set of vertices
\begin{equation*}
\begin{aligned}
V(Y_{m+4,n+4}) \cap \phi^{-1} \left( \ell \right),	
\end{aligned}
\end{equation*}
if broadcasting at strength two, dominate the vertices of $G_{m,n}$. Let 
\begin{center}
\begin{tabular}{ccc}
$B_2 = V(G_{m,n}) \cap \phi^{-1} \left( \ell \right)$ & and & $ B_1 = \phi^{-1} \left( \ell \right) \cap \left( V(Y_{m+4,n+4}) \setminus V(G_{m,n})\right) $.
\end{tabular}
\end{center}
For each vertex  $v \in B_1$ whose broadcast is heard by a vertex in $G_{m,n}$, let $v'$ be the vertex of $G_{m,n}$ with minimum distance to $v$. Note that this choice of $v'$ is unique and that $v'$ necessarily hears a broadcast from $v$. The vertices dominated in $G_{m,n}$ by a broadcast of strength two at $v$ are a subset of the vertices dominated in $G_{m,n}$ by a broadcast of strength one at $v'$. Define the set
\begin{equation*}
\begin{aligned}
B_1' &= \left\{
v' \in V(G_{m,n})\biggm| \begin{array}{c}
v' \textnormal{ is the vertex of } G_{m,n} \textnormal{ which hears a broadcast from}  \\
\textnormal{ and is at minimum distance to } v \textnormal{ for some }  v \in B_1
\end{array}
\right\}. \\ 	
\end{aligned}
\end{equation*}
Informally, $B_1'$ is the resulting collection of vertices formed by moving each vertex in $B_1$ (whose broadcast is heard by a vertex in $G_{m,n}$) to the nearest vertex in $G_{m,n}$. Each vertex in the plane hears only one broadcast under $\phi^{-1} \left(\ell\right) $ \cite[Lemma V.7]{undergrad2}, hence each $v' \in B_1'$ hears only the broadcast from some $v \in B_1$.  Therefore, $v' \not \in B_2$ and $B_2 \cap B_1' = \emptyset$. The broadcast $f: V \left( G_{m,n} \right) \mapsto \left\{ 0,1,2 \right\} $ defined by 
\begin{equation*}
\begin{aligned}
f \left( v \right)  &= \begin{cases}
2 & \textnormal{if } v\in B_2, \\
1 & \textnormal{if } v \in B_1', \textnormal{ and} \\
0 & \textnormal{otherwise} 
\end{cases}
\\ 	
\end{aligned}
\end{equation*}
is a $2$-limited dominating broadcast of $G_{m,n}$. The cost of $f$ is
\begin{equation*}
\begin{aligned}
2 \left| B_2 \right| + \left| B_1' \right| \leq \left|  V(G_{m,n}) \cap \phi^{-1} \left( \ell \right)  \right| + \left|  V(Y_{m+4, n+4})\cap \phi^{-1} \left( \ell \right) \right| . 
\end{aligned}
\end{equation*}
Note that equality holds in the above expression if and only if the broadcast from each vertex $v\in B_1$ is heard by some vertex $v' \in V(G_{m,n})$. Given this construction, it follows that
\begin{equation*}
\begin{aligned}
\gamma_{b,2} \left( P_m \square P_n \right) &\leq \min_{\ell \in \mathbb{Z}_{13}} \left\{ \left|  V(G_{m,n}) \cap \phi^{-1} \left( \ell \right) \right| + \left|  V(Y_{m+4, n+4}) \cap \phi^{-1} \left( \ell \right)  \right| \right\}. 
\end{aligned}
\end{equation*}

For a fixed $\ell \in \mathbb{Z}^2$, $ \left| V(G_{m,n}) \cap \phi^{-1} \left( \ell \right)\right|$ is computed as follows. Let $m_{13}$ and $n_{13}$ be the least residue of $m$ and $n$ modulo 13, respectively. Partition $V\left( G_{m,n}\right)$ into the following subsets:
\begin{center}
\begin{tabular}{l}
$G_1 =  \left\{ (i,j) \in \mathbb{Z}^2 : 0 \leq i \leq n-1-n_{13} \textnormal{ and } 0 \leq j \leq m-1  \right\}$, \\
$G_2 =  \left\{ (i,j) \in \mathbb{Z}^2 : n-n_{13} \leq i \leq n-1 \textnormal{ and } 0 \leq j \leq m-1 - m_{13}   \right\} $, and \\            $G_3 =  \left\{ (i,j) \in \mathbb{Z}^2 :  n-n_{13} \leq i \leq n-1 \textnormal{ and } m-m_{13} \leq j \leq m-1  \right\} $. \\
\end{tabular}
\end{center}
Figure \ref{fig:partitionGmn} (Left) depicts $G_{16,18}$ partitioned into $G_1,G_2,$ and $G_3$.

% !!!
% Colour required for figure
% !!!
\begin{figure}[htbp]
\centering
\begin{tikzpicture}[scale = 0.29,
baseline={([yshift=-.5ex]current bounding box.center)},
vertex/.style = {circle, fill, inner sep=1.4pt, outer sep=0pt},
every edge quotes/.style = {auto=left, sloped, font=\scriptsize, inner sep=1pt}
]

\path
(-1,-1.5) edge[thick]  (-1,16.5)
(0,-1.5) edge [thick] (0,16.5)
(1,-1.5) edge [thick] (1,16.5)
(2,-1.5) edge [thick] (2,16.5)
(3,-1.5) edge [thick] (3,16.5)
(4,-1.5) edge [thick] (4,16.5)
(5,-1.5) edge [thick] (5,16.5)
(6,-1.5) edge [thick] (6,16.5)
(7,-1.5) edge [thick] (7,16.5)
(8,-1.5) edge [thick] (8,16.5)
(9,-1.5) edge [thick] (9,16.5)
(10,-1.5) edge[thick]  (10,16.5)
(11,-1.5) edge[thick]  (11,16.5)
(12,-1.5) edge[thick]  (12,16.5)
(13,-1.5) edge[thick]  (13,16.5)
(14,-1.5) edge[thick]  (14,16.5)
(15,-1.5) edge[thick]  (15,16.5)
(16,-1.5) edge[thick]  (16,16.5)
(17,-1.5) edge[thick]  (17,16.5)
(18,-1.5) edge[thick]  (18,16.5)
(-1.5,-1) edge[thick]  (18.5,-1)
(-1.5,0) edge [thick] (18.5,0)
(-1.5,1) edge [thick] (18.5,1)
(-1.5,2) edge [thick] (18.5,2)
(-1.5,3) edge [thick] (18.5,3)
(-1.5,4) edge [thick] (18.5,4)
(-1.5,5) edge [thick] (18.5,5)
(-1.5,6) edge [thick] (18.5,6)
(-1.5,7) edge [thick] (18.5,7)
(-1.5,8) edge [thick] (18.5,8)
(-1.5,9) edge [thick] (18.5,9)
(-1.5,10) edge[thick]  (18.5,10)
(-1.5,11) edge[thick]  (18.5,11)
(-1.5,12) edge[thick]  (18.5,12)
(-1.5,13) edge[thick]  (18.5,13)
(-1.5,14) edge[thick]  (18.5,14)
(-1.5,15) edge[thick]  (18.5,15)
(-1.5,16) edge[thick]  (18.5,16)

(-1.5, 0) edge [ ultra thick, black] (18.5,0)
(0, -1.5) edge [ultra thick, black] (0, 16.5)
;
\draw[ultra thick, fill =white, opacity =0.8] (0,0) rectangle (12,15);
\draw[ultra thick, fill =white, , opacity =0.8] (13,0) rectangle (17,12);
\draw[ultra thick, fill =white, , opacity =0.8] (13,13) rectangle (17,15);
\node (G_1) at (6,7) [label = : $G_1$] {}; 
\node (G_2) at (15,5) [label = : $G_2$] {}; 
\node (G_4) at (15,12.6) [label = : $G_3$] {}; 
\end{tikzpicture}
\hspace{1cm}
\begin{tikzpicture}[scale = 0.29,
baseline={([yshift=-.5ex]current bounding box.center)},
vertex/.style = {circle, fill, inner sep=1pt, outer sep=0pt},
every edge quotes/.style = {auto=left, sloped, font=\scriptsize, inner sep=1pt}
]
\path

(0,-1) edge  (0,20)
(1,-1) edge  (1,20)
(2,-1) edge  (2,20)
(3,-1) edge  (3,20)
(4,-1) edge  (4,20)
(5,-1) edge  (5,20)
(6,-1) edge  (6,20)
(7,-1) edge  (7,20)
(8,-1) edge  (8,20)
(9,-1) edge  (9,20)
(10,-1) edge  (10,20)
(11,-1) edge  (11,20)
(12,-1) edge  (12,20)
(13,-1) edge  (13,20)
(14,-1) edge  (14,20)
(15,-1) edge  (15,20)
(16,-1) edge  (16,20)
(17,-1) edge  (17,20)
(18,-1) edge  (18,20)
(19,-1) edge  (19,20)
(20,-1) edge  (20,20)
(21,-1) edge  (21,20)
(-1,0) edge  (22,0)
(-1,1) edge  (22,1)
(-1,2) edge  (22,2)
(-1,3) edge  (22,3)
(-1,4) edge  (22,4)
(-1,5) edge  (22,5)
(-1,6) edge  (22,6)
(-1,7) edge  (22,7)
(-1,8) edge  (22,8)
(-1,9) edge  (22,9)
(-1,10) edge  (22,10)
(-1,11) edge  (22,11)
(-1,12) edge  (22,12)
(-1,13) edge  (22,13)
(-1,14) edge  (22,14)
(-1,15) edge  (22,15)
(-1,16) edge  (22,16)
(-1,17) edge  (22,17)
(-1,18) edge  (22,18)
(-1,19) edge  (22,19)
(2, -1) edge [ultra thick]  (2, 20)
(-1, 2) edge [ultra thick] (22,2)
;
\draw[ultra thick,fill = Gray] (2,2) rectangle (14,17);
\draw[ultra thick,fill =Gray ] (15,2) rectangle (19,14);
\draw[ultra thick,fill =Gray ] (15,15) rectangle (19,17);

\draw[ultra thick,fill = white, opacity=0.8] (0,0) rectangle (12,19);
\draw[ultra thick,fill = white, opacity=0.8] (13,0) rectangle (21,12);
\draw[ultra thick,fill = white, opacity=0.8] (13,13) rectangle (21,19);
\node (Y_1) at (6,9) [label = : $Y_1$] {}; 
\node (Y_2) at (17,5) [label = : $Y_2$] {}; 
\node (Y_4) at (17,14.5) [label = : $Y_3$] {}; 

\end{tikzpicture}
\caption{(Left) $G_{16,18}$ partitioned into $G_1,G_2,$ and $G_3$. (Right) $G_{16,18}$ partitioned into $G_1,G_2,$ and $G_3$  overlaid by $Y_{20,22}$ partitioned into $Y_1,Y_2,$ and $Y_3$.}
\label{fig:partitionGmn}
\end{figure}
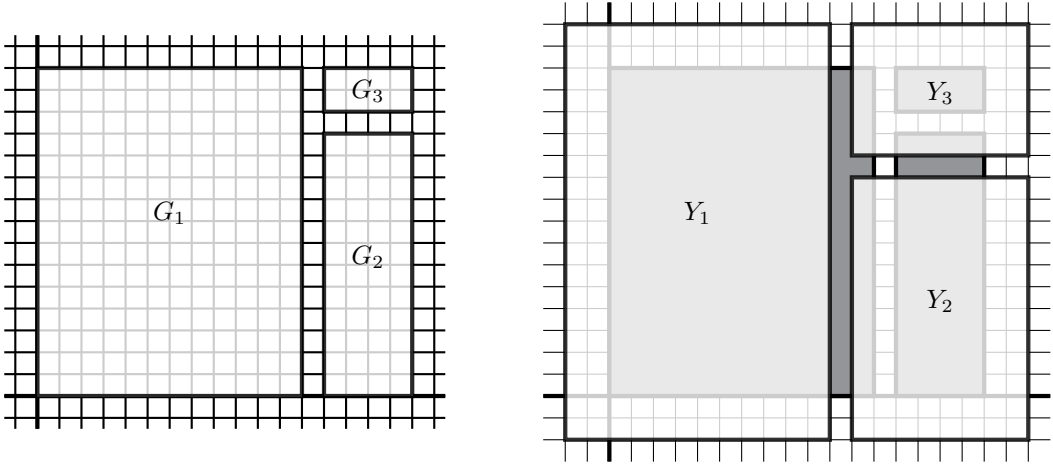   
% !!!
% Colour required for figure
% !!!

For each $\ell \in \mathbb{Z}_{13}$,
\begin{equation*}
\begin{aligned}
\left| V(G_{m,n}) \cap \phi^{-1} \left( \ell \right)  \right|  &= \left| G_1 \cap \phi^{-1} \left( \ell \right)  \right| +\left| G_2 \cap \phi^{-1} \left( \ell \right)  \right| +\left| G_3 \cap \phi^{-1} \left( \ell \right)  \right| .\\ 	
\end{aligned}
\end{equation*}
By Lemma \ref{lem:undergrad2},
\begin{equation}\label{eq:G4}
\begin{aligned}
\left| G_1 \cap \phi^{-1} \left( \ell \right)\right| = \frac{  m  \left( n-n_{13} \right)  }{ 13 } \hspace{0.5cm} \textnormal{and} \hspace{0.5cm}  \left| G_2 \cap \phi^{-1} \left( \ell \right)\right| = \frac{ n_{13}  \left( m-m_{13} \right)  }{ 13 }   \\  
\end{aligned}
\end{equation}
regardless of $\ell \in \mathbb{Z}_{13}$. Hence, to determine $ \left| V(G_{m,n}) \cap \phi^{-1} \left( \ell \right)  \right| $, what remains is to find $ \left| G_3 \cap \phi^{-1} \left( \ell \right)\right| $. The same methodology can be applied to $Y_{m+4,n+4}$ by defining $m_{13}^Y$ and $n_{13}^Y$ as the least residues of $(m_{13}+4)$ and $(n_{13} +4)$ modulo 13, respectively, and partitioning $V \left( Y_{m+4,n+4}\right)$ into the following subsets:
\begin{center}
\begin{tabular}{l}
$Y_1 =  \left\{ (i,j) \in \mathbb{Z}^2 : -2 \leq i \leq n+1-n_{13}^Y \textnormal{ and }  -2 \leq j \leq m+1 \right\}$, \\
$Y_2 =  \left\{ (i,j) \in \mathbb{Z}^2 :  n+2-n_{13}^Y\leq i \leq n+1 \textnormal{ and }  -2 \leq j \leq m+1-m_{13}^Y\right\} $, and  \\            $Y_3 =  \left\{ (i,j) \in \mathbb{Z}^2 :  n+2-n_{13}^Y \leq i \leq n+1  \textnormal{ and }  m+2-m_{13}^Y \leq j \leq m+1\right\} $. \\
\end{tabular}
\end{center}
Figure \ref{fig:partitionGmn} (Right) depicts $G_{16,18}$ partitioned into $G_1,G_2,$ and $G_3,$ (all in grey) overlaid by $Y_{20,22}$ partitioned into $Y_1,Y_2,$ and $Y_3$ (all opaque white). By Lemma \ref{lem:undergrad2},
\begin{equation}\label{eq:Y4}
    \begin{aligned}
        \left| Y_1 \cap \phi^{-1} \left( \ell \right)\right| = \frac{ \left( m+4 \right) \left( n+4-n_{13}^Y \right)  }{ 13 }  \hspace{0.5cm} \textnormal{and} \hspace{0.5cm} \left| Y_2 \cap \phi^{-1} \left( \ell \right)\right| = \frac{  n_{13}^Y  \left( m+4-m_{13}^Y \right)  }{ 13 }  	
    \end{aligned}
\end{equation}
regardless of $\ell \in \mathbb{Z}_{13}$. Finding a best choice of $\ell$ for the construction therefore reduces to determining, for each least residue of $m$ and $n$ modulo $13$, the minimum values of $ \left| G_3 \cap \phi^{-1} \left( \ell \right)  \right|  + \left| Y_3 \cap \phi^{-1} \left( \ell \right)\right|$ for $\ell \in \mathbb{Z}_{13}$. Said least values are given in Table \ref{table:pathcorner} and the corresponding $\ell \in \mathbb{Z}_{13}$ is given in Table \ref{table:pathpositin} (see Section \ref{appendix}). In the cases where there are several such values of $\ell \in \mathbb{Z}_{13}$ for a given value in  Table \ref{table:pathcorner}, we state the lexicographically largest in Table \ref{table:pathpositin} (see Section \ref{appendix}).

\begin{table}[htbp]
\centering
\begin{tabular}{cc|ccccccccccccc}
& & \multicolumn{13}{c}{Least residue $n_{13}$ of $n$ modulo $13$} \tabularnewline
& & $0$ & $1$ & $2$ & $3$ & $4$ & $5$ & $6$ &$ 7$ &$ 8$ &$ 9$ & $10$ & $11$ &$ 12$\\ \hline
\parbox[t]{2mm}{\multirow{13}{*}{\rotatebox[origin=c]{90}{Least residue $m_{13}$ of $m$ modulo $13$}}} &$0$ &  $0$ & $1$ & $1$ & $1$ & $2$ & $2$ & $2$ & $3$ & $3$ & $0$ & $0$ & $0$ & $0$ \tabularnewline
&$1$ & $1$ & $2$ & $2$ & $2$ & $3$ & $3$ & $4$ & $4$ & $4$ & $0$ & $0$ & $1$ & $1$ \tabularnewline
&$2$ & $1$ &$2$ & $2$ & $2$ & $3$ & $3$ & $4$ & $4$ & $6$ & $1$ & $1$ & $2$ & $2$ \tabularnewline
&$3$ & $1$ &$2$ &$2$ & $4$ & $4$ & $5$ & $6$ & $6$ & $8$ & $1$ & $2$ & $3$ & $3$ \tabularnewline
&$4$ & $2$ &$3$ &$3$ &$4$ & $4$ & $6$ & $7$ & $7$ & $9$ & $2$ & $2$ & $4$ & $4$ \tabularnewline
&$5$ & $2$ &$3$ &$3$ &$5$ &$6$ & $6$ & $8$ & $9$ & $11$ & $3$ & $3$ & $5$ & $5$ \tabularnewline
&$6$ & $2$ &$4$ &$4$ &$6$ &$7$ &$8$ & $10$ & $10$ & $12$ & $3$ & $4$ & $6$ & $7$ \tabularnewline
&$7$ & $3$ &$4$ &$4$ &$6$ &$7$ &$9$ &$10$ & $12$ & $14$ & $4$ & $5$ & $7$ & $8$ \tabularnewline
&$8$ & $3$ &$4$ &$6$ &$8$ &$9$ &$11$ &$12$ &$14$ & $16$ & $5$ & $6$ & $8$ & $9$ \tabularnewline
&$9$ & $0$ &$0$ &$1$ &$1$ &$2$ &$3$ &$3$ &$4$ &$5$ & $5$ & $6$ & $7$ & $8$ \tabularnewline
&$10$ & $0$ &$0$ &$1$ &$2$ &$2$ &$3$ &$4$ &$5$ &$6$ &$6$ & $7$ & $8$ & $9$ \tabularnewline
&$11$ & $0$ &$1$ &$2$ &$3$ &$4$ &$5$ &$6$ &$7$ &$8$ &$7$ &$8$ & $9$ & $10$ \tabularnewline
&$12$ & $0$ &$1$ &$2$ &$3$ &$4$ &$5$ &$7$ &$8$ &$9$ &$8$ &$9$ &$10$ & $11$ \tabularnewline
\end{tabular}
\caption{Minimum values of $ \left| G_3 \cap \phi^{-1} \left( \ell \right)  \right|  + \left| Y_3 \cap \phi^{-1} \left( \ell \right)\right|$ for $\ell \in \mathbb{Z}_{13}$ for $m,n \geq 13$.}
\label{table:pathcorner}
\end{table}

The values in Table \ref{table:pathcorner} and Table \ref{table:pathpositin} (see Section \ref{appendix}) can be verified by an exhaustive search. Let $ c'\left( m_{13},n_{13} \right) $ correspond with the values given in Table \ref{table:pathcorner}. For a general $ \ell \in \mathbb{Z}_{13}$, by Equations \ref{eq:G4} and \ref{eq:Y4}, $\left| \phi^{-1} \left( \ell \right)  \cap V(G_{m,n}) \right| + \left| \phi^{-1} \left( \ell \right)  \cap V(Y_{m+4, n+4}) \right| $ is simply
 \begin{equation*}
     \begin{aligned}
         \frac{  m  \left( n-n_{13} \right)  }{ 13 }+ \frac{ n_{13}\left( m-m_{13} \right)   }{ 13 }+\frac{ \left( m+4 \right) \left( n+4-n_{13}^Y \right)  }{ 13 } + \frac{ n_{13}^Y\left( m+4-m_{13}^Y \right)   }{ 13 }+ c'\left( m_{13},n_{13} \right). \\ 	
     \end{aligned}
 \end{equation*}
As $n_{13}^Y$ and $m_{13}^Y$ are the least residues of $(m_{13}+4)$ and $(n_{13} +4)$ modulo 13, respectively, the previous equation can be simplified for each least residue of $m$ and $n$ modulo $13$ by computation. The result gives
 \begin{equation*}
     \begin{aligned}
         2 \left( \frac{ mn }{ 13 }  \right) + 4 \left( \frac{ m+n }{ 13 }  \right) + \frac{ c \left( m_{13},n_{13} \right) }{ 13}   
     \end{aligned}
 \end{equation*}
where the values $ c \left( m_{13},n_{13} \right) $ are given in Table \ref{table:pathconstant}. This proves the theorem.
\end{pf}

\section{Upper Bounds for $\gamma_{b,2} \left(P_m \square C_n\right)$}\label{sec:pathsandcycle}
This section establishes upper bounds for the $2$-limited broadcast domination number of the Cartesian product of a path and a cycle. Section \ref{sec:P-C3} states bounds for $\gamma_{b,2}(P_{m}\square C_{n})$ for $2 \leq m \leq 22$ and $n \geq 3$. Section \ref{sec:P23-C13} states  bounds for $\gamma_{b,2}(P_{m}\square C_{n})$ for $m \geq 23$ and $n \geq 13$. Section \ref{PCfinal} states  bounds for $\gamma_{b,2}(P_{m}\square C_{n})$ for $m \geq 23$ and $3 \leq n \leq 12$. These results require additional case work, when compared to  Section \ref{sec:paths}, since $ \gamma_{b,2} \left(P_m \square P_n \right) =\gamma_{b,2} \left(P_n \square P_m \right)$ but $ \gamma_{b,2} \left(P_m \square C_n \right)$ is not necessarily equal to $\gamma_{b,2} \left(P_n \square C_m \right)$. All results in this section are established using the tiling method described in Section \ref{sec:upperbound}. 

\subsection{$P_{ 2 \leq m \leq 22}\square C_{n \geq 3}$}\label{sec:P-C3}
The tilings used to establish Theorem \ref{thm:pathcycle} are in \cite[Sections 3.2 through 3.13]{slobodin}.

\begin{thm}\label{thm:pathcycle}
Fix $2 \leq m \leq 22$ and $n \geq 3$. Let $x$ be the value in Table \ref{table:pathxcycle} dependant upon $m$ and define $n_x$ as the least residue of $n$ modulo $x$. By construction, 
\begin{equation*}
\begin{aligned}
\gamma_{b,2} \left(P_m \square C_n\right) \leq b(m) + c(m,n,n_x),
\end{aligned}
\end{equation*}
where $b(m)$ corresponds with the $O(n)$ terms in Table \ref{table:pathxcyclelinearterm} and $c(m,n,n_x)$ corresponds with the $O(1)$ terms in Table \ref{table:pathxcycleconstanterm} (see Section \ref{appendix}).
\end{thm}

\begin{table}[htbp]
\centering
\begin{tabular}{c|ccccccccccccccc}
$m$ &  $2,3$ &$4$  & $5$ & $6$ & $7$ & $8$ & $9$ & $10$ & $11$ & $12$ & $13, \dots,  22$\\ \hline
$x $ & $1$  & $10$ & $2$ & $16$ & $14$ & $22$ & $10$ & $18$ & $26$ & $24$ & $13$ \\
\end{tabular}
\caption{Value of $x$ in the upper bound of $ \gamma_{b,2} \left(P_m \square C_n\right)$ for $2 \leq m \leq 22$ and $n \geq 3$.}
\label{table:pathxcycle}
\end{table}

\begin{table}[htbp]
\renewcommand\arraystretch{1.25}
\begin{tabular}{c|ccccccccccc}
$m$ & $2$ & $3$ & $4$ & $5$ & $6$ & $7$ & $8$ & $9$ & $10$ & $11$ & $12$\\ \hline
$b(m) $&        $ \left\lceil \frac{n}{2}  \right\rceil   $ & $ \left\lceil \frac{2n}{3}  \right\rceil $ & $ 8 \left\lfloor \frac{n}{10}  \right\rfloor $ & $n$ & $ 18 \left\lfloor \frac{n}{16}  \right\rfloor $ & $ 18 \left\lfloor \frac{n}{14}  \right\rfloor $ & $32 \left\lfloor \frac{n}{22}  \right\rfloor  $ & $ 16 \left\lfloor \frac{n}{10}  \right\rfloor  $ & $ 32 \left\lfloor  \frac{n}{18}  \right\rfloor $ & $ 50 \left\lfloor \frac{n}{26}  \right\rfloor  $ & $ 50 \left\lfloor \frac{n}{ 24}  \right\rfloor $ \\
\multicolumn{10}{c}{ }\\
\end{tabular}

\begin{tabular}{c|ccccccccccc}
$m$ & $13$ & $14$ & $15$ & $16$ & $17$ & $18$ & $19$ & $20$ & $21$ & $22$ \\ \hline
$b(m) $& $30 \left\lfloor \frac{n}{13}  \right\rfloor $ &   $32 \left\lfloor \frac{n}{13}  \right\rfloor $ &  $34 \left\lfloor \frac{n}{13}  \right\rfloor $ & $36 \left\lfloor \frac{n}{13}  \right\rfloor $ & $38 \left\lfloor \frac{n}{13}  \right\rfloor $ & $40 \left\lfloor \frac{n}{13}  \right\rfloor $ & $42 \left\lfloor \frac{n}{13}  \right\rfloor $  & $44 \left\lfloor \frac{n}{13}  \right\rfloor $ & $46 \left\lfloor \frac{n}{13}  \right\rfloor $ & $48 \left\lfloor \frac{n}{13}  \right\rfloor $\\
\end{tabular}
\caption{Value of $b(m)$ in the upper bound of $ \gamma_{b,2} \left(P_m \square C_n\right)$ for $2 \leq m \leq 22$ and $n \geq 3$.}
\label{table:pathxcyclelinearterm}
\end{table}
% default
\renewcommand\arraystretch{1}

\subsection{$P_{m \geq 23}\square C_{n \geq 13}$}\label{sec:P23-C13}
The tilings used to establish Theorem \ref{thm:pathxcyclegeneral} are in \cite[Section  3.14]{slobodin}.

\begin{thm}\label{thm:pathxcyclegeneral}
If $m \geq 23$ and $n \geq 13$ then
\begin{equation*}
\begin{aligned}
\gamma_{b,2} \left( P_{m} \square C_{n} \right)  &\leq 2 \left(  \frac{ mn }{ 13 } \right)  + \frac{ 4m }{13  }+  \frac{ b(n)}{ 13 }  + \frac{  c(m_{13}',n_{13}) }{13  }  \\
\end{aligned}
\end{equation*}
where
\begin{equation*}
\begin{aligned}
b(n) &= \begin{cases}
0 & \textnormal{ for }  n \equiv 0 \pmod{ 13 }, \\
2n & \textnormal{ for } n  \equiv  4,7,11, \textnormal{or } 12 \pmod{ 13 }, \\
3n & \textnormal{ for } n \equiv  1,2,5,6,8,9, \textnormal{or } 10\pmod{ 13 }, \textnormal{ and}  \\
4n & \textnormal{ for }  n \equiv 3 \pmod{ 13 }\\
\end{cases}
\\
\end{aligned}
\end{equation*}
and $c(m_{13}',n_{13})$ corresponds with the $O(1)$ terms in Table \ref{table:pathxcyclegeneralcornerforforumla} (see Section \ref{appendix}) where $m_{13}'$ and $n_{13}$ are the least residues of $(m-10)$ and $n$ modulo 13, respectively.
\end{thm}

\subsection{$P_{ m \geq 23} \square C_{3 \leq n \leq 12}$}\label{PCfinal}
The tilings used to establish Theorem \ref{thm:pathcyclefinal} are in \cite[Sections 3.16 through 3.25]{slobodin}.

\begin{thm}\label{thm:pathcyclefinal}%
Fix $m \geq 23$ and $3 \leq n \leq 12$. Let $x$ be the value in Table \ref{table:pathxcyclefinal} dependant upon $n$ and define $m_x$ as the least residue of $m$ modulo $x$. By construction, 
\begin{equation*}
\begin{aligned}
\gamma_{b,2} \left(P_m \square C_n\right) \leq b(n) + c(n,m_x),
\end{aligned}
\end{equation*}
where $b(n)$ corresponds with the $O(m)$ terms in Table \ref{table:pathxcyclefinallineartermfinal} and $c(n,m_x)$ corresponds with the $O(1)$ terms in Table \ref{table:pathxcyclefinalconstanterm} (see Section \ref{appendix}).
\end{thm}

\begin{table}[htbp]
\centering
\begin{tabular}{c|ccccccccccccc}
$n$ & $3$ & $4$ & $5,6$ & $7$ & $8$ & $9,10$ & $11, 12$\\ \hline
$x $ & $1$ & $6$ & $1$  & $35$ & $6$ & $10$  & $13$  \\
\end{tabular}
\caption{Value of $x$ in the upper bound of $ \gamma_{b,2} \left(P_m \square C_n\right)$ for  $m \geq 23$ and $3 \leq n \leq 12$.}
\label{table:pathxcyclefinal}
\end{table}

\renewcommand\arraystretch{1.25}
\begin{table}[htbp]
\centering
\begin{tabular}{c|ccccccccccc}
$n$ & $3$ & $4$ & $5, 6$ & $7$ & $8$ & $9,10$ & $11$ & $12$\\ \hline
$b(n)$ & $ \left\lceil \frac{2m}{3}  \right\rceil $ & $ 4 \left\lfloor \frac{m}{6}  \right\rfloor  $ & $m$  & $ 42 \left\lfloor \frac{m}{35}  \right\rfloor  $ & $ 8 \left\lfloor \frac{m}{6}  \right\rfloor $ & $16 \left\lfloor \frac{m}{10}  \right\rfloor $  & $24 \left\lfloor \frac{m}{13}  \right\rfloor $ & $ 26 \left\lfloor \frac{m}{13}  \right\rfloor $ \\
\end{tabular}
\caption{Value of $b(n)$ in the upper bound of $ \gamma_{b,2} \left(P_m \square C_n\right)$ for $m \geq 23$ and $3 \leq n \leq 12$.}
\label{table:pathxcyclefinallineartermfinal}
\end{table}
% default
\renewcommand\arraystretch{1}

\section{Upper Bounds for $\gamma_{b,2} \left(C_m \square C_n\right)$}\label{sec:cycles}

This section establishes upper bounds for the $2$-limited broadcast domination number of the Cartesian product of two cycles. Section \ref{cycle:part1} states bounds for $\gamma_{b,2}(C_{m}\square C_{n})$ for $3 \leq m \leq 25$ and $n \geq m$. Section \ref{cycle:part2} states  bounds for $\gamma_{b,2}(C_{m}\square C_{n})$ for $m \geq 26$ and $n \geq m$. All results in this section are established using the tiling method described in Section \ref{sec:upperbound}. 

\subsection{$C_{3 \leq m \leq 25} \square C_{n \geq m}$}\label{cycle:part1}
The tilings used to establish Theorem \ref{thm:pathcyclefinal} are in \cite[Sections 4.2 through 4.10]{slobodin}.

\begin{thm}\label{thm:cycle}    
Fix $3 \leq m \leq 25$ and $n \geq m$. Let $x$ be the value in Table \ref{table:cyclexcyclex} dependant upon $m$ and define $n_x$ as the least residue of $n$ modulo $x$. By construction, 
\begin{equation*}
\begin{aligned}
\gamma_{b,2} \left(C_m \square C_n\right) \leq b(m) + c(m,n_x),
\end{aligned}
\end{equation*}
where $b(m)$ corresponds with the $O(n)$ terms in Table \ref{table:cyclexcyclelinearterm} and $c(m,n_x)$ corresponds with the $O(1)$ terms in Table \ref{table:cyclexcycleconstanterm} (see Section \ref{appendix}).
\end{thm}

\begin{table}[htbp]
\centering
\begin{tabular}{c|ccccccccccccc}
$m$ & $3$ & $4$ & $5$ &  $6$ & $7$ & $8$ & $9,10$ & $11, \dots,  25$\\ \hline
$x $ & $1$ & $6$ & $1$ & $4$ & $35$ & $6$ & $10$  & $13$  \\
\end{tabular}
\caption{Value of $x$ in the upper bound of $ \gamma_{b,2} \left(C_m \square C_n\right)$ for $3 \leq m \leq 25$ and $n \geq m$.}
\label{table:cyclexcyclex}
\end{table}

\renewcommand\arraystretch{1.25}
\begin{table}[htbp]
\centering
\begin{tabular}{c|ccccccccccc}
$m$ & $3$ & $4$ & $5, 6$ & $7$ & $8$ & $9,10$ & $11$ & $12, 13$ & $14$ & $15$\\ \hline
$b(m)$ & $ \left\lceil \frac{2n}{3}  \right\rceil $ & $ 4 \left\lfloor \frac{n}{6}  \right\rfloor  $ & $n$  & $ 42 \left\lfloor \frac{n}{35}  \right\rfloor  $ & $ 8 \left\lfloor \frac{n}{6}  \right\rfloor $ & $16 \left\lfloor \frac{n}{10}  \right\rfloor $  & $24 \left\lfloor \frac{n}{13}  \right\rfloor $ & $ 26 \left\lfloor \frac{n}{13}  \right\rfloor $ & $ 31 \left\lfloor \frac{n}{13}  \right\rfloor $ & $33 \left\lfloor \frac{n}{13}  \right\rfloor $\\
\multicolumn{12}{c}{ }\\
$m$ & $16,17$ & $18$ & $19$ & $20$ & $21$ & $22$ & $23$ & $24$ & $25$ \\ \hline
$b(m)$ &  $ 36 \left\lfloor \frac{n}{13}  \right\rfloor $ & $39 \left\lfloor \frac{n}{13}  \right\rfloor $ & $41 \left\lfloor \frac{n}{13}  \right\rfloor $& $42 \left\lfloor \frac{n}{13}  \right\rfloor $& $45 \left\lfloor \frac{n}{13}  \right\rfloor $& $48 \left\lfloor \frac{n}{13}  \right\rfloor $& $49 \left\lfloor \frac{n}{13}  \right\rfloor $& $50 \left\lfloor \frac{n}{13}  \right\rfloor $& $52 \left\lfloor \frac{n}{13}  \right\rfloor $\\
\end{tabular}
\caption{Value of $b(m)$ in the upper bound of $ \gamma_{b,2} \left(C_m \square C_n\right)$ for $3 \leq m \leq 25$ and $n \geq m$.}
\label{table:cyclexcyclelinearterm}
\end{table}
% default
\renewcommand\arraystretch{1}

\subsection{$C_{m \geq 26} \square C_{n \geq m}$}\label{cycle:part2}
The tilings used to establish Theorem \ref{thm:pathcyclefinal} are in \cite[Section 4.11]{slobodin}.

\begin{thm}\label{thm:cyclexcyclegeneral}
If $m,n\geq 26$ then
\begin{equation*}
\begin{aligned}
\gamma_{b,2} \left( C_{m} \square C_{n} \right)  &\leq 2 \left(  \frac{ mn }{ 13 } \right)  + \frac{ b(m) + b(n)}{ 13 }  -c(m_{13},n_{13})  \\ 	
\end{aligned}
\end{equation*}
where
\begin{equation*}
\begin{aligned}
b(x) &= \begin{cases}
0 & \textnormal{ if } 0 \equiv x \pmod{ 13 }, \\
4x & \textnormal{ if } 3 \equiv x \pmod{ 13 }, \\
2x & \textnormal{ if } 4,7,11, \textnormal{or } 12  \equiv x \pmod{ 13 }, \textnormal{ and}  \\
3x & \textnormal{ if } 1,2,5,6,8,9, \textnormal{or } 10  \equiv x \pmod{ 13 } \\
\end{cases}
\\ 	
\end{aligned}
\end{equation*}
and $c(m_{13},n_{13})$ corresponds with the $O(1)$ terms in Table \ref{table:cyclxcyclegeneralcornerfinal} (see Section \ref{appendix}) where $m_{13}$ and $n_{13}$ are the least residues of $m$ and $n$ modulo 13, respectively.
\end{thm}

\section{Multipacking}\label{sec:multipack}

Sections \ref{sec:paths}, \ref{sec:pathsandcycle}, and \ref{sec:cycles} established upper bounds for the $2$-limited broadcast domination numbers of the Cartesian products of two paths, a path and a cycle, and two cycles. This section is devoted to determining corresponding lower bounds for the $2$-limited broadcast domination numbers of these graphs. These bounds are obtained via LP duality by finding the lower bounds for the fractional $2$-limited multipacking numbers  of these graphs. All results, with the exception of those stated in Section \ref{sec:multipack:cycle}, were found by using an exact LP solver \cite{soplexagainGleixnerSteffyWolter2012, soplexGleixnerSteffyWolter2015} as a part of the SoPlex distribution \cite{soplexmainGleixnerEiflerGallyetal2017}. 

This section focuses on \textit{fractional $2$-limited multipackings}, which is defined by setting $k=2$ in \ref{LP:multipackc}. Section \ref{sec:multipack:cycle} establishes optimal values for $mp_{f,2}( C_{m} \square C_{n} )$ for all $m,n \geq 3$. Sections \ref{sec:multipack:pathxcycle} and \ref{sec:multipack:path} establish lower bounds for $mp_{f,2}( P_{m} \square C_{n} )$ and $mp_{f,2}( P_{m} \square P_{n} ) $, respectively, via multipacking constructions.

\subsection{The Cartesian Product of Two Cycles}\label{sec:multipack:cycle}

This section describes a general fractional $2$-limited multipacking on $ C_{m} \square C_{n} $ for $m,n \geq 3$. These results establish Theorem \ref{prop:mulitpackcycle}. The proof of Theorem \ref{prop:mulitpackcycle} utilizes the LP relaxation of $2$-limited broadcast domination. The optimal solution of the LP relaxation of \ref{ILP} for $k=2$ defines the \textit{fractional $2$-limited broadcast domination number} $ \gamma_{f,b,2}(G)$ of a graph $G$.

\begin{thm}\label{prop:mulitpackcycle}
For $m,n\geq 3$,
\begin{equation*}
mp_{f,2}( C_{m} \square C_{n} ) = \begin{cases}
\frac{2 mn }{ 13 }  & \textnormal{ if } m, n \geq 5,\\
\frac{4n}{6} & \textnormal{ if } m = 4 \textnormal{ and } n \geq 5,\\
\frac{6n}{11} & \textnormal{ if } m = 3 \textnormal{ and } n \geq 5,\\
\frac{32}{11}  & \textnormal{ if } m = n = 4,\\
\frac{12}{5}  & \textnormal{ if } m = 4 \textnormal{ and }  n = 3, \textnormal{ and} \\
\frac{9}{5}  & \textnormal{ if } m = n = 3.
\end{cases}
\end{equation*}
\end{thm}

\begin{pf}
Fix $m,n \geq 5$. For each vertex $i \in V(G)$, define the variable $y_i$. Define the fractional $2$-limited multipacking on $ C_{m} \square C_{n} $ by $y_i = 2/13 $ for all $i \in V( C_{m} \square C_{n} )$.  As $ C_{m} \square C_{n} $ is vertex-transitive, for each vertex $i \in V ( C_{m} \square C_{n} )$, there are five and thirteen vertices within distance one and two, respectively, of $i$.  This assignment therefore satisfies the constraints of \ref{LP:multipackc} with $k=2$. The cost of this multipacking is $2mn/13 $ which establishes $2mn/13  \leq mp_{f,2}( C_{m} \square C_{n} )$. To prove the opposite inequality, consider the following fractional $2$-limited dominating broadcast. For $i \in V \left( C_{m} \square C_{n}  \right) $ and $k \in \left\{ 1,2 \right\} $ define $x_{i,1}$ and $x_{i,2}$ as the fractional broadcasts of vertex $i$ at strengths one and two, respectively. Define the fractional $2$-limited broadcast by $x_{i,1} =0$ and $x_{i,2} =1/13 $ for all $ i \in V ( C_{m} \square C_{n} )$. As $m,n \geq 5$, for each $i \in V ( C_{m} \square C_{n} )$, there are thirteen vertices within distance two of $i$. This assignment therefore satisfies the constraints of the LP relaxation of \ref{ILP} for $k=2$. The cost of this broadcast is $2mn/13 $ which establishes $\gamma_{f,b,2} \left( C_{m} \square C_{n} \right) \leq 2mn/13$. By the duality theorem of linear programming, $ mp_{f,2} ( C_{m} \square C_{n} ) = \gamma_{f,b,2} \left( C_{m} \square C_{n}  \right) $.

For $m=4$ and $n \geq 5$, define the fractional $2$-limited multipacking on $ C_{4} \square C_{n} $ by $y_i = 1/6 $ for all $i \in V( C_{4} \square C_{n} )$ and the fractional $2$-limited dominating broadcast by $x_{i,1} =0$ and $x_{i,2} = 1/12 $ for all $ i \in V ( C_{4} \square C_{n} )$. As $n \geq 5$, for each $i \in V ( C_{4} \square C_{n} )$, there are five and twelve vertices within distance one and two, respectively, of $i$. The cost of this multipacking and this broadcast are both equal to $ 4n/6 $.

For $m=3$ and $n \geq 5$, define the fractional $2$-limited multipacking on $ C_{3} \square C_{n} $ by $y_i = 2/11 $ for all $i \in V( C_{3} \square C_{n} )$ and the fractional $2$-limited dominating broadcast by $x_{i,1} =0$ and $x_{i,2} = 1/11 $ for all $ i \in V ( C_{3} \square C_{n} )$. For each $i \in V ( C_{3} \square C_{n} )$, there are five and eleven vertices within distance one and two, respectively, of $i$. The cost of this multipacking and this broadcast are both equal to $ 6n/11 $.

For $m =n=4$, use the same fractional $2$-limited multipacking and fractional $2$-limited dominating broadcast as defined for $m=3$ and $n \geq 5$. The cost of this multipacking and this broadcast are both equal to $ 32/11 $.

For $m = 4$ and $n = 3$, define the fractional $2$-limited multipacking on $ C_{4} \square C_{3} $ by $y_i = 1/5 $ for all $i \in V( C_{4} \square C_{3} )$ and the fractional $2$-limited dominating broadcast by $x_{i,1} =1/5 $ and $x_{i,2} = 0 $ for all $ i \in V ( C_{4} \square C_{3} )$. For each $i \in V ( C_{4} \square C_{3} )$, there are five and ten vertices within distance one and two, respectively, of $i$.  The cost of this multipacking and this broadcast are both equal to $ 12/5$.

For $m =n=3$, use the same fractional $2$-limited multipacking and fractional $2$-limited dominating broadcast as defined for $m=4$ and $n =3$. The cost of this multipacking and this broadcast are both equal to $9/5$.
\end{pf}

Combined, the results found in Section \ref{sec:cycles} and Theorem \ref{prop:mulitpackcycle} establish the upper and lower bounds, respectively, in Table \ref{table:multipackcycle} (excluding $m=n=3$, $m=n=4$, and $m=4$ and $n=3$, the values of which are given in Theorem \ref{prop:mulitpackcycle}) for the $2$-limited broadcast domination number of  Cartesian products of two  cycles. The values of the constant terms in the upper bounds in Table \ref{table:multipackcycle} can be found in their respective theorems in Section \ref{sec:cycles} and are omitted here. To easily compare the upper and lower bounds in Table \ref{table:multipackcycle} we also, for some bounds, include a simpler lower bound. Observe that the bounds in Table \ref{table:multipackcycle} give optimal values for $\gamma_{b,2} \left( C_{m} \square C_{n} \right)$ for all $m,n \geq 3$ such that $m,n \equiv 0 \pmod{ 13 } $ and $\gamma_{b,2}(C_4 \square C_n)$ for all $n \geq 5$ such that $n \equiv 0, 2, 4, \textnormal{ and } 5 \pmod{6}$.

\begin{table}[htbp]
\centering
\renewcommand\arraystretch{1.1}
\begin{tabular}{rrcl}
&Lower Bound &  $ \gamma_{b,2} \left( C_{m} \square C_{n} \right) $  & Upper Bound\\ \hline
$1.6 \left( \frac{ n }{ 3 }  \right) \leq$ & $ \left\lceil 6\left(\frac{ n }{ 11 }\right)\right\rceil\leq $ & $ \gamma_{b,2} \left( C_{3} \square C_{n} \right) $ & $\leq \left\lceil \frac{ 2n }{ 3 }  \right\rceil$\\
& $ \left\lceil 4\left(\frac{ n }{ 6 }\right)\right\rceil\leq $ & $ \gamma_{b,2} \left( C_{4} \square C_{n} \right) $ & $\leq 4 \left\lfloor \frac{ n }{ 6 }  \right\rfloor + O(1)$\\
&$\left\lceil 10\left(\frac{ n }{ 13}\right)\right\rceil\leq $ & $ \gamma_{b,2} \left( C_{5} \square C_{n} \right) $ & $\leq n$\\
&$ \left\lceil 12\left(\frac{ n }{ 13}\right)\right\rceil\leq $ & $ \gamma_{b,2} \left( C_{6} \square C_{n} \right) $ & $\leq n + O(1)$\\
$37.6\left( \frac{ n }{ 35 }  \right) \leq$ & $ \left\lceil 14\left(\frac{ n }{ 13}\right)\right\rceil\leq $ & $ \gamma_{b,2} \left( C_{7} \square C_{n} \right) $ & $\leq 42 \left\lfloor \frac{ n }{ 35 }  \right\rfloor  + O(1)$\\
$7.3 \left( \frac{ n }{ 6 }  \right) \leq$ &$\left\lceil 16\left(\frac{ n }{ 13}\right)\right\rceil\leq $ & $ \gamma_{b,2} \left( C_{8} \square C_{n} \right) $ & $\leq 8 \left\lfloor \frac{ n }{ 6 }  \right\rfloor  + O(1)$\\
$13.8 \left( \frac{ n }{ 10 }  \right) \leq$ & $ \left\lceil 18\left(\frac{ n }{ 13}\right)\right\rceil\leq $ & $ \gamma_{b,2} \left( C_{9} \square C_{n} \right) $ & $\leq 16 \left\lfloor \frac{ n }{ 10 }  \right\rfloor  + O(1)$\\
$15.3 \left( \frac{ n }{ 10 }  \right) \leq$ & $ \left\lceil 20\left(\frac{ n }{ 13}\right)\right\rceil\leq $ & $ \gamma_{b,2} \left( C_{10} \square C_{n} \right) $ & $\leq 16 \left\lfloor \frac{ n }{ 10 }  \right\rfloor  + O(1)$\\
&$\left\lceil 22\left(\frac{  n}{ 13}\right)\right\rceil\leq $ & $ \gamma_{b,2} \left( C_{11} \square C_{n} \right) $ & $\leq 24 \left\lfloor \frac{ n }{ 13 }  \right\rfloor  + O(1)$\\
&$\left\lceil 24\left(\frac{  n}{ 13}\right)\right\rceil\leq $ & $ \gamma_{b,2} \left( C_{12} \square C_{n} \right) $ & $\leq 26 \left\lfloor \frac{ n }{ 13 }  \right\rfloor  + O(1)$\\
&$\left\lceil 26\left(\frac{ n }{ 13}\right)\right\rceil\leq $ & $ \gamma_{b,2} \left( C_{13} \square C_{n} \right) $ & $\leq 26 \left\lfloor \frac{ n }{ 13 }  \right\rfloor  + O(1)$\\
&$\left\lceil 28\left(\frac{ n }{ 13}\right)\right\rceil\leq $ & $ \gamma_{b,2} \left( C_{14} \square C_{n} \right) $ & $\leq 31 \left\lfloor \frac{ n }{ 13 }  \right\rfloor  + O(1)$\\
&$\left\lceil 30\left(\frac{ n}{ 13 }\right)\right\rceil\leq $ & $ \gamma_{b,2} \left( C_{15} \square C_{n} \right) $ & $\leq 33 \left\lfloor \frac{ n }{ 13 }  \right\rfloor  + O(1)$\\
&$\left\lceil 32\left(\frac{ n}{ 13 }\right)\right\rceil\leq $ & $ \gamma_{b,2} \left( C_{16} \square C_{n} \right) $ & $\leq 36 \left\lfloor \frac{ n }{ 13 }  \right\rfloor  + O(1)$\\
&$\left\lceil 34\left(\frac{ n}{ 13 }\right)\right\rceil\leq $ & $ \gamma_{b,2} \left( C_{17} \square C_{n} \right) $ & $\leq 36 \left\lfloor \frac{ n }{ 13 }  \right\rfloor  + O(1)$\\
&$\left\lceil 36\left(\frac{ n}{ 13 }\right)\right\rceil\leq $ & $ \gamma_{b,2} \left( C_{18} \square C_{n} \right) $ & $\leq 39 \left\lfloor \frac{ n }{ 13 }  \right\rfloor  + O(1)$\\
&$\left\lceil 38\left(\frac{ n}{ 13 }\right)\right\rceil\leq $ & $ \gamma_{b,2} \left( C_{19} \square C_{n} \right) $ & $\leq 41 \left\lfloor \frac{ n }{ 13 }  \right\rfloor  + O(1)$\\
&$\left\lceil 40\left(\frac{ n}{ 13 }\right)\right\rceil\leq $ & $ \gamma_{b,2} \left( C_{20} \square C_{n} \right) $ & $\leq 42 \left\lfloor \frac{ n }{ 13 }  \right\rfloor  + O(1)$\\
&$\left\lceil 42\left(\frac{ n}{ 13 }\right)\right\rceil\leq $ & $ \gamma_{b,2} \left( C_{21} \square C_{n} \right) $ & $\leq 45 \left\lfloor \frac{ n }{ 13 }  \right\rfloor  + O(1)$\\
&$\left\lceil 44\left(\frac{ n}{ 13 }\right)\right\rceil\leq $ & $ \gamma_{b,2} \left( C_{22} \square C_{n} \right) $ & $\leq 48 \left\lfloor \frac{ n }{ 13 }  \right\rfloor  + O(1)$\\
&$\left\lceil 46 \left( \frac{ n}{ 13 } \right)\right\rceil \leq $ & $ \gamma_{b,2} \left( C_{23} \square C_{n} \right) $ & $\leq 49 \left\lfloor \frac{ n }{ 13 }  \right\rfloor  + O(1)$\\
&$\left\lceil48 \left(  \frac{ n}{ 13 } \right)\right\rceil \leq $ & $ \gamma_{b,2} \left( C_{24} \square C_{n} \right) $ & $\leq 50 \left\lfloor \frac{ n }{ 13 }  \right\rfloor  + O(1)$\\
&$\left\lceil 50 \left( \frac{ n}{ 13 }  \right)\right\rceil \leq $ & $ \gamma_{b,2} \left( C_{25} \square C_{n} \right) $ & $\leq 52 \left\lfloor \frac{ n }{ 13 }  \right\rfloor  + O(1)$\\
&$\left\lceil2  \left(  \frac{ mn}{ 13 }\right)\right\rceil \leq $ & $ \gamma_{b,2} \left( C_{m\geq 26} \square C_{n \geq 26} \right) $ & $\leq 2 \left( \frac{ mn }{ 13 }  \right) + O(m+n)$\\
\end{tabular}
\caption{Upper and lower bounds for $ \gamma_{b,2} \left( C_{m} \square C_{n} \right) $ for $m,n \geq 3$ excluding $m=n=3$, $m=n=4$, and $m=4$ and $n=3$, the values of which are given in Theorem \ref{prop:mulitpackcycle}.}
\label{table:multipackcycle}
\end{table}
% default
\renewcommand\arraystretch{1}

\subsection{The Cartesian Product of a Path and a Cycle}\label{sec:multipack:pathxcycle}
This section describes the construction and associated costs of fractional $2$-limited multipackings on $ P_{m} \square C_{n} $ for $m \geq 2$ and $n \geq 3$. For each $2 \leq m \leq 22$, an explicit fractional $2$-limited multipacking is given on $ P_m \square C_n$ for all $n \geq 3$. These multipackings were found by examining fractional $2$-limited multipackings on $ P_{m} \square C_{5} $ for $2 \leq m \leq 22$. These constructions are used to prove Theorem \ref{thm:multipack222}. For $m \geq 23$ and all $n \geq 3$ a fractional $2$-limited multipacking is given on $ P_m \square C_n$. This multipacking was determined by examining fractional $2$-limited multipackings on $ P_{23} \square C_{5} $. This construction is used to prove Theorem \ref{thm:multipack:pathxcyclegeneral}.

\begin{thm}\label{thm:multipack222}
For $2 \leq m \leq 22$ and all $n \geq 3$, $f(n)  \leq mp_{f,2}(P_m \square C_n)$ where $f(n)$ is given in Table \ref{table:multipack222}.
\end{thm}

\begin{table}[htbp]
\centering
\renewcommand\arraystretch{1.1}
% used to be dfrac
\begin{tabular}{c|ccccccccc}
$m$ & $2$ & $3$ & $4$ & $5$ & $6$ &$ 7$ &$ 8$  \\ \hline
$f(n)$ &    
$\frac{ n }{ 2 }  $ &
$\frac{ 2n }{ 3 } $ &
$\frac{ 4n }{ 5 } $ & 
$\frac{ 26n }{ 27 }  $ & 
$\frac{ 29n }{ 26 }  $ &
$\frac{ 19n }{ 15 }  $ & 
$\frac{ 212n }{ 149 }  $  \\
\multicolumn{8}{c}{}\\ 
$m$ & $9$ & $10$  & $11$ &$ 12$ & $13$ & $14$ & $15$ \\ \hline
$f(n)$  &
$\frac{ 52n }{ 33 }  $ &
$\frac{ 780n }{ 451 }$ &
$\frac{ 81n }{ 43 }  $ &
$\frac{ 273n }{ 134 }  $&
$\frac{ 5042n }{ 2301 }$ &
$\frac{ 2324n }{ 991 }$ &
$\frac{ 5690n}{ 2277 }$ \\
\multicolumn{8}{c}{}\\
$m$& $16$ & $17$   & $18$ & $19$ & $20$ & $21$ &$22$  \\ \hline
$f(n)$&
$\frac{ 15593n }{ 5878 } $ &
$\frac{ 28417n }{ 10125 } $& 
$\frac{ 103240n }{ 34873 }$ &
$\frac{ 3896n}{ 1251 }  $ &
$\frac{ 337976n }{ 103415 }$ &
$\frac{ 304705n }{ 89043 }  $& 
$\frac{ 548313n }{ 153338 } $ 
\end{tabular}
\caption{Values of $f(n)$ for lower bounds for $ mp_{2} \left( P_{m} \square C_{n} \right) $ where $2 \leq m \leq 22$ and $n \geq 3$.}
\label{table:multipack222}
\end{table}
% default
\renewcommand\arraystretch{1}

\begin{pf}
Fix $ 2 \leq m \leq 22$. Let $ \textbf{v} $ be the $m \times 1$ vector as given in either Table \ref{table:multipackvector213} (for $2 \leq m \leq 13$) or Table \ref{table:multipackvectoer1422} (for $m \geq 14$, see Section \ref{appendix}). Let $ \textbf{J}_{1,n}$ be the all ones matrix of size $1 \times n$. 
\begin{table}[htbp]
\renewcommand\arraystretch{1.2}
\[
\arraycolsep=3.5pt
\begin{array}{ccccccccccccc}
\multicolumn{12}{c}{m}\\
$2$ & $3$ & $4$ & $5$ & $6$ & $7$ & $8$ & $9$ & $10$ & $11$ & $12$ & $13$\\
\hline
\belowbaseline[-8pt]{$\begin{bmatrix}\frac{1}{4}\\  \frac{1}{4}\\  \end{bmatrix}$}
&
\belowbaseline[-8pt]{$\begin{bmatrix}\frac{1}{3}\\  0\\  \frac{1}{3}\\  \end{bmatrix}$}
&
\belowbaseline[-8pt]{$\begin{bmatrix}\frac{3}{10}\\  \frac{1}{10}\\  \frac{1}{10}\\  \frac{3}{10}\\  \end{bmatrix}$}
&
\belowbaseline[-8pt]{$\begin{bmatrix}\frac{8}{27}\\  \frac{1}{9}\\  \frac{4}{27}\\  \frac{1}{9}\\  \frac{8}{27}\\  \end{bmatrix}$}
&
\belowbaseline[-8pt]{$\begin{bmatrix}\frac{4}{13}\\  \frac{1}{13}\\  \frac{9}{52}\\  \frac{9}{52}\\  \frac{1}{13}\\  \frac{4}{13}\\  \end{bmatrix}$}
&
\belowbaseline[-8pt]{$\begin{bmatrix}\frac{3}{10}\\  \frac{1}{10}\\  \frac{2}{15}\\  \frac{1}{5}\\  \frac{2}{15}\\  \frac{1}{10}\\  \frac{3}{10}\\  \end{bmatrix}$}
& 
\belowbaseline[-8pt]{$\begin{bmatrix}\frac{45}{149}\\  \frac{14}{149}\\  \frac{23}{149}\\  \frac{24}{149}\\  \frac{24}{149}\\  \frac{23}{149}\\  \frac{14}{149}\\  \frac{45}{149}\\  \end{bmatrix}$}
&
\belowbaseline[-8pt]{$\begin{bmatrix}\frac{10}{33}\\  \frac{1}{11}\\  \frac{5}{33}\\  \frac{2}{11}\\  \frac{4}{33}\\  \frac{2}{11}\\  \frac{5}{33}\\  \frac{1}{11}\\  \frac{10}{33}\\  \end{bmatrix}$}
&
\belowbaseline[-8pt]{$\begin{bmatrix}\frac{136}{451}\\  \frac{43}{451}\\  \frac{6}{41}\\  \frac{81}{451}\\  \frac{64}{451}\\  \frac{64}{451}\\  \frac{81}{451}\\  \frac{6}{41}\\  \frac{43}{451}\\  \frac{136}{451}\\  \end{bmatrix}$}
&
\belowbaseline[-8pt]{$\begin{bmatrix}\frac{13}{43}\\  \frac{4}{43}\\  \frac{13}{86}\\  \frac{15}{86}\\  \frac{6}{43}\\  \frac{7}{43}\\  \frac{6}{43}\\  \frac{15}{86}\\  \frac{13}{86}\\  \frac{4}{43}\\  \frac{13}{43}\\  \end{bmatrix}$}
&
\belowbaseline[-8pt]{$\begin{bmatrix}\frac{81}{268}\\  \frac{25}{268}\\  \frac{10}{67}\\  \frac{12}{67}\\  \frac{9}{67}\\  \frac{43}{268}\\  \frac{43}{268}\\  \frac{9}{67}\\  \frac{12}{67}\\  \frac{10}{67}\\  \frac{25}{268}\\  \frac{81}{268}\\  \end{bmatrix}$}
&
\belowbaseline[-8pt]{$\begin{bmatrix}\frac{695}{2301}\\  \frac{72}{767}\\  \frac{343}{2301}\\  \frac{136}{767}\\  \frac{320}{2301}\\  \frac{119}{767}\\  \frac{28}{177}\\  \frac{119}{767}\\  \frac{320}{2301}\\  \frac{136}{767}\\  \frac{343}{2301}\\  \frac{72}{767}\\  \frac{695}{2301}\\  \end{bmatrix}$}
\end{array}
\]
\caption{Vectors used in the proof of the lower bound for $mp_{f,2}(P_m \square C_n)$ for $2 \leq m \leq 13$}
\label{table:multipackvector213}
\end{table}
% default
\renewcommand\arraystretch{1}

The product $\textbf{v}\textbf{J}_{1,n}$ defines an $m \times n$ matrix $ \textbf{M} $ where the elements within each row, respectively, have identical entries. For example if $m=3$ and $n=6$ then
\renewcommand\arraystretch{1.2}
$$ 
\textbf{M} =\begin{bmatrix}
    \frac{1}{3} & \frac{1}{3} & \frac{1}{3} & \frac{1}{3} & \frac{1}{3} & \frac{1}{3}\\
    0 & 0 & 0 & 0 & 0 & 0\\
    \frac{1}{3} & \frac{1}{3} & \frac{1}{3} & \frac{1}{3} & \frac{1}{3} & \frac{1}{3}\\
\end{bmatrix}.
$$
% default
\renewcommand\arraystretch{1}
For $1 \leq i \leq m$ and $1 \leq j \leq n$, define the fractional $2$-limited multipacking $f$ on $ P_{m} \square C_{n} $ by letting the variable $x_{i,j}$ (corresponding to the vertex in row $i$ and column $j$ of $ P_{m} \square C_{n} $) be the value  in the $i$th row and $j$th column of $ \textbf{M} $. As the elements within each row are identical, due to symmetry, to check that $f$ is a valid fractional $2$-limited multipacking it suffices to check the entries in the first column of $ P_{m} \square C_{n} $. The check can be done by computation. Let $v_1, \dots v_m$ denote the elements of $ \textbf{v} $. The cost of $f$ is $n\sum_{i=1}^{ m } v_i.$ For each $2 \leq m \leq 22$ the cost of $f$ corresponds with the values in Table \ref{table:multipack222}. This completes the proof.
\end{pf}

\begin{thm}\label{thm:multipack:pathxcyclegeneral}
For $m \geq 23$ and all $n \geq 3$, 
\begin{equation*}
\begin{aligned}
\frac{ 2mn }{ 13 } + \frac{ 2620n}{ 13767}  \leq mp_{f,2}( P_{m} \square C_{n}).
\end{aligned}
\end{equation*}
\end{thm}

\begin{pf}
The general case follows from the case  $m=23$. Letting $ \textbf{v} $ be the vector in Figure \ref{Figure:multipackvector23} for $m=23$ and proceeding with as in the proof of Theorem \ref{thm:multipack222} (with $ \textbf{v}^T$) yields a fractional $2$-limited multipacking on $ P_{23} \square C_{n} $ for all $n \geq 3$. 
\begin{figure}[htbp]
\setcounter{MaxMatrixCols}{19}
\setlength\arraycolsep{1.5pt}
$$\begin{bmatrix}    
    \frac{16639}{55068} & \frac{1717}{18356} &\frac{686}{4589} &\frac{814}{4589} &\frac{1895}{13767} &\frac{2903}{18356}&\frac{8573}{55068}&\frac{697}{4589} & \undermat{\substack{m-16\\ \textnormal{entries of } \frac{2}{13}}}{\frac{2}{13} & \dots & \frac{2}{13}} &  \frac{697}{4589} &\frac{8573}{55068} &\frac{2903}{18356} &\frac{1895}{13767} &\frac{814}{4589} &\frac{686}{4589} &\frac{1717}{18356} &\frac{16639}{55068}\\   
\end{bmatrix}$$
\caption{Vector used in the proof of the lower bound for $mp_{f,2}(P_m \square C_n)$ for $ m \geq 23$}
\label{Figure:multipackvector23}
\end{figure}

The cost of such a multipacking is
\begin{equation*}
\begin{aligned}
\left( m-16 \right) \left( \frac{ 2n }{ 13 }  \right) + \frac{ 36508n }{ 13767 } =7\left( \frac{ 2n }{ 13 }  \right) + \frac{ 36508n }{ 13767 }.
\end{aligned}
\end{equation*}
Observe that the center most row of this multipacking on $ P_{23} \square C_{n} $ is composed solely of vertices with weight $2/13$. Moreover, all vertices within distance 2 of the center row also have weight $2/13$. As this row does not violate the fractional $2$-limited multipacking on $ P_{23} \square C_{n} $, adding additional center rows, all of which contain vertices with weight $2/13$, will yield valid fractional $2$-limited multipacking  on $ P_{m \geq 23} \square C_{n} $. The cost of such a multipacking will be
\begin{equation*}
\begin{aligned}
\left( m-16 \right) \left( \frac{ 2n }{ 13 }  \right) + \frac{ 36508n }{ 13767 },  
\end{aligned}
\end{equation*}
which when simplified proves the theorem.
\end{pf}

Combined, the results found in Section \ref{sec:pathsandcycle} and Theorems \ref{prop:mulitpackcycle}, \ref{thm:multipack222}, and \ref{thm:multipack:pathxcyclegeneral} establish the upper and lower bounds in Table \ref{table:multipackpathxcycle} for the $2$-limited broadcast domination number of the Cartesian products of a path and a cycle. Note that Theorem \ref{prop:mulitpackcycle} is used in the case of $P_{m} \square C_3$ and  $P_{m} \square C_4$ for $m \geq 23$ since $mp_{f,2}(C_m\square C_n) \leq mp_{f,2}(P_m\square C_n)$. The value of the constant terms in the upper bounds in Table \ref{table:multipackpathxcycle} can be found in the respective theorems in Section \ref{sec:pathsandcycle} and are omitted here. To easily compare the upper and lower bounds in Table \ref{table:multipackpathxcycle} we also, for some bounds, include a simpler lower bound.

\begin{table}[htbp]
\centering
\renewcommand\arraystretch{1.1}
\begin{tabular}{rrcl}
&    Lower Bound &  $ \gamma_{b,2} \left( P_{m} \square C_{n} \right) $  & Upper Bound\\ \hline
& $\left\lceil \frac{ n }{ 2 }  \right\rceil\leq $ & $ \gamma_{b,2} \left( P_{2} \square C_{n} \right) $ & $\leq \left\lceil \frac{ n }{ 2 }  \right\rceil $\\
& $\left\lceil \frac{ 2n }{ 3 }   \right\rceil\leq $ & $ \gamma_{b,2} \left( P_{3} \square C_{n} \right) $ & $\leq \left\lceil \frac{ 2n }{ 3 }  \right\rceil $\\
&    $\left\lceil \frac{ 4n }{ 5 }  \right\rceil\leq $ & $ \gamma_{b,2} \left( P_{4} \square C_{n} \right) $ & $\leq 8 \left\lfloor \frac{ n }{ 10 }  \right\rfloor + O(1)$\\
&    $\left\lceil \frac{ 26n }{ 27 }  \right\rceil\leq $ & $ \gamma_{b,2} \left( P_{5} \square C_{n} \right) $ & $\leq n + O(1)$\\
$17.8 \left( \frac{ n }{ 16 }  \right) \leq $ & $\left\lceil \frac{ 29n }{ 26 }  \right\rceil\leq $ & $ \gamma_{b,2} \left( P_{6} \square C_{n} \right) $ & $\leq 18 \left\lfloor \frac{ n }{ 16 }  \right\rfloor  + O(1)$\\
$17.7 \left( \frac{ n }{ 14 }  \right) \leq $ & $\left\lceil \frac{ 19n }{ 15 }  \right\rceil\leq $ & $ \gamma_{b,2} \left( P_{7} \square C_{n} \right) $ & $\leq 18 \left\lfloor \frac{ n }{ 14 }  \right\rfloor   + O(1)$\\
$31.3 \left( \frac{ n }{ 22 }  \right) \leq $ & $\left\lceil \frac{ 212n }{ 149 }   \right\rceil\leq $ & $ \gamma_{b,2} \left( P_{8} \square C_{n} \right) $ & $\leq 32 \left\lfloor \frac{ n }{ 22 }  \right\rfloor   + O(1)$\\
$15.7 \left( \frac{ n }{ 10 }  \right) \leq $ & $\left\lceil \frac{ 52n }{ 33 }  \right\rceil\leq $ & $ \gamma_{b,2} \left( P_{9} \square C_{n} \right) $ & $\leq 16 \left\lfloor \frac{ n }{ 10 }  \right\rfloor  + O(1)$\\
$31.1 \left( \frac{ n }{ 18 }  \right) \leq$ & $\left\lceil \frac{ 780n }{ 451 }  \right\rceil\leq $ & $ \gamma_{b,2} \left( P_{10} \square C_{n} \right) $ & $\leq 32 \left\lfloor \frac{ n }{ 18 }  \right\rfloor  + O(1)$\\
$48.9 \left( \frac{ n }{ 26 }  \right) \leq $ & $\left\lceil \frac{ 81n }{ 43 }  \right\rceil\leq $ & $ \gamma_{b,2} \left( P_{11} \square C_{n} \right) $ & $\leq 50 \left\lfloor \frac{ n }{ 26 }  \right\rfloor  + O(1)$\\
$48.8 \left( \frac{ n }{ 24 }  \right) \leq$ & $\left\lceil \frac{ 273n }{ 134 }  \right\rceil\leq $ & $ \gamma_{b,2} \left( P_{12} \square C_{n} \right) $ & $\leq 50 \left\lfloor \frac{ n }{ 24 }  \right\rfloor  + O(1)$\\
$28.4 \left( \frac{ n }{ 13 }  \right) \leq $ & $\left\lceil \frac{ 5042n }{ 2301 }  \right\rceil\leq $ & $ \gamma_{b,2} \left( P_{13} \square C_{n} \right) $ & $\leq 30 \left\lfloor \frac{ n }{ 13 }  \right\rfloor  + O(1)$\\
$30.4 \left( \frac{ n }{ 13 }  \right) \leq$ & $ \left\lceil  \frac{ 2324n }{ 991 } \right\rceil\leq $ & $ \gamma_{b,2} \left( P_{14} \square C_{n} \right) $ & $\leq 32 \left\lfloor \frac{ n }{ 13 }  \right\rfloor  + O(1)$\\
$32.4 \left( \frac{ n }{ 13 }  \right) \leq$ & $\left\lceil  \frac{ 5690n}{ 2277 } \right\rceil\leq $ & $ \gamma_{b,2} \left( P_{15} \square C_{n} \right) $ & $\leq 34 \left\lfloor \frac{ n }{ 13 }  \right\rfloor  + O(1)$\\
$34.4 \left( \frac{ n }{ 13 }  \right) \leq$ & $\left\lceil  \frac{ 15593n }{ 5878 } \right\rceil\leq $ & $ \gamma_{b,2} \left( P_{16} \square C_{n} \right) $ & $\leq 36 \left\lfloor \frac{ n }{ 13 }  \right\rfloor  + O(1)$\\
$36.4 \left( \frac{ n }{ 13 }  \right) \leq$ & $\left\lceil \frac{ 28417n }{ 10125 }  \right\rceil\leq $ & $ \gamma_{b,2} \left( P_{17} \square C_{n} \right) $ & $\leq 38 \left\lfloor \frac{ n }{ 13 }  \right\rfloor  + O(1)$\\
$38.4 \left( \frac{ n }{ 13 }  \right) \leq$ & $\left\lceil \frac{ 103240n }{ 34873 }  \right\rceil\leq $ & $ \gamma_{b,2} \left( P_{18} \square C_{n} \right) $ & $\leq 40 \left\lfloor \frac{ n }{ 13 }  \right\rfloor  + O(1)$\\
$40.4 \left( \frac{ n }{ 13 }  \right) \leq$ & $\left\lceil \frac{ 3896n}{ 1251 }  \right\rceil\leq $ & $ \gamma_{b,2} \left( P_{19} \square C_{n} \right) $ & $\leq 42 \left\lfloor \frac{ n }{ 13 }  \right\rfloor  + O(1)$\\
$42.4 \left( \frac{ n }{ 13 }  \right) \leq$ & $\left\lceil  \frac{ 337976n }{ 103415 }  \right\rceil\leq $ & $ \gamma_{b,2} \left( P_{20} \square C_{n} \right) $ & $\leq 44 \left\lfloor \frac{ n }{ 13 }  \right\rfloor  + O(1)$\\
$40.4 \left( \frac{ n }{ 13 }  \right) \leq$ & $\left\lceil \frac{ 304705n }{ 89043 }  \right\rceil\leq $ & $ \gamma_{b,2} \left( P_{21} \square C_{n} \right) $ & $\leq 46 \left\lfloor \frac{ n }{ 13 }  \right\rfloor  + O(1)$\\
$46.4 \left( \frac{ n }{ 13 }  \right) \leq$ & $\left\lceil \frac{ 548313n }{ 153338 }   \right\rceil\leq $ & $ \gamma_{b,2} \left( P_{22} \square C_{n} \right) $ & $\leq 48 \left\lfloor \frac{ n }{ 13 }  \right\rfloor  + O(1)$\\
                         &    $\left\lceil2  \left(  \frac{ mn}{ 13 }\right) + \frac{ 2620n }{ 13767 } \right\rceil \leq $ & $ \gamma_{b,2} \left( P_{m\geq 23} \square C_{n\geq 13} \right) $ & $\leq 2 \left( \frac{ mn }{ 13 }  \right) + \frac{ 4m }{ 13 } + O(n)$\\
$1.6 \left( \frac{ m }{ 3 }  \right) \leq$ & $ \left\lceil 6\left(\frac{ m }{ 11 }\right) \right\rceil\leq $ & $ \gamma_{b,2} \left( P_{m \geq 23} \square C_{3} \right) $ & $\leq \left\lceil \frac{ 2m }{ 3 }  \right\rceil $\\
%$1.3 \left( \frac{ m }{ 3 }  \right) \leq$ & $ \left\lceil \frac{ 6m }{ 13 } + \frac{ 2620 }{ 4586 }  \right\rceil\leq $ & $ \gamma_{b,2} \left( P_{m \geq 23} \square C_{3} \right) $ & $\leq \left\lceil \frac{ 2m }{ 3 }  \right\rceil $\\
& $ \left\lceil 4 \left(\frac{m}{6} \right)  \right\rceil \leq$ & $ \gamma_{b,2} \left( P_{m \geq 23} \square C_{4} \right) $ & $\leq 4 \left\lfloor \frac{ m }{ 6 }  \right\rfloor + O(1)$\\
%$3.6 \left( \frac{ m }{ 3 }  \right) \leq $ & $\left\lceil \frac{ 8m }{ 13 } + \frac{10480}{ 13767 }  \right\rceil\leq $ & $ \gamma_{b,2} \left( P_{m \geq 23} \square C_{4} \right) $ & $\leq 4 \left\lfloor \frac{ m }{ 6 }  \right\rfloor + O(1)$\\
                        &$\left\lceil \frac{ 10m }{13  }+ \frac{ 13100 }{ 13767 }    \right\rceil\leq $ & $ \gamma_{b,2} \left( P_{m \geq 23} \square C_{5} \right) $ & $\leq m+1$\\
                        &$\left\lceil \frac{ 12m }{ 13 } + \frac{ 5240 }{ 4589 }  \right\rceil\leq $ & $ \gamma_{b,2} \left( P_{m \geq 23} \square C_{6} \right) $ & $\leq m + 1$\\
$37.6 \left( \frac{ m }{ 35 }  \right) \leq $ & $\left\lceil \frac{ 14m}{ 13 } + \frac{ 18340 }{ 13767 }  \right\rceil\leq $ & $ \gamma_{b,2} \left( P_{m \geq 23} \square C_{7} \right) $ & $\leq 42 \left\lfloor \frac{ m }{ 35 }  \right\rfloor  + O(1)$\\
$7.3 \left( \frac{ n }{ 6 }  \right) \leq $ & $\left\lceil \frac{ 16m }{ 13 } + \frac{ 20960 }{ 13767 }  \right\rceil\leq $ & $ \gamma_{b,2} \left( P_{m \geq 23} \square C_{8} \right) $ & $\leq 8 \left\lfloor \frac{ m }{ 6 }  \right\rfloor  + O(1)$\\
$13.8 \left( \frac{ n }{ 10 }  \right) \leq $ & $\left\lceil \frac{ 18m }{ 13 } + \frac{ 7860 }{ 4589 }  \right\rceil\leq $ & $ \gamma_{b,2} \left( P_{m \geq 23} \square C_{9} \right) $ & $\leq 16 \left\lfloor \frac{ m }{ 10 }  \right\rfloor  + O(1)$\\
$15.3 \left( \frac{ n }{ 10 }  \right) \leq$ & $\left\lceil \frac{ 20m }{ 13 } + \frac{ 26200 }{ 13767 }  \right\rceil\leq $ & $ \gamma_{b,2} \left( P_{m \geq 23} \square C_{10} \right) $ & $\leq 16 \left\lfloor \frac{ m }{ 10 }  \right\rfloor  + O(1)$\\
                         & $ \left\lceil \frac{ 22m }{ 13 } + \frac{ 28820 }{ 13767 }  \right\rceil\leq $ & $ \gamma_{b,2} \left( P_{m \geq 23} \square C_{11} \right) $ & $\leq 24 \left\lfloor \frac{ m }{ 13 }  \right\rfloor  + O(1)$\\
                         & $  \left\lceil \frac{ 24m }{ 13 } + \frac{ 10480 }{ 4589 }   \right\rceil\leq $ & $ \gamma_{b,2} \left( P_{m \geq 23} \square C_{12} \right) $ & $\leq 26 \left\lfloor \frac{ m }{ 13 }  \right\rfloor  + O(1)$\\
\end{tabular}
\caption{Upper and lower bounds for $ \gamma_{b,2} \left( P_{m} \square C_{n} \right) $ for $m,n \geq 3$.}
\label{table:multipackpathxcycle}
\end{table}
% default
\renewcommand\arraystretch{1}

Observe that the bounds in Table \ref{table:multipackpathxcycle} give exact values for $ \gamma_{b,2} \left( P_{2} \square C_{n} \right) $ and $ \gamma_{b,2} \left( P_{3} \square C_{n} \right) $. These bounds also give optimal values for $ \gamma_{b,2} \left( P_{4} \square C_{n} \right) $ for all $n \geq 3$ such that $n \not \equiv 1 \textnormal{ or } 5  \pmod{ 10 } $ and $ \gamma_{b,2} \left( P_{m} \square C_{4} \right) $ for all $m \geq 23$ such that $m \equiv 2 \textnormal{ and } 5  \pmod{ 6 } $.

\subsection{The Cartesian Product of Two Paths}\label{sec:multipack:path}
This section determines lower bounds for $\gamma_{b,2}(P_m \square P_n)$ via fractional $2$-limited multipackings. For $2 \leq m\leq 22$, Section \ref{sec:multipack:pathxcycle} provides lower bounds for $ \gamma_{b,2} \left( P_{m} \square P_{n} \right) $ since $ mp_{f,2} \left( P_{m} \square C_{n}  \right) \leq mp_{f,2} \left( P_{m} \square P_{n}  \right) $. For $m,n \geq 23$, a construction very similar to the proof of Theorem \ref{thm:multipack:pathxcyclegeneral} can be found in  \cite[Section 5.4]{slobodin} and is omitted here for the sake of brevity. This construction establishes Theorem \ref{thm:mutipack:pathxpathgeneral}.

\begin{thm}\label{thm:mutipack:pathxpathgeneral}
For $m,n \geq 23$, 
\begin{equation*}
\begin{aligned}
\frac{2mn}{ 13 } + \frac{ 14287568 }{75254411  } \left( m+n \right) - \frac{ 177612468 }{ 978307343 }  \leq mp_{f,2}( P_{m} \square P_{n} ).
\end{aligned}
\end{equation*}
\end{thm}

Combined, the results found in Section \ref{sec:paths} and Theorems \ref{thm:multipack222} and \ref{thm:mutipack:pathxpathgeneral} establish the upper and lower bounds, respectively, in Table \ref{table:multipackpathxpath} for the $2$-limited broadcast domination numbers of the Cartesian products of two  paths. The value of the constant terms in the upper bounds in Table \ref{table:multipackpathxpath} can be found in the respective theorems in Section \ref{sec:paths} and are omitted here. To easily compare the upper and lower bounds in Table \ref{table:multipackpathxpath} we also, for some bounds, include a simpler lower bound.

\begin{table}[htbp]
\centering
\renewcommand\arraystretch{1.1}
\begin{tabular}{rrcl}
& Lower Bound &  $ \gamma_{b,2} \left( P_{m} \square P_{n} \right) $  & Upper Bound\\ \hline
&$\left\lceil \frac{ n }{ 2 }  \right\rceil\leq $ & $ \gamma_{b,2} \left( P_{2} \square P_{n} \right) $ & $\leq \left\lfloor \frac{ n }{ 2 }  \right\rfloor +1$\\
&$\left\lceil \frac{ 2n }{ 3 }   \right\rceil\leq $ & $ \gamma_{b,2} \left( P_{3} \square P_{n} \right) $ & $\leq \left\lceil \frac{ 2n }{ 3 }  \right\rceil$\\
&$\left\lceil \frac{ 4n }{ 5 }  \right\rceil\leq $ & $ \gamma_{b,2} \left( P_{4} \square P_{n} \right) $ & $\leq 8 \left\lfloor \frac{ n }{ 10 }  \right\rfloor + O(1)$\\
&$\left\lceil \frac{ 26n }{ 27 }  \right\rceil\leq $ & $ \gamma_{b,2} \left( P_{5} \square P_{n} \right) $ & $\leq n + 1$\\
$17.8 \left( \frac{ n }{ 16 }  \right) \leq$ & $ \left\lceil \frac{ 29n }{ 26 }  \right\rceil\leq $ & $ \gamma_{b,2} \left( P_{6} \square P_{n} \right) $ & $\leq 18 \left\lfloor \frac{ n }{ 16 }  \right\rfloor  + O(1)$\\
$17.7 \left( \frac{ n }{ 14 }  \right) \leq$& $ \left\lceil \frac{ 19n }{ 15 }  \right\rceil\leq $ & $ \gamma_{b,2} \left( P_{7} \square P_{n} \right) $ & $\leq 18 \left\lfloor \frac{ n }{ 14 }  \right\rfloor   + O(1)$\\
$31.3 \left( \frac{ n }{ 22 }  \right) \leq$& $ \left\lceil \frac{ 212n }{ 149 }   \right\rceil\leq $ & $ \gamma_{b,2} \left( P_{8} \square P_{n} \right) $ & $\leq 32 \left\lfloor \frac{ n }{ 22 }  \right\rfloor   + O(1)$\\
$15.7 \left( \frac{ n }{ 10 }  \right) \leq$& $ \left\lceil \frac{ 52n }{ 33 }  \right\rceil\leq $ & $ \gamma_{b,2} \left( P_{9} \square P_{n} \right) $ & $\leq 16 \left\lfloor \frac{ n }{ 10 }  \right\rfloor  + O(1)$\\
$31.1 \left( \frac{ n }{ 18 }  \right) \leq$& $\left\lceil \frac{ 780n }{ 451 }  \right\rceil\leq $ & $ \gamma_{b,2} \left( P_{10} \square P_{n} \right) $ & $\leq 32 \left\lfloor \frac{ n }{ 18 }  \right\rfloor  + O(1)$\\
$48.9 \left( \frac{ n }{ 26 }  \right) \leq$& $ \left\lceil \frac{ 81n }{ 43 }  \right\rceil\leq $ & $ \gamma_{b,2} \left( P_{11} \square P_{n} \right) $ & $\leq 50 \left\lfloor \frac{ n }{ 26 }  \right\rfloor  + O(1)$\\
$48.8 \left( \frac{ n }{ 24 }  \right) \leq$& $\left\lceil \frac{ 273n }{ 134 }  \right\rceil\leq $ & $ \gamma_{b,2} \left( P_{12} \square P_{n} \right) $ & $\leq 50 \left\lfloor \frac{ n }{ 24 }  \right\rfloor  + O(1)$\\
$28.4 \left( \frac{ n }{ 13 }  \right) \leq$& $ \left\lceil \frac{ 5042n }{ 2301 }  \right\rceil\leq $ & $ \gamma_{b,2} \left( P_{13} \square P_{n} \right) $ & $\leq 30 \left\lfloor \frac{ n }{ 13 }  \right\rfloor  + O(1)$\\
$30.4 \left( \frac{ n }{ 13 }  \right) \leq$& $ \left\lceil  \frac{ 2324n }{ 991 } \right\rceil\leq $ & $ \gamma_{b,2} \left( P_{14} \square P_{n} \right) $ & $\leq 32 \left\lfloor \frac{ n }{ 13 }  \right\rfloor  + O(1)$\\
$32.4 \left( \frac{ n }{ 13 }  \right) \leq$& $\left\lceil  \frac{ 5690n}{ 2277 } \right\rceil\leq $ & $ \gamma_{b,2} \left( P_{15} \square P_{n} \right) $ & $\leq 34 \left\lfloor \frac{ n }{ 13 }  \right\rfloor  + O(1)$\\
$34.4 \left( \frac{ n }{ 13 }  \right) \leq$& $\left\lceil  \frac{ 15593n }{ 5878 } \right\rceil\leq $ & $ \gamma_{b,2} \left( P_{16} \square P_{n} \right) $ & $\leq 36 \left\lfloor \frac{ n }{ 13 }  \right\rfloor  + O(1)$\\
$36.4 \left( \frac{ n }{ 13 }  \right) \leq$& $\left\lceil \frac{ 28417n }{ 10125 }  \right\rceil\leq $ & $ \gamma_{b,2} \left( P_{17} \square P_{n} \right) $ & $\leq 38 \left\lfloor \frac{ n }{ 13 }  \right\rfloor  + O(1)$\\
$38.4 \left( \frac{ n }{ 13 }  \right) \leq$& $\left\lceil \frac{ 103240n }{ 34873 }  \right\rceil\leq $ & $ \gamma_{b,2} \left( P_{18} \square P_{n} \right) $ & $\leq 40 \left\lfloor \frac{ n }{ 13 }  \right\rfloor  + O(1)$\\
$40.4 \left( \frac{ n }{ 13 }  \right) \leq$& $\left\lceil \frac{ 3896n}{ 1251 }  \right\rceil\leq $ & $ \gamma_{b,2} \left( P_{19} \square P_{n} \right) $ & $\leq 42 \left\lfloor \frac{ n }{ 13 }  \right\rfloor  + O(1)$\\
$42.4 \left( \frac{ n }{ 13 }  \right) \leq$& $\left\lceil  \frac{ 337976n }{ 103415 }  \right\rceil\leq $ & $ \gamma_{b,2} \left( P_{20} \square P_{n} \right) $ & $\leq 44 \left\lfloor \frac{ n }{ 13 }  \right\rfloor  + O(1)$\\
$40.4 \left( \frac{ n }{ 13 }  \right) \leq$& $\left\lceil \frac{ 304705n }{ 89043 }  \right\rceil\leq $ & $ \gamma_{b,2} \left( P_{21} \square P_{n} \right) $ & $\leq 46 \left\lfloor \frac{ n }{ 13 }  \right\rfloor  + O(1)$\\
$46.4 \left( \frac{ n }{ 13 }  \right) \leq$& $\left\lceil \frac{ 548313n }{ 153338 }   \right\rceil\leq $ & $ \gamma_{b,2} \left( P_{22} \square P_{n} \right) $ & $\leq 48 \left\lfloor \frac{ n }{ 13 }  \right\rfloor  + O(1)$\\
\multicolumn{2}{r}{$\left\lceil  \frac{2mn}{ 13 } + \frac{ 14287568 }{75254411  } \left( m+n \right) -O(1)  \right\rceil \leq $} & $ \gamma_{b,2} \left( P_{m\geq 23} \square P_{n} \right) $ & $\leq  \frac{2 mn }{ 13 } + \frac{ 4 }{ 13 } \left( m+n \right)  + O(1)$\\
\end{tabular}
\caption{Upper and lower bounds for $ \gamma_{b,2} \left( P_{m} \square P_{n} \right) $ for $m,n \geq 2$.}
\label{table:multipackpathxpath}
\end{table}
% default
\renewcommand\arraystretch{1}

Observe that the bounds in Table \ref{table:multipackpathxpath} give optimal values for $ \gamma_{b,2} \left( P_{2} \square P_{n} \right) $ for odd $n$, optimal values for $ \gamma_{b,2} \left( P_{3} \square P_{n} \right) $, and optimal values for $ \gamma_{b,2} \left( P_{4} \square P_{n} \right) $ when $n  \equiv 4 \textnormal{ or } 9 \pmod{10}$. The periodically optimal values for $ \gamma_{b,2} \left( P_{2} \square P_{n} \right) $ allow for an easy proof of Proposition \ref{prop:multipackP2Pn}. 

\begin{prop}\label{prop:multipackP2Pn}
For $n \geq 2$, $\gamma_{b,2} \left( P_{2} \square P_{n} \right)  =  \left\lfloor \frac{ n }{ 2 }  \right\rfloor +1.$
\end{prop}

\begin{pf}
If $n\in \left\{ 2,3,4,5 \right\} $ the claim is easily verified. Fix $n \geq 6$. Suppose $n$ is odd. As    
\begin{equation*}
\begin{aligned}
\left\lceil \frac{ n }{ 2 }  \right\rceil  & =   \frac{ n+1 }{ 2 } = \frac{ n-1 }{2 } +1 = \left\lfloor \frac{ n }{ 2 }  \right\rfloor +1,
\end{aligned}
\end{equation*}
by the bounds for $ \gamma_{b,2} \left( P_{2} \square P_{n} \right) $ in Theorem \ref{thm:multipack222} and Theorem \ref{thm:path}, the claim holds. 

Suppose now that $n$ is even and that there exists a $2$-limited dominating broadcast $f$ of $ P_{2} \square P_{n}$ of cost  $\leq \left\lfloor n/2 \right\rfloor $. Observe that a vertex broadcasting at strength one or two on $ P_{2} \square P_n$ can be heard by at most four or eight vertices, respectively. Therefore, given  $v \in V(P_{2} \square P_{n}) $, the broadcast from $v$ can be heard by at most $ 4f(v)$ vertices. Consider the leftmost column of $ P_{2} \square P_{n} $. Any $2$-limited broadcast of $ P_{2} \square P_{n} $ which dominates this column necessarily ``wastes'' a portion of a  broadcast from a vertex. By ``wastes'' we mean that, in dominating the leftmost column, there exists a vertex $v$ whose broadcast either overlaps with the broadcast from some other vertex $v'$ or is heard by at least one less than $4f(v)$ vertices. In either case, a broadcast from such a vertex is heard by at least one less than $4f(v)$ vertices. The same logic applies to the rightmost column of $ P_{2} \square P_{n} $. This follows since $n \geq 6$. Therefore, there cannot exist a vertex which dominates a vertex in the left and rightmost columns. By assumption therefore, $f$ can dominate at most 
\begin{equation*}
\begin{aligned}
4 \left\lfloor \frac{ n }{ 2 }  \right\rfloor -2 =4 \left( \frac{ n }{ 2 }  \right) -2 = 2n-2 < 2n = \left| V( P_{2} \square P_{n} ) \right| 
\end{aligned}
\end{equation*}
vertices. This is a contradiction.
\end{pf}

\section{Conclusion}\label{sec:conclusion}
This paper summarizes the results from Chapters 2 through 5 of \cite{slobodin}. These results establish upper and lower bounds for the $2$-limited broadcast domination number of the Cartesian product of two paths, a path and a cycle, and two cycles.  Although we have bounded the $2$-limited broadcast domination number for the Cartesian products of two paths, a path and a cycle, and two cycles, few of our bounds are tight. Additionally, given a graph $G$, since $\gamma_{b,k}(G) \leq \gamma_{b,k-1}(G)$, our bounds provide upper bounds for the $k$-limited broadcast domination numbers for these graphs where $k \geq 2$. However, these bounds are very likely far from the truth with respect to large grids. Given this, we present the following problems.

\begin{problem}
Determine $\gamma_{b,k} \left(P_m \square P_n\right) $ for all $m, n \geq 2$ for $k \geq 2$.
\end{problem}
\begin{problem}
Determine $\gamma_{b,k} \left(P_m \square C_n\right) $ for all $m \geq 2$ and $n \geq 3$ for $k \geq 2$.
\end{problem}
\begin{problem}
Determine $\gamma_{b,k} \left(C_m \square C_n\right) $ for all $m, n \geq 3$ for $k \geq 2$.
\end{problem}

\section*{Acknowledgements} All computations for this paper we run on machines purchased by Myrvold using NSERC funding. We would also like to thank the referees for their thoughtful and instructive comments.

\bibliographystyle{plain} 
\bibliography{mybib} 

\section{Tables 18 -- 25} \label{appendix}

\begin{table}[htbp]
\centering
\begin{tabular}{cc|ccccccccccccc}
& & \multicolumn{13}{c}{Least residue $n_{13}$ of $n$ modulo $13$} \tabularnewline
& & $0$ & $1$ & $2$ & $3$ & $4$ & $5$ & $6$ &$ 7$ &$ 8$ &$ 9$ & $10$ & $11$ &$ 12$\\ \hline
\parbox[t]{2mm}{\multirow{13}{*}{\rotatebox[origin=c]{90}{Least residue $m_{13}$ of $m$ modulo $13$}}} 
&0 & $4$ &$10$ &$10$ &$10$ &$12$ &$10$ &$10$ &$12$ &$10$ &$12$ &$11$ &$10$ &$4$ \tabularnewline
& 1 & $11$ &$12$ &$11$ &$10$ &$12$ &$11$ &$12$ &$11$ &$7$ &$11$ &$11$ &$11$ &$11$ \tabularnewline
& 2 & $10$ &$12$ &$10$ &$10$ &$12$ &$10$ &$12$ &$1$ &$10$ &$12$ &$11$ &$12$ &$11$ \tabularnewline
& 3 & $4$ &$4$ &$4$ &$12$ &$11$ &$10$ &$12$ &$4$ &$11$ &$11$ &$11$ &$12$ &$11$ \tabularnewline
& 4 & $10$ &$12$ &$10$ &$12$ &$1$ &$10$ &$12$ &$7$ &$9$ &$11$ &$7$ &$12$ &$4$ \tabularnewline
& 5 & $10$ &$10$ &$10$ &$12$ &$10$ &$10$ &$12$ &$10$ &$12$ &$11$ &$10$ &$12$ &$4$ \tabularnewline
& 6 & $0$ &$11$ &$10$ &$12$ &$11$ &$10$ &$12$ &$7$ &$9$ &$11$ &$10$ &$12$ &$11$ \tabularnewline
& 7 & $12$ &$10$ &$10$ &$12$ &$10$ &$12$ &$1$ &$10$ &$12$ &$11$ &$10$ &$12$ &$10$ \tabularnewline
& 8 & $6$ &$4$ &$10$ &$12$ &$10$ &$12$ &$4$ &$10$ &$12$ &$11$ &$0$ &$11$ &$4$ \tabularnewline
& 9 & $12$ &$10$ &$12$ &$1$ &$10$ &$12$ &$7$ &$9$ &$11$ &$7$ &$9$ &$11$ &$11$ \tabularnewline
& 10 & $12$ &$10$ &$12$ &$12$ &$10$ &$12$ &$12$ &$12$ &$1$ &$10$ &$12$ &$12$ &$12$ \tabularnewline
& 11 & $11$ &$10$ &$12$ &$12$ &$12$ &$12$ &$12$ &$12$ &$11$ &$11$ &$12$ &$12$ &$12$ \tabularnewline
& 12 & $4$ &$10$ &$12$ &$12$ &$12$ &$1$ &$12$ &$12$ &$7$ &$12$ &$12$ &$12$ &$12$ \tabularnewline
\end{tabular}
\caption{Lexicographically largest $\ell \in \mathbb{Z}_{13}$ corresponding with the minimum values of $ \left| G_3 \cap \phi^{-1} \left( \ell \right)  \right|  + \left| Y_3 \cap \phi^{-1} \left( \ell \right)\right|$  for $m,n \geq 13$.}
\label{table:pathpositin}
\end{table}

\begin{table}[htbp]
\renewcommand\arraystretch{0.9}
\centering
\begin{tabular}{c|c|cccccccccccc}
\multicolumn{1}{c|}{} &   & \multicolumn{10}{c}{$m$} \\
\multicolumn{1}{c|}{$n_x$}&  & $2$ & $3$ & $4$ & $5$ & $6$ & $7$ & $8$ & $9$ & $10$ & $11$ & $12$\\ \hline
      $0$  & & $1$&$0$ &$1$ &$1$& $2$ &$2$ &$2$ &$2$  &$2$ &$3$ &$3$ \\
      $1$  & &    &    &$2$ &   & $2$ &$3$ &$4$ &$3$  &$4$ &$5$ &$4$\\
      $2$  & &    &    &$3$ &   & $4$ &$4$ &$5$ &$5$  &$6$ &$6$ &$7$\\
\multirow{2}{*}{$6$} &$n=21$ &    &    &$4$ &   & $5$ &$5$ &$7$ &$6$  &$7$ &$8$ &$9$\\
        &$n \neq 21$ &    &    &$4$ &   & $5$ &$5$ &$7$ &$6$  &$8$ &$8$ &$9$\\
      $4$  & &    &    &$4$ &   & $6$ &$7$ &$8$ &$8$  &$9$ &$10$&$11$\\
      $5$  & &    &    &$5$ &   & $7$ &$8$ &$9$ &$10$ &$11$&$12$&$13$\\
\multirow{2}{*}{$6$} & $n=6$&    &    &$6$ &   & $8$ &$9$ &$11$&$11$ &$13$&$14$&$15$\\
      & $n>6$&    &    &$6$ &   & $9$ &$9$ &$11$&$11$ &$13$&$14$&$15$\\
      $7$  & &    &    &$7$ &   & $9$ &$11$&$12$&$13$ &$15$&$16$&$17$\\
      $8$  & &    &    &$8$ &   & $11$&$12$&$14$&$15$ &$16$&$18$&$20$\\
      $9$  & &    &    &$8$ &   & $11$&$13$&$15$&$16$ &$18$&$19$&$21$\\
      $10$ & &    &    &    &   & $13$&$14$&$16$&     &$20$&$22$&$24$\\
      $11$ & &    &    &    &   & $14$&$16$&$18$&     &$22$&$24$&$25$\\
      $12$ & &    &    &    &   & $15$&$17$&$20$&     &$24$&$25$&$28$\\
      $13$ & &    &    &    &   & $16$&$18$&$21$&     &$25$&$28$&$29$\\
      $14$ & &    &    &    &   & $18$&    &$23$&     &$27$&$30$&$32$\\
      $15$ & &    &    &    &   & $18$&    &$24$&     &$29$&$31$&$34$\\
      $16$ & &    &    &    &   &     &    &$25$&     &$31$&$33$&$36$\\
      $17$ & &    &    &    &   &     &    &$27$&     &$32$&$35$&$38$\\
      $18$ & &    &    &    &   &     &    &$28$&     &    &$37$&$40$\\
      $19$ & &    &    &    &   &     &    &$30$&     &    &$39$&$42$\\
      $20$ & &    &    &    &   &     &    &$31$&     &    &$41$&$45$\\
      $21$ & &    &    &    &   &     &    &$32$&     &    &$43$&$46$\\
      $22$ & &    &    &    &   &     &    &    &     &    &$44$&$49$\\
      $23$ & &    &    &    &   &     &    &    &     &    &$47$&$50$\\
      $24$ & &    &    &    &   &     &    &    &     &    &$49$&\\
      $25$ & &    &    &    &   &     &    &    &     &    &$50$&\\%\cline{3-13}
%& & \multicolumn{11}{c}{$c(m,n,n_x)$} \\
\end{tabular}
\caption{Value of $c(m,n,n_x)$ in the upper bound of $ \gamma_{b,2} \left(P_m \square P_n\right)$ for $2 \leq m \leq 12 $ and $n \geq m$.}
\label{table:pathxpathconstanterm}
\end{table}
% default
\renewcommand\arraystretch{1}

\begin{table}[htbp]
\begin{center}
\begin{tabular}{ccc|ccccccccccccc}
&& \multicolumn{13}{c}{Least residue $n_{13}$ of $n$ modulo 13} \tabularnewline
&& & 0 & 1 & 2 & 3 & 4 & 5 & 6 & 7 & 8 & 9 & 10 & 11 & 12\\ \hline
\parbox[t]{2mm}{\multirow{13}{*}{\rotatebox[origin=c]{90}{Least residue $m_{13}'$}}} & \parbox[t]{2mm}{\multirow{13}{*}{\rotatebox[origin=c]{90}{of $m-10$ modulo 13}}}
&$0$ & $0$ & -$2$ & $13$ & -$21$ & $1$ & -$7$ & -$18$ & -$6$ & -$14$ & -$12$ & -$10$ & $2$ & -$9$\tabularnewline 
&&$1$ & $0$ & $6$ & $6$ & -$18$ & -$9$ & -$7$ & -$20$ & -$9$ & -$7$ & -$7$ & -$7$ & $4$ & $4$\tabularnewline 
&&$2$ & $0$ & $1$ & $12$ & -$15$ & -$6$ & -$7$ & -$22$ & $1$ & $0$ & -$15$ & -$17$ & $6$ & -$9$\tabularnewline 
&&$3$ & $0$ & $9$ & $5$ & -$12$ & -$3$ & -$7$ & -$24$ & -$2$ & -$6$ & -$10$ & -$14$ & -$5$ & -$9$\tabularnewline 
&&$4$ & $0$ & $4$ & -$2$ & -$22$ & -$13$ & -$20$ & -$26$ & -$5$ & -$12$ & -$18$ & -$11$ & -$3$ & $4$\tabularnewline 
&&$5$ & $0$ & -$1$ & $4$ & -$19$ & $3$ & -$7$ & -$15$ & -$8$ & -$5$ & -$13$ & -$21$ & -$1$ & -$9$\tabularnewline 
&&$6$ & $0$ & $7$ & $10$ & -$16$ & -$7$ & -$7$ & -$17$ & -$11$ & -$11$ & -$8$ & -$18$ & $1$ & -$9$\tabularnewline 
&&$7$ & $0$ & $2$ & $16$ & -$13$ & -$4$ & -$7$ & -$19$ & -$1$ & -$4$ & -$16$ & -$15$ & $3$ & -$9$\tabularnewline 
&&$8$ & $0$ & $10$ & $9$ & -$23$ & -$1$ & -$7$ & -$21$ & -$4$ & -$10$ & -$11$ & -$12$ & -$8$ & -$9$\tabularnewline 
&&$9$ & $0$ & $5$ & $2$ & -$20$ & -$11$ & -$7$ & -$23$ & -$7$ & -$16$ & -$6$ & -$9$ & $7$ & -$9$\tabularnewline 
&&$10$ & $0$ & $0$ & $8$ & -$17$ & -$8$ & -$7$ & -$25$ & -$10$ & -$9$ & -$14$ & -$19$ & -$4$ & -$9$\tabularnewline 
&&$11$ & $0$ & $8$ & $1$ & -$14$ & -$5$ & -$7$ & -$27$ & -$13$ & -$15$ & -$9$ & -$16$ & -$2$ & -$9$\tabularnewline 
&&$12$ & $0$ & $3$ & $7$ & -$24$ & -$2$ & -$7$ & -$29$ & -$3$ & -$8$ & -$17$ & -$13$ & $0$ & $4$\tabularnewline 
\end{tabular}
\end{center}
\caption{Constant term $c(m_{13}',n_{13})$ for the upper bound for $ \gamma_{b,2} \left( P_{m} \square C_{n} \right)$ for $m \geq 23$ and $n \geq 13$. }
\label{table:pathxcyclegeneralcornerforforumla}
\end{table}

\begin{table}[htbp]
\renewcommand\arraystretch{0.9}
\setlength{\tabcolsep}{4pt}
\centering
\begin{tabular}{c|c|ccccccccccccccccccccc}
\multicolumn{1}{c|}{} &  &     \multicolumn{20}{c}{$m$}\\
\multicolumn{1}{c|}{$n_x$}&  & $2$ & $3$ & $4$ & $5$ & $6$ & $7$ & $8$ & $9$ & $10$ & $11$ & $12$ & $13$ & $14$ & $15$ & $16$ & $17$ & $18$ & $19$ & $20$ & $21$ & $22$\\ \hline
      $0$  & & $0$&$0$ &$0$ &$0$& $0$ &$0$ &$0$ &$0$  &$0$   &$0$   &$0$  &$0$ &$0$ &$0$& $0$ &$0$ &$0$ &$0$  &$0$   &$0$   &$0$\\ %&\multirow{26}{*}{$c(m,n_x)$}\\
      $1$  & & & &$2$&$1$&$2$ &$3$ &$3$ &$3$  &$3$ &$4$ &$4$ &$6 $&$6 $&$7 $&$7 $&$7 $&$8 $&$8 $&$8 $&$9 $&$9 $\\ 
      $2$  & & & &$2$&   &$3$ &$3$ &$4$ &$4$  &$4$ &$5$ &$5$ &$8 $&$8 $&$9 $&$9 $&$9 $&$10$&$11$&$11$&$12$&$12$\\ 
\multirow{2}{*}{$3$} &$n=3$  & & &$3$&   &$4$ &$5$ &$6$ &$6$  &$7$ &$8$ &$8$ &$9 $&$10$&$10$&$11$&$11$&$12$&$13$&$14$&$14$&$15$\\ 
     &$n>3$  & & &$3$&   &$4$ &$5$ &$6$ &$6$  &$7$ &$8$ &$8$ &$9 $&$11$&$11$&$12$&$13$&$14$&$15$&$15$&$17$&$16$\\ 
\multirow{2}{*}{$4$} &$n=4$  & & &$4$&   &$5$ &$6$ &$6$ &$7$  &$8$ &$8 $&$9$ &$10$&$10$&$10$&$12$&$12$&$13$&$14$&$14$&$15$&$16$\\ 
     &$n>4$  & & &$4$&   &$5$ &$6$ &$7$ &$7$  &$8$ &$9 $&$9$ &$11$&$11$&$12$&$12$&$14$&$14$&$15$&$16$&$17$&$18$\\ 
\multirow{2}{*}{$5$} &$n=5$  & & &$5$&   &$7$ &$8$ &$9$ &$10$ &$11$&$12$&$13$&$14$&$15$&$16$&$17$&$18$&$19$&$20$&$21$&$22$&$23$\\
     &$n>5$  & & &$5$&   &$7$ &$8$ &$9$ &$10$ &$10$&$12$&$12$&$14$&$15$&$16$&$17$&$18$&$19$&$20$&$21$&$22$&$23$\\
\multirow{2}{*}{$6$} &$n=6$  & & &$5$&   &$7$ &$8$ &$9$ &$10$ &$11$&$12$&$13$&$14$&$15$&$16$&$17$&$18$&$19$&$20$&$21$&$22$&$23$\\
     &$n>6$  & & &$5$&   &$7$ &$8$ &$9$ &$10$ &$11$&$12$&$13$&$15$&$16$&$18$&$19$&$20$&$21$&$22$&$23$&$24$&$26$\\
\multirow{2}{*}{$7$} &$n=7$  & & &$6$&   &$9$ &$10$&$11$&$12$ &$14$&$15$&$16$&$17$&$18$&$20$&$21$&$22$&$23$&$24$&$26$&$27$&$28$\\
     &$n>7$  & & &$6$&   &$9$ &$10$&$11$&$13$ &$14$&$15$&$16$&$17$&$18$&$20$&$22$&$22$&$24$&$25$&$26$&$27$&$28$\\
\multirow{2}{*}{$8$} &$n=8$  & & &$7$&   &$10$&$11$&$12$&$14$ &$15$&$16$&$18$&$19$&$20$&$22$&$23$&$24$&$26$&$27$&$28$&$30$&$31$\\
     &$n>8$  & & &$7$&   &$10$&$11$&$13$&$14$ &$15$&$17$&$18$&$20$&$21$&$23$&$24$&$25$&$27$&$28$&$30$&$31$&$33$\\
\multirow{2}{*}{$9$} &$n=9$  & & &$8$&   &$11$&$13$&$14$&$16$ &$18$&$19$&$21$&$22$&$24$&$26$&$27$&$29$&$30$&$32$&$34$&$34$&$37$\\
     &$n>9$  & & &$8$&   &$11$&$13$&$14$&$16$ &$17$&$19$&$20$&$22$&$24$&$26$&$28$&$29$&$31$&$33$&$34$&$36$&$37$\\
\multirow{2}{*}{$10$}&$n=10$ & & &   &   &$12$&$13$&$15$&     &$18$&$20$&$21$&$23$&$24$&$26$&$28$&$29$&$31$&$32$&$34$&$38$&$37$\\
     &$n>10$ & & &   &   &$12$&$13$&$15$&     &$18$&$20$&$21$&$24$&$26$&$27$&$29$&$31$&$32$&$35$&$36$&$38$&$40$\\
      $11$ & & & &   &   &$14$&$15$&$18$&     &$21$&$23$&$25$&$27$&$28$&$30$&$32$&$33$&$36$&$37$&$39$&$41$&$32$\\
\multirow{2}{*}{$12$}&$n=12$ & & &   &   &$14$&$16$&$18$&     &$22$&$24$&$26$&$28$&$30$&$32$&$34$&$36$&$38$&$40$&$42$&$44$&$46$\\
     &$n>12$ & & &   &   &$14$&$16$&$18$&     &$22$&$24$&$26$&$29$&$30$&$32$&$34$&$36$&$38$&$40$&$42$&$44$&$46$\\
      $13$ & & & &   &   &$16$&$18$&$20$&     &$24$&$27$&$29$  \\
\multirow{2}{*}{$14$}&$n=14$ & & &   &   &$17$&    &$21$&     &$25$&$28$&$30$  \\
     &$n>14$ & & &   &   &$17$&    &$22$&     &$25$&$28$&$30$  \\
      $15$ & & & &   &   &$18$&    &$23$&     &$28$&$31$&$33$  \\
      $16$ & & & &   &   &    &    &$24$&     &$29$&$31$&$34$  \\
      $17$ & & & &   &   &    &    &$26$&     &$32$&$34$&$37$  \\
      $18$ & & & &   &   &    &    &$27$&     &    &$36$&$39$  \\
      $19$ & & & &   &   &    &    &$29$&     &    &$38$&$41$  \\
      $20$ & & & &   &   &    &    &$30$&     &    &$39$&$42$  \\
      $21$ & & & &   &   &    &    &$32$&     &    &$42$&$46$  \\
      $22$ & & & &   &   &    &    &    &     &    &$43$&$47$  \\
      $23$ & & & &   &   &    &    &    &     &    &$46$&$50$  \\
      $24$ & & & &   &   &    &    &    &     &    &$48$&\\
      $25$ & & & &   &   &    &    &    &     &    &$50$&\\ %\cline{4-24}
\end{tabular}
\caption{Value of $c(m,n,n_x)$ in the upper bound of $ \gamma_{b,2} \left(P_m \square C_n\right)$ for $2 \leq m \leq 22 $ and $n \geq 3$.}
\label{table:pathxcycleconstanterm}
\end{table}
% default
\setlength{\tabcolsep}{6pt}
% default
\renewcommand\arraystretch{1}

\renewcommand\arraystretch{0.9}
\begin{table}[htbp]
\centering
\begin{tabular}{c|cccccccccc}
&     \multicolumn{10}{c}{$n$}\\
$m_x$ & $3$ & $4$ & $5$ &$6$ & $7$ & $8$ & $9$ & $10$ & $11$ & $12$\\ \hline
$0$ &$0$& $1$ &$1$& $1$ &$2$ &$2$ &$2$  &$2$ &$3$ &$2$ \\%&\multirow{26}{*}{$c(m,n_x)$}\\
$1$ &   & $2$ &   & $ $ &$4$ &$3$ &$3$  &$4$ &$5$ &$2$ \\
$2$ &   & $2$ &   & $ $ &$5$ &$4$ &$5$  &$5$ &$6$ &$7$ \\
$3$ &   & $3$ &   & $ $ &$6$ &$6$ &$6$  &$7$ &$9$ &$9$ \\
$4$ &   & $4$ &   &     &$7$ &$7$ &$8$  &$8$ &$10$&$11$\\
$5$ &   & $4$ &   &     &$8$ &$8$ &$10$ &$10$&$12$&$13$\\
$6$ &   &     &   &     &$10$ &    &$11$ &$12$&$14$&$15$\\
$7$ &   &     &   &     &$11$&    &$13$ &$13$&$15$&$17$\\
$8$ &   &     &   &     &$12$&    &$15$ &$15$&$18$&$19$\\
$9$ &   &     &   &     &$13$&    &$16$ &$16$&$19$&$21$\\
$10$&   &     &   &     &$14$&    &     &    &$21$&$23$\\
$11$&   &     &   &     &$16$&    &     &    &$23$&$25$\\
$12$&   &     &   &     &$17$&    &     &    &$24$&$27$\\
$13$&   &     &   &     &$18$& &&&&   \\
$14$&   &     &   &     &$19$& &&&&   \\
$15$&   &     &   &     &$20$& &&&&   \\
$16$&   &     &   &     &$22$& &&&&   \\
$17$&   &     &   &     &$23$& &&&&   \\
$18$&   &     &   &     &$24$& &&&&   \\
$19$&   &     &   &     &$25$& &&&&   \\
$20$&   &     &   &     &$26$& &&&&   \\
$21$&   &     &   &     &$28$& &&&&   \\
$22$&   &     &   &     &$29$& &&&&   \\
$23$&   &     &   &     &$30$& &&&&   \\
$24$&   &     &   &     &$31$& &&&&   \\
$25$&   &     &   &     &$32$& &&&&   \\
$26$&   &     &   &     &$34$& &&&&   \\
$27$&   &     &   &     &$35$& &&&&   \\
$28$&   &     &   &     &$36$& &&&&   \\
$29$&   &     &   &     &$37$& &&&&   \\
$30$&   &     &   &     &$38$& &&&&   \\
$31$&   &     &   &     &$40$& &&&&   \\
$32$&   &     &   &     &$41$& &&&&   \\
$33$&   &     &   &     &$42$& &&&&   \\
$34$&   &     &   &     &$43$& &&&&   \\%\cline{2-11}
\end{tabular}
\caption{Value of $c(n,m_x)$ in the upper bound of $ \gamma_{b,2} \left(P_m \square C_n\right)$ for $m \geq 23$ and $3 \leq n \leq 12$.}
\label{table:pathxcyclefinalconstanterm}
\end{table}
% default
\renewcommand\arraystretch{1}

\renewcommand\arraystretch{0.9}
\setlength{\tabcolsep}{4pt}
\begin{table}[htbp]
\centering
\begin{tabular}{c|ccccccccccccccccccccccc}
&     \multicolumn{20}{c}{$m$}\\
$n_x$ & $3$ & $4$ & $5$ &$6$ & $7$ & $8$ & $9$ & $10$ & $11$ & $12$& $13$ & $14$ & $15$& $16$ & $17$& $18$ & $19$ & $20$& $21$ &$22$& $23$& $24$& $25$\\ \hline
$0$ &$0$& $0$ &$0$& $0$ &$0$ &$0$ &$0$  &$0$ &$0$ &$0$ &$0$ &$0$ &$0$ &$0$ &$0$ &$0$ &$0$ &$0$ &$0$ &$0$ &$0$ &$0$ &$0$ \\%&\multirow{26}{*}{$c(m,n_x)$}\\
$1$ &   & $2$ &   & $1$ &$3$ &$2$ &$2$  &$4$ &$5$ &$4$ &$5$ &$5$ &$6$ &$6$ &$7$ &$7$ &$7$ &$8$ &$8$ &$9$ &$9$ &$9$ &$9$ \\
$2$ &   & $2$ &   & $1$ &$4$ &$4$ &$4$  &$5$ &$6$ &$6$ &$7$ &$8$ &$8$ &$8$ &$9$ &$9$ &$10$&$11$&$11$&$12$&$12$&$13$&$12$\\
$3$ &   & $3$ &   & $1$ &$5$ &$6$ &$5$  &$7$ &$8$ &$8$ &$10$&$10$&$11$&$12$&$13$&$13$&$14$&$15$&$15$&$16$&$17$&$18$&$18$\\
$4$ &   & $3$ &   &     &$6$ &$6$ &$7$  &$7$ &$9 $&$10$&$10$&$11$&$12$&$12$&$13$&$13$&$14$&$16$&$15$&$18$&$17$&$19$&$19$\\
$5$ &   & $4$ &   &     &$7$ &$8$ &$9 $ &$10$&$11$&$11$&$13$&$14$&$15$&$16$&$17$&$18$&$19$&$20$&$21$&$22$&$23$&$24$&$25$\\
$6$ &   &     &   &     &$8$ &    &$10$ &$11$&$13$&$14$&$15$&$16$&$17$&$18$&$19$&$21$&$22$&$23$&$23$&$25$&$26$&$27$&$28$\\
$7$ &   &     &   &     &$9 $&    &$12$ &$13$&$14$&$16$&$16$&$18$&$19$&$20$&$21$&$23$&$24$&$25$&$25$&$27$&$29$&$30$&$32$\\
$8$ &   &     &   &     &$10$&    &$13$ &$14$&$16$&$16$&$19$&$21$&$22$&$22$&$25$&$26$&$28$&$29$&$30$&$32$&$33$&$34$&$34$\\
$9$ &   &     &   &     &$12$&    &$15$ &$16$&$18$&$20$&$21$&$23$&$25$&$26$&$28$&$29$&$31$&$33$&$34$&$36$&$37$&$39$&$40$\\
$10$&   &     &   &     &$13$&    &     &    &$20$&$22$&$23$&$25$&$27$&$28$&$30$&$32$&$34$&$35$&$36$&$39$&$39$&$43$&$45$\\
$11$&   &     &   &     &$14$&    &     &    &$21$&$24$&$24$&$27$&$29$&$30$&$32$&$35$&$36$&$37$&$39$&$41$&$43$&$45$&$47$\\
$12$&   &     &   &     &$16$&    &     &    &$24$&$26$&$26$&$30$&$32$&$33$&$35$&$36$&$40$&$42$&$42$&$46$&$48$&$50$&$51$\\
$13$&   &     &   &     &$16$&    \\
$14$&   &     &   &     &$18$&    \\
$15$&   &     &   &     &$19$&    \\
$16$&   &     &   &     &$20$&    \\
$17$&   &     &   &     &$21$&    \\
$18$&   &     &   &     &$22$&    \\
$19$&   &     &   &     &$24$&    \\
$20$&   &     &   &     &$25$&    \\
$21$&   &     &   &     &$26$&    \\
$22$&   &     &   &     &$27$&    \\
$23$&   &     &   &     &$29$&    \\
$24$&   &     &   &     &$30$&    \\
$25$&   &     &   &     &$31$&    \\
$26$&   &     &   &     &$32$&    \\
$27$&   &     &   &     &$33$&    \\
$28$&   &     &   &     &$34$&    \\
$29$&   &     &   &     &$36$&    \\
$30$&   &     &   &     &$37$&    \\
$31$&   &     &   &     &$38$&    \\
$32$&   &     &   &     &$39$&    \\
$33$&   &     &   &     &$40$&    \\
$34$&   &     &   &     &$42$&    \\
\end{tabular}
\caption{Value of $c(m,n_x)$ in the upper bound of $ \gamma_{b,2} \left(C_m \square C_n\right)$ for $3 \leq m \leq 25$ and $n \geq m$.}
\label{table:cyclexcycleconstanterm}
\end{table}
% default
\setlength{\tabcolsep}{6pt}
% default
\renewcommand\arraystretch{1}

\begin{table}[htbp]
\centering
\begin{tabular}{cc|ccccccccccccc} 
& \multicolumn{13}{c}{Least residue $n_{13}$ of $n$ modulo 13} \tabularnewline 
& & 0 & 1 & 2 & 3 & 4 & 5 & 6 & 7 & 8 & 9 & 10 & 11 & 12\\ \hline 
\parbox[t]{2mm}{\multirow{13}{*}{\rotatebox[origin=c]{90}{Least residue $m_{13}$ of $m$ modulo $13$}}} 
&$0 $& $0$ & $0$ & $0$ & $0$ & $0$ & $0$ & $0$ & $0$ & $0$ & $0$ & $0$ & $0$ & $0$\tabularnewline 
&$1 $& $0$ & $21$ & $26$ & $58$ & $22$ & $15$ & $33$ & $24$ & $30$ & $35$ & $40$ & $18$ & $23$\tabularnewline 
&$2 $& $0$ & $26$ & $20$ & $55$ & $19$ & $15$ & $35$ & $14$ & $23$ & $43$ & $37$ & $16$ & $10$\tabularnewline 
&$3 $& $0$ & $71$ & $55$ & $81$ & $33$ & $46$ & $82$ & $37$ & $63$ & $47$ & $70$ & $51$ & $35$\tabularnewline 
&$4 $& $0$ & $22$ & $19$ & $46$ & $35$ & $36$ & $59$ & $13$ & $27$ & $24$ & $47$ & $14$ & $11$\tabularnewline 
&$5 $& $0$ & $28$ & $28$ & $59$ & $23$ & $54$ & $28$ & $10$ & $28$ & $54$ & $41$ & $10$ & $23$\tabularnewline 
&$6 $& $0$ & $33$ & $48$ & $82$ & $46$ & $28$ & $30$ & $26$ & $47$ & $49$ & $64$ & $34$ & $36$\tabularnewline 
&$7 $& $0$ & $11$ & $14$ & $37$ & $13$ & $10$ & $26$ & $9$ & $32$ & $22$ & $25$ & $21$ & $-2$\tabularnewline 
&$8 $& $0$ & $17$ & $23$ & $50$ & $27$ & $28$ & $34$ & $32$ & $33$ & $26$ & $45$ & $30$ & $23$\tabularnewline 
&$9 $& $0$ & $35$ & $43$ & $60$ & $37$ & $41$ & $49$ & $22$ & $52$ & $34$ & $68$ & $28$ & $36$\tabularnewline 
&$10 $& $0$ & $27$ & $37$ & $70$ & $47$ & $41$ & $51$ & $25$ & $58$ & $42$ & $65$ & $26$ & $10$\tabularnewline 
&$11 $& $0$ & $18$ & $16$ & $51$ & $14$ & $10$ & $34$ & $21$ & $43$ & $28$ & $26$ & $26$ & $11$\tabularnewline 
&$12 $& $0$ & $23$ & $23$ & $35$ & $24$ & $36$ & $23$ & $-2$ & $75$ & $23$ & $23$ & $11$ & $-2$\tabularnewline 
\end{tabular} 
\caption{Constant term $c(m_{13},n_{13})$ in the upper bound for $ \gamma_{b,2} \left( C_{m} \square C_{n} \right)$ for $m, n\geq 26$. }
\label{table:cyclxcyclegeneralcornerfinal}
\end{table}

\begin{table}[htbp]
\renewcommand\arraystretch{1.2}
\[
\arraycolsep=3.5pt
\begin{array}{cccccccccc}
\multicolumn{9}{c}{m}\\
$14$ & $15$ & $16$ & $17$ & $18$ & $19$ & $20$ & $21$ & $22$ \\
\hline
\belowbaseline[-8pt]{$\begin{bmatrix}\frac{599}{1982}\\  \frac{185}{1982}\\  \frac{297}{1982}\\  \frac{351}{1982}\\  \frac{136}{991}\\  \frac{317}{1982}\\  \frac{303}{1982}\\  \frac{303}{1982}\\  \frac{317}{1982}\\  \frac{136}{991}\\  \frac{351}{1982}\\  \frac{297}{1982}\\  \frac{185}{1982}\\  \frac{599}{1982}\\  \end{bmatrix}$}
&
\belowbaseline[-8pt]{$\begin{bmatrix}\frac{688}{2277}\\  \frac{71}{759}\\  \frac{340}{2277}\\  \frac{45}{253}\\  \frac{104}{759}\\  \frac{40}{253}\\  \frac{359}{2277}\\  \frac{112}{759}\\  \frac{359}{2277}\\  \frac{40}{253}\\  \frac{104}{759}\\  \frac{45}{253}\\  \frac{340}{2277}\\  \frac{71}{759}\\  \frac{688}{2277}\\  \end{bmatrix}$}
&
\belowbaseline[-8pt]{$\begin{bmatrix}\frac{888}{2939}\\  \frac{275}{2939}\\  \frac{1757}{11756}\\  \frac{2085}{11756}\\  \frac{405}{2939}\\  \frac{464}{2939}\\  \frac{458}{2939}\\  \frac{1791}{11756}\\  \frac{1791}{11756}\\  \frac{458}{2939}\\  \frac{464}{2939}\\  \frac{405}{2939}\\  \frac{2085}{11756}\\  \frac{1757}{11756}\\  \frac{275}{2939}\\  \frac{888}{2939}\\  \end{bmatrix}$}
&
\belowbaseline[-8pt]{$\begin{bmatrix}\frac{6119}{20250}\\  \frac{631}{6750}\\  \frac{1514}{10125}\\  \frac{599}{3375}\\  \frac{278}{2025}\\  \frac{119}{750}\\  \frac{3151}{20250}\\  \frac{508}{3375}\\  \frac{1591}{10125}\\  \frac{508}{3375}\\  \frac{3151}{20250}\\  \frac{119}{750}\\  \frac{278}{2025}\\  \frac{599}{3375}\\  \frac{1514}{10125}\\  \frac{631}{6750}\\  \frac{6119}{20250}\\  \end{bmatrix}$}
& 
\belowbaseline[-8pt]{$\begin{bmatrix}\frac{10537}{34873}\\  \frac{3262}{34873}\\  \frac{5211}{34873}\\  \frac{144}{811}\\  \frac{4792}{34873}\\  \frac{5515}{34873}\\  \frac{5454}{34873}\\  \frac{5241}{34873}\\  \frac{5416}{34873}\\  \frac{5416}{34873}\\  \frac{5241}{34873}\\  \frac{5454}{34873}\\  \frac{5515}{34873}\\  \frac{4792}{34873}\\  \frac{144}{811}\\  \frac{5211}{34873}\\  \frac{3262}{34873}\\  \frac{10537}{34873}\\  \end{bmatrix}$}
&
\belowbaseline[-8pt]{$\begin{bmatrix}\frac{42}{139}\\  \frac{13}{139}\\  \frac{187}{1251}\\  \frac{74}{417}\\  \frac{172}{1251}\\  \frac{22}{139}\\  \frac{65}{417}\\  \frac{21}{139}\\  \frac{194}{1251}\\  \frac{64}{417}\\  \frac{194}{1251}\\  \frac{21}{139}\\  \frac{65}{417}\\  \frac{22}{139}\\  \frac{172}{1251}\\  \frac{74}{417}\\  \frac{187}{1251}\\  \frac{13}{139}\\  \frac{42}{139}\\  \end{bmatrix}$}
&
\belowbaseline[-8pt]{$\begin{bmatrix}\frac{31248}{103415}\\  \frac{9671}{103415}\\  \frac{15458}{103415}\\  \frac{18357}{103415}\\  \frac{384}{2795}\\  \frac{16376}{103415}\\  \frac{1241}{7955}\\  \frac{3114}{20683}\\  \frac{16119}{103415}\\  \frac{15848}{103415}\\  \frac{15848}{103415}\\  \frac{16119}{103415}\\  \frac{3114}{20683}\\  \frac{1241}{7955}\\  \frac{16376}{103415}\\  \frac{384}{2795}\\  \frac{18357}{103415}\\  \frac{15458}{103415}\\  \frac{9671}{103415}\\  \frac{31248}{103415}\\  \end{bmatrix}$}
&
\belowbaseline[-8pt]{$\begin{bmatrix}\frac{26905}{89043}\\  \frac{2776}{29681}\\  \frac{26617}{178086}\\  \frac{10537}{59362}\\  \frac{12238}{89043}\\  \frac{4697}{29681}\\  \frac{13898}{89043}\\  \frac{8945}{59362}\\  \frac{27665}{178086}\\  \frac{4572}{29681}\\  \frac{13625}{89043}\\  \frac{4572}{29681}\\  \frac{27665}{178086}\\  \frac{8945}{59362}\\  \frac{13898}{89043}\\  \frac{4697}{29681}\\  \frac{12238}{89043}\\  \frac{10537}{59362}\\  \frac{26617}{178086}\\  \frac{2776}{29681}\\  \frac{26905}{89043}\\  \end{bmatrix}$}
&
\belowbaseline[-8pt]{$\begin{bmatrix}\frac{2155}{7132}\\  \frac{667}{7132}\\  \frac{11460}{76669}\\  \frac{13608}{76669}\\  \frac{10537}{76669}\\  \frac{1129}{7132}\\  \frac{47835}{306676}\\  \frac{11559}{76669}\\  \frac{11920}{76669}\\  \frac{11770}{76669}\\  \frac{47169}{306676}\\  \frac{47169}{306676}\\  \frac{11770}{76669}\\  \frac{11920}{76669}\\  \frac{11559}{76669}\\  \frac{47835}{306676}\\  \frac{1129}{7132}\\  \frac{10537}{76669}\\  \frac{13608}{76669}\\  \frac{11460}{76669}\\  \frac{667}{7132}\\  \frac{2155}{7132}\\  \end{bmatrix}$}
\end{array}
\]
\caption{Vectors used in the proof of the lower bound for $mp_{f,2}(P_m \square C_n)$ for $14 \leq m \leq 22$}
\label{table:multipackvectoer1422}
\end{table}
% default
\renewcommand\arraystretch{1}

\end{document}